\documentclass[12pt,a4paper,asmart]{memoir}

\usepackage[OT4]{polski}
\usepackage[utf8]{inputenc} 

\usepackage{a4wide}
\usepackage{amssymb}
\usepackage{amsmath}
\usepackage{amsthm}
\usepackage{wrapfig,epsfig}
\usepackage{color}
\usepackage{t1enc}
\usepackage{txfonts}
\usepackage{listings}
\usepackage{indentfirst}

\newtheorem{definicja}{Definicja}[chapter]
\newtheorem{wniosek}{Wniosek}[chapter]
\newtheorem{uwaga}{Uwaga}[chapter]
\newtheorem{lemat}{Lemat}[chapter]
\newtheorem{stwierdzenie}{Stwierdzenie}[chapter]

\newtheorem{twierdzenie}{Twierdzenie}[chapter]

\makechapterstyle{mychapterstyle}{%

}

\chapterstyle{mychapterstyle}

\setsecheadstyle{\large\sffamily\bfseries}
\setsubsecheadstyle{\large\sffamily\bfseries}
\setsubsubsecheadstyle{\normalfont\sffamily\bfseries}
\setparaheadstyle{\normalfont\sffamily}

\makeevenhead{headings}{\thepage}{}{\small\slshape\leftmark}
\makeoddhead{headings}{\small\slshape\rightmark}{}{\thepage}

\settocdepth{subsection}
\setsecnumdepth{subsection}
\maxsecnumdepth{subsection}
\settocdepth{subsection}
\maxtocdepth{subsection}

\begin{document}

\pagestyle{plain}
\pagenumbering{arabic}
\sffamily

\thispagestyle{empty}
\begin{center}
    \Large
    Uniwersytet Gdański\\
    Wydział Matematyki, Fizyki i Informatyki\\
    Instytut Matematyki
\end{center}
\hrule

\vfill\vfill\vfill\vfill

\begin{center}
    \Huge\bfseries
    Ergodyczne własności pewnych stochastycznych układów dynamicznych
\end{center}

\vfill
\begin{center}
    \Large
    Hanna Wojewódka
\end{center}

\vfill\vfill\vfill
\begin{center}
    \Large
    Rozprawa doktorska napisana pod kierunkiem\\
		prof. dr. hab. Tomasza Szarka
\end{center}

\vfill\vfill

\hrule
\begin{center}
\large
    Gdańsk, 2015
\end{center}

\newpage

\vspace*{\fill}
{
\sffamily\itshape\small Serdeczne podziękowania kieruję do mojego Promotora - Pana prof. dr. hab. Tomasza Szarka. Dziękuję za zasugerowanie tematu niniejszej rozprawy oraz nieocenioną pomoc w~jej przygotowaniu, a~także za zarażenie pasją do pracy naukowej. 

Za owocną współpracę dziękuję również Pani dr hab. Katarzynie Horbacz oraz Panu prof.~Sandrowi Hille.}
\newpage

\rmfamily
\normalfont


\tableofcontents*
\thispagestyle{empty}

\newpage
\chapter*{Wstęp}

Rozprawa dotyczy badania ergodycznych własności pewnych stochastycznych układów dynamicznych generowanych przez łańcuchy Markowa o~wartościach w~przestrzeni stanów będącej przestrzenią polską.

Rozwój teorii operatorów Markowa zapoczątkował w 1907 roku A.A. Markow \cite{markow}. Rosyjski matematyk zdefiniował łańcuch jako nieskończony ciąg zmiennych losowych $x_0,x_1,\ldots,x_k,\ldots$ o~takiej własności, że zmienna $x_{k+1}$ jest niezależna od $x_0,x_1,\ldots,x_{k-1}$, gdy $x_k$ jest znana, dla dowolnej liczby naturalnej $k$. Łańcuchem jednorodnym A.A. Markow nazwał każdy łańcuch, dla którego rozkład zmiennej $x_{k+1}$ pod warunkiem $x_k$ jest niezależny od $k$. W~1952 roku istotny krok w~rozwoju teorii zrobił W.~Feller, wprowadzając operatory Markowa działające na~miarach. Operatory te mogą opisywać ewolucję rozkładów dla stochastycznych układów dynamicznych, w szczególności ciągłych iterowanych układów funkcyjnych, które często służą do modelowania procesów w~biologii, fizyce czy naukach społecznych. 

W rozprawie rozważane będą operatory Markowa działające na miarach borelowskich określonych na przestrzeniach polskich. Analizie zostanie poddany pewien stochastyczny układ dynamiczny opisujący proces podziału komórki. Jeden z pierwszych modeli cyklu komórkowego zaproponowali w~1988 roku J.J. Tyson i K.B. Hannsgen \cite{tysonhannsgen}, a pełne biologiczne tło badań nad procesem podziału komórki zostało zestawione w pracy A. Murraya i T. Hunta \cite{murrayhunt}. W~1999~roku ukazał się ważny wynik A. Lasoty i M.C. Mackeya \cite{lasotam}. Autorzy rozważyli model cyklu komórkowego wyrażony prostym stochastycznym układem dynamicznym. Udowodnili, że jest on asymptotycznie stabilny, co okazało się istotne nie tylko z punktu widzenia matematyki, ale i~biologii. S. Hille  i współautorzy \cite{hille} zaproponowali uogólnienie modelu A. Lasoty i~M.C.~Mackeya. 

Ergodyczny opis uogólnionego modelu cyklu komórkowego jest celem tej rozprawy. Chcąc badać ergodyczne własności, w pierwszej kolejności należy pytać o istnienie miary niezmienniczej. Następnie  można dowodzić asymptotycznej stabilności. Kryterium istnienia miar niezmienniczych dla nierozszerzających operatorów Markowa na przestrzeniach polskich wyznaczył T.~Szarek \cite{szarek}. Badania nad asymptotyczną stabilnością operatorów, najpierw na lokalnie zwartych oraz $\sigma$-zwartych przestrzeniach metrycznych, a następnie na przestrzeniach polskich, prowadzili m.in. A. Lasota i J.A. Yorke \cite{lasotay} oraz T.~Szarek \cite{szarek_stability}. W niniejszej pracy pokażemy asymptotyczną stabilność oraz oszacujemy tempo zbieżności do miary niezmienniczej dla pewnych, niekoniecznie nierozszerzających, operatorów Markowa określonych na miarach na~przestrzeniach polskich. W tym celu, wzorując się na wynikach M.~Hairera \cite{hairer}, zastosujemy technikę opartą na odpowiedniej konstrukcji miary sprzęgającej. Wykorzystana będzie również zupełność rodziny miar probabilistycznych z metryką Fortet-Mouriera. Geometryczne tempo zbieżności pozwoli dowieść kolejnych ergodycznych własności, takich jak centralne twierdzenie graniczne (CTG) czy prawo iterowanego logarytmu (PIL). 

Na przełomie XX i XXI wieku ukazało się wiele prac dotyczących CTG dla stacjonarnych i~ergodycznych procesów Markowa, m.in. \cite{gordin}, \cite{kv}, \cite{dl} czy \cite{woodr}. Później zaczęto dowodzić CTG dla asymptotycznie stabilnych, ale niestacjonarnych łańcuchów Markowa (zob. np.  \cite{komorowskiwalczuk}). 

Funkcjonalna forma PIL została zdefiniowana przez V. Strassena w 1964 roku. Jego wynik, omówiony w pracy \cite{strassen}, powszechnie nazywa się zasadą niezmienniczości Strassena. Ważny rezultat dla martyngałów opublikowali w 1973 roku C.C. Heyde i D.J. Scott \cite{heyde_scott}. 
Aby dowieść zasady niezmienniczości PIL dla szerszej klasy procesów stochastycznych, zaczęto korzystać z~pewnych metod martyngałowych (zob. \cite{bms}, \cite{kom_szar}).

Niniejsza praca składa się z siedmiu rozdziałów. W pierwszym rozdziale omówimy podstawowe definicje i~oznaczenia, a także stwierdzenia i uwagi wykorzystywane w~dalszej części rozprawy.

Rozdział drugi poświęcony będzie opisowi zarówno prostego, jak i uogólnionego modelu cyklu komórkowego. Ponadto zestawione zostaną założenia przyjęte w rozprawie.

W rozdziale trzecim będziemy badać własności uogólnionego modelu cyklu komórkowego. Wprowadzimy~również model pomocniczy generowany przez niejednorodne łańcuchy Markowa i odpowiadające im iterowane układy funkcyjne, dla których skonstruujemy miarę sprzęgającą na całych trajektoriach.

Miara sprzęgająca będzie pełnić kluczową rolę w dowodach głównych twierdzeń tej rozprawy. Jej konstrukcja, wzorowana na pracy M. Hairera \cite{hairer}, zostanie opisana w rozdziale czwartym. W~podobny sposób miarę sprzęgającą dla klasycznych iterowanych układów funkcyjnych zbudował m.in. M. Ślęczka \cite{sleczka}.

W rodziale piątym pokażemy istnienie jedynej miary niezmienniczej, 
a także dowiedziemy asymptotycznej stabilności modelu.  Ponadto oszacujemy tempo zbieżności kolejnych iteracji badanego operatora Markowa do miary niezmienniczej w normie Fortet-Mouriera (twierdzenie~\ref{GTZ}). Twierdzenie \ref{GTZ} stanie się punktem wyjścia do badania CTG i PIL.

W rozdziale szóstym sformułowane i wykazane zostanie CTG dla uogólnionego modelu cyklu komórkowego (twierdzenie \ref{CTG}). 
Zaproponowany dowód będzie bazował na gotowych wynikach M. Maxwella i M. Woodroofa \cite{woodr} oraz odpowiednich własnościach miary sprzęgającej, dzięki czemu będzie mniej techniczny i znacznie krótszy niż standardowe dowody, które wymagają wyprowadzenia pełnej metody martyngałowej. 

Przedmiotem ostatniego rozdziału będzie dowód PIL (twierdzenie~\ref{PIL}). Zaczniemy od przytoczenia wyników C.C. Heyde'go i D.J. Scotta \cite{heyde_scott} dla martyngałów, z~których korzystali m.in. autorzy prac \cite{bms} czy \cite{kom_szar}. Zastosowanie własności miary sprzęgającej pozwoli uprościć znaną dotychczas adaptację metody martyngałowej. 

Fragmenty badań zaprezentowanych w rozprawie zostały spisane w pracach \cite{hw}, \cite{hhsw} i~ \cite{hsw}.

\chapter{Podstawowe definicje i oznaczenia}

Zacznijmy od uzgodnienia podstawowych oznaczeń.

Niech $N$ oznacza zbiór liczb naturalnych $\left\{1,2,\ldots\right\}$ oraz niech $N_0=N\cup\left\{0\right\}$. Zbiór liczb rzeczywistych oznaczmy przez $R$.

Niech przestrzeń metryczna $(X,\varrho)$ będzie przestrzenią polską, tzn. ośrodkową i zupełną. Symbolem $B_X$ oznaczmy $\sigma$-ciało borelowskich podzbiorów przestrzeni $X$. W przestrzeni $(X,\varrho)$ określmy kulę otwartą $B(x,r)=\left\{y\in X: \varrho(x,y)<r\right\}$ o~środku w~punkcie $x\in X$ i~promieniu $r>0$ oraz kulę domkniętą $\bar{B}(x,r)=\left\{y\in X:\varrho(x,y)\leq r\right\}$.

Ustalmy $\bar{x}\in X$. 
Jeśli przestrzeń $(X,\varrho)$ jest nieograniczona, to dowolną funkcję ciągłą $V:X\to[0,\infty)$ spełniającą warunek
\begin{align*}
\lim_{\varrho(x,\bar{x})\to\infty}V(x)=\infty
\end{align*} 
nazywamy funkcją Lapunowa.

Rodzinę $(\Omega,\mathcal{F},\mathbb{P})$ składającą się z niepustego zbioru $\Omega$, określonego na nim $\sigma$-ciała $\mathcal{F}$ oraz miary probabilistycznej $\mathbb{P}:\mathcal{F}\to[0,1]$ nazywamy przestrzenią probabilistyczną. 
Gdy dane zdarzenie $A$ zachodzi z prawdopodobieństwem $\mathbb{P}(A)=1$, piszemy krótko, że zachodzi $\mathbb{P}$-p.n.

Dla dowolnego zbioru $A\in B_X$, symbolem $1_A$ określmy indykator zbioru $A$, to jest funkcję
\begin{align*}
1_A(x)=\left\{\begin{array}{ll}
1&\quad \textrm{dla }x\in A,\\
0&\quad \textrm{dla }x\in X\backslash A,
\end{array}\right.
\end{align*}
z kolei przez $\delta_x$ oznaczmy się miarę Diraca w punkcie $x\in X$, czyli 
\begin{align*}
\delta_x(A)=\left\{\begin{array}{ll}
1&\quad \textrm{dla }x\in A,\\
0&\quad \textrm{dla }x\in X\backslash A.
\end{array}\right.
\end{align*}

\section{Przestrzenie funkcyjne i miary na przestrzeniach polskich}

Niech $B(X)$ oznacza przestrzeń wszystkich funkcji $f:X\to R$ mierzalnych i ograniczonych, z~normą supremum $\|f\|_{\infty}=\sup_{x\in X}|f(x)|$. Podprzestrzeń  
$B(X)$ składającą się z funkcji ciągłych określmy symbolem $C(X)$. W pracy będziemy również wyróżniać przestrzeń $\tilde{B}(X)$ funkcji $f:X\to R$ mierzalnych i ograniczonych z dołu.

Przez $M(X)$ oznaczmy rodzinę wszystkich miar borelowskich na $B_X$. Rozważmy jej następujące podrodziny:
\begin{align*}
&M_{{fin}}(X)=\left\{\mu\in M(X):\;\mu(X)<\infty\right\},\\
&M_{1}(X)=\left\{\mu\in M(X):\;\mu(X)=1\right\}.
\end{align*}
Elementy $M_{{fin}}(X)$ oraz $M_1(X)$ nazywa się odpowiednio miarami skończonymi i probabilistycznymi. Mówimy, że miara $\mu$ jest sub-probabilistyczna, jeśli $\mu(X)\leq 1$. 
Nośnik miary $\mu\in M_{{fin}}(X)$ zdefiniujmy jako zbiór
\begin{align*}
\textrm{supp }\mu=\left\{x\in X:\;\mu(B(x,r))>0\;\textrm{ dla }r>0\right\}.
\end{align*}

Ustalmy punkt $\bar{x}\in X$. Niech $M_1^1(X)$ będzie rodziną miar probabilistycznych z~pierwszym momentem skończonym, a więc zbiorem postaci
\begin{align*}
M_1^1(X)=\left\{\mu\in M_1(X):\:\int_X\varrho(x,\bar{x})\mu(dx)<\infty\right\}.
\end{align*}
Zauważmy, że definicja rodziny $M_1^1(X)$ nie zależy od wyboru $\bar{x}\in X$. W analogiczny sposób możemy definiować przestrzenie miar z~$r$-tym momentem skończonym
\begin{align*}
M_1^r(X)=\left\{\mu\in M_1(X):\:\int_X\varrho^r(x,\bar{x})\mu(dx)<\infty\right\}\quad \text{dla }r>0.
\end{align*}
Symbolem $M_{{sig}}(X)$ określmy przestrzeń miar znakozmiennych, czyli zbiór postaci
\begin{align*}
M_{{sig}}(X)=\left\{\mu_1-\mu_2:\:\mu_1,\mu_2\in M_{{fin}}(X)\right\}.
\end{align*}
W przestrzeni $M_{{sig}}(X)$ rozważmy normę $\|\cdot\|$ zadaną wzorem
\begin{align*}
\|\mu\|=\mu^+(X)+\mu^-(X)\quad \text{dla }\mu\in M_{sig}(X),
\end{align*}
gdzie $\mu^+$ oraz $\mu^-$ są odpowiednio wahaniem: górnym oraz dolnym. Istnienie miar $\mu^+$, $\mu^-$ wynika z twierdzeń Hahna i Jordana o rozkładzie (zob. \S 29, \cite{halmos}). W szczególności, jeśli miara $\mu$ jest nieujemna, to $\|\mu\|=\mu(X)$.

Dla wygody będziemy korzystać z notacji iloczynu skalarnego
\begin{align*}
\left\langle f,\mu \right\rangle =\int_Xf(x)\mu(dx)\quad\textrm{dla $f\in \tilde{B}(X)$ oraz $\mu\in M_{sig}(X)$.}
\end{align*}

W przestrzeni $M_{{sig}}(X)$ możemy określić normę Fourtet-Mouriera
\begin{align}
\begin{aligned}
\|\mu\|_{\mathcal{L}}=\sup_{f\in\mathcal{L}}|\left\langle f,\mu\right\rangle |,
\end{aligned}
\end{align}
gdzie 
\begin{align}\label{def:L}
\begin{aligned}
\mathcal{L}=\left\{f\in C(X):\;|f(x)-f(y)|\leq\varrho(x,y),\;|f(x)|\leq 1\;\textrm{dla }x,y\in X\right\}
\end{aligned}
\end{align}
(por. \cite{lasota}, \cite{rachev} lub \cite{villani}).

\begin{stwierdzenie}\label{stwierdzenie:zupelność}
Jeśli $(X,\varrho)$ jest przestrzenią polską, to przestrzeń $(M_1(X),\|\cdot\|_{\mathcal{L}})$ z metryką indukowaną przez normę Fortet-Mouriera jest zupełna (twierdzenie~1.48,~\cite{lasota}).
\end{stwierdzenie}

\section{Zbieżność ciągów miar}
Mówimy, że ciąg miar borelowskich $(\mu_n)_{n\in N_0}\subset M_{fin}(X)$ jest słabo zbieżny do $\mu\in M_{fin}(X)$ (por. np. \cite{billingsley} lub \cite{kurtz}) i piszemy $\mu_n\xrightarrow{w}\mu$, gdy $n\to\infty$, jeśli
\begin{align*}
\lim_{n\to\infty}\left\langle f,\mu_n\right\rangle=\left\langle f,\mu\right\rangle\quad\textrm{dla wszystkich }f\in C(X).
\end{align*}
Granica słabo zbieżnego ciągu miar jest określona w sposób jednoznaczny (wniosek~1.39,~\cite{lasota}).

\begin{stwierdzenie}\label{tw_aleksandrowa}
Niech dana będzie miara $\mu\in M_{fin}(X)$ oraz ciąg miar $(\mu_n)_{n\in N_0}\subset M_{fin}(X)$. Następujące warunki są równoważne (twierdzenie 11.3.3, \cite{dudley}):
\begin{itemize}
\item[(1)] ciąg $(\mu_n)_{n\in N_0}$ zbiega słabo do miary $\mu$;
\item[(2)] $\lim_{n\to\infty}\left\langle f,\mu_n\right\rangle=\left\langle f,\mu\right\rangle\;$ dla każdej funkcji $f\in \mathcal{L}$;
\item[(3)] $\lim_{n\to\infty}\|\mu_n-\mu\|_{\mathcal{L}}=0$.
\end{itemize}
\end{stwierdzenie}

\section{Operatory Markowa i Fellera}

Mówimy, że operator $P:M_{{fin}}(X)\to M_{{fin}}(X)$ jest operatorem Markowa, jeśli spełnione są następujące warunki:
\begin{itemize}
\item $P(\lambda_1\mu_1+\lambda_2\mu_2)=\lambda_1P(\mu_1)+\lambda_2P(\mu_2)\quad\textrm{dla } \lambda_1,\lambda_2\geq 0,\:\mu_1,\mu_2\in M_{{fin}}(X)$;
\item $P\mu(X)=\mu(X)\quad\textrm{dla } \mu\in M_{{fin}}(X)$.
\end{itemize}
Operator Markowa $P$, dla którego istnieje liniowy operator $U:B(X)\to B(X)$ o~własności
\begin{align}\label{dualny}
\left\langle Uf,\mu \right\rangle =\left\langle f,P\mu \right\rangle \quad\textrm{dla } f\in B(X),\: \mu\in M_{{fin}}(X),
\end{align}
nazywamy operatorem regularnym. 
O operatorze $U$ spełniającym warunek (\ref{dualny}) mówimy, że jest dualny do operatora $P$.
\begin{uwaga}\label{uwaga:rozszerzenie_U}
Operator $U$ można rozszerzyć do przestrzeni $\tilde{B}(X)$ (wniosek 3.7, \cite{lasota}).
\end{uwaga}

Niech dana będzie funkcja $\Pi:X\times B_X\to R$ o własnościach:
\begin{itemize}
\item[(a)] przy ustalonym $x\in X$, funkcja $\Pi(x,\cdot):B_X\to R$ jest miarą probabilistyczną,
\item[(b)] przy ustalonym $A\in B_X$, funkcja $\Pi(\cdot,A):X\to R$ jest borelowsko mierzalna.
\end{itemize}
Funkcję $\Pi$ spełniającą powyższe warunki nazywamy funkcją przejścia.
\begin{stwierdzenie}\label{stwierdzenie:operatory_regularne_z_funkcji_przejścia}
Każda funkcja $\Pi:X\times B_X\to R$ spełniająca warunki (a) i (b) generuje jedyny regularny operator Markowa $P$, dla którego jest funkcją przejścia. Operator $P$ jest postaci
\begin{align*}
P\mu(A)=\int_X\Pi(x,A)\mu(dx)\quad \text{dla $A\in B_X$, $\mu\in M_{{fin}}(X)$.}
\end{align*}
Ponadto operator dualny do $P$ wyraża się wzorem
\begin{align*}
Uf(x)=\int_Xf(y)\Pi(x,dy)\quad\text{dla $f\in B(X)$, $x\in X$.} 
\end{align*}
Spełniona jest więc następująca równość
\begin{align*}
\Pi(x,A)=U1_A(x)=P\delta_x(A)
\end{align*}
dla dowolnych $x\in X$ oraz $A\in B_X$ (por. rozdz. 1.1, \cite{z}).
\end{stwierdzenie}

Każdy operator Markowa $P$ zdefiniowany na $M_{{fin}}(X)$ możemy jednoznacznie rozszerzyć na~przestrzeń miar znakozmiennych $M_{{sig}}(X)$. Wystarczy przyjąć, że 
\begin{align*}
P\mu=P\mu_1-P\mu_2\quad\text{dla }\mu=\mu_1-\mu_2,\:\mu_1,\mu_2\in M_{{fin}}(X).
\end{align*}

Mówimy, że regularny operator Markowa jest operatorem Fellera lub fellerowskim, jeśli
\begin{align*}
U(C(X))\subset C(X).
\end{align*}
Odwzorowanie $\kappa:X\to M_{fin}(X)$ nazywamy słabo ciągłym, jeśli $\kappa(x_n)\xrightarrow{w}\kappa(x)$ dla wszystkich $x\in X$ oraz $(x_n)_{n\in N_0}\subset X$ takich, że $\lim_{n\to\infty}x_n=x$.
\begin{stwierdzenie}\label{stwierdzenie:kiedy_Markowa_jest_Fellera}
Jeśli $P$ jest operatorem regularnym, indukowanym przez funkcję przejścia $\Pi$, to $P$ jest fellerowski wtedy i tylko wtedy, gdy odwzorowanie
\begin{align*}
X\ni x\mapsto \Pi(x,\cdot)\in M_1(X)
\end{align*}
jest słabo ciągłe (zob. rozdz. 6, \cite{tweedie}).
\end{stwierdzenie}

\section{Proces Markowa i własność Markowa}

Niech $(X,B_X,\mu)$ będzie przestrzenią z nieujemną miarą $\mu\in M_{fin}(X)$. O~zbiorach postaci $A_1\times A_2$ dla $A_1,A_2\in B_X$ mówimy, że są prostokątami. Oznaczmy przez $B_X\otimes B_X$ $\sigma$-ciało generowane przez wszystkie prostokąty. Zdefiniujmy miarę $\mu\times \mu$ na $B_X\otimes B_X$ jako jedyną miarę o~własności 
\[(\mu\times\mu)(A_1\times A_2)=\mu(A_1)\mu(A_2)\quad\text{dla }A_1,A_2\in B_X.\]
Istnienie takiej miary wynika np. z twierdzenia 3.3.1, \cite{bogachev}.

Zauważmy, że indukcyjnie możemy rozszerzyć definicję miar produktowych na dowolną skończenie wymiarową przestrzeń $X^n=X\times\ldots\times X$ (por. \S 37, \cite{halmos}). Niech $\otimes_{i=1}^n B_X$ będzie $\sigma$-ciałem generowanym przez wszystkie prostokąty postaci $A_1\times\ldots\times A_n$ dla $A_1,\ldots,A_n\in B_X$. Miara produktowa $\mu\times \ldots\times \mu$ określona na $\otimes_{i=1}^nB_X$ jest jedyną miarą taką, że 
\begin{align*}
\left(\mu\times\ldots\times\mu\right)\left(A_1\times\ldots\times A_n\right)=\Pi_{i=1}^n\mu\left(A_i\right) \quad\text{dla }A_1,\ldots,A_n\in B_X.
\end{align*}

Weźmy przestrzeń nieskończenie wymiarową $X^{\infty}=X\times X\times\ldots$. Niech $(X,B_X,\mu)$ będzie przestrzenią mierzalną oraz $\mu(X)=1$. Zdefiniujmy zbiór $\otimes_{i=1}^{\infty}B_X$ jako $\sigma$-ciało generowane przez cylindry, czyli zbiory postaci $A_1\times\ldots\times A_n\times X\times X\times\ldots$ dla $A_1\times\ldots\times A_n\in \otimes_{i=1}^n B_X$ oraz $n\in N$. Istnieje jedyna miara produktowa $\bar{\mu}$ określona na $\otimes_{i=1}^{\infty}B_X$ o własności
\[\bar{\mu}(E)=\left(\mu\times\ldots\times\mu\right)\left(A_1\times\ldots\times A_n\right)=\Pi_{i=1}^n\mu\left(A_i\right)\quad\text{dla }E=A_1\times\ldots\times A_n\times X\times X\times \ldots\]
(twierdzenie B z \S 38, \cite{halmos}). Pełna konstrukcja miar produktowych na przestrzeniach nieskończenie wymiarowych została opisana w rozdziale VII, \cite{halmos}, a także w rozdziale III, \cite{bogachev}.

\begin{stwierdzenie}\label{stwierdzenie:istnienie_procesu_Markowa}
Jeśli $\mu\in M_{fin}(X)$ jest rozkładem zmiennej losowej $x_0$ oraz $\Pi:X\times B_X\to R$ jest funkcją przejścia, to istnieje proces stochastyczny $(x_n)_{n\in N_0}$ określony na $X^{\infty}$ i mierzalny względem $\sigma$-ciała $\otimes_{i=1}^{\infty}B_X$ oraz taka miara probabilistyczna  $\mathbb{P}_{\mu}$ na $\otimes_{i=1}^{\infty}B_X$, że $\mathbb{P}_{\mu}(A)$ wyraża prawdopodobieństwo zajścia zdarzenia $\left\{(x_0,x_1,\ldots)\in A\right\}$ dla $A\in \otimes_{i=1}^{\infty}B_X$. Ponadto
\begin{align*}
\begin{aligned}
\mathbb{P}_{\mu}(x_0\in A_0,\ldots,x_n\in A_n)=\int_{A_0}\int_{A_1}\ldots\int_{A_{n-1}}\Pi(s_{n-1},A_n)\Pi(s_{n-2},ds_{n-1})\ldots \Pi(s_0,ds_1)\mu(ds_0)
\end{aligned}
\end{align*} 
dla $n\in N_0$ oraz wszystkich $A_0,\ldots, A_n\in B_X$ (twierdzenie 3.4.1, \cite{tweedie}).
\end{stwierdzenie}
Symbolem $E_{\mu}$ określmy wartość oczekiwaną względem miary $\mathbb{P}_{\mu}$. 
Jeśli $\mu=\delta_x$, to użyjemy oznaczeń $\mathbb{P}_x$ oraz ${E}_x$.

Niech $(X^{\infty}, \otimes_{i=1}^{\infty}B_X,\mathbb{P}_{\mu})$ będzie przestrzenią probabilistyczną zadaną w tezie stwierdzenia \ref{stwierdzenie:istnienie_procesu_Markowa}. 
Zdefiniujmy $\sigma$-ciała $\mathcal{F}_n=\sigma(x_0,\ldots, x_n)$ dla $n\in N_0$. Ciąg $(\mathcal{F}_n)_{n\in N_0}$ nazywamy filtracją naturalną procesu $(x_n)_{n\in N_0}$. 

Niech $\upsilon:X^{\infty}\to R$ będzie funkcją ograniczoną i mierzalną. 
Rozważmy zmienną losową 
$\upsilon(x_0,x_1,\ldots)$, 
gdzie $(x_n)_{n\in N_0}$ jest pewnym procesem startującym z miary $\mu\in M_{fin}(X)$, czyli takim procesem, że rozkład zmiennej $x_0$ jest zadany miarą $\mu$. 
Wprowadźmy operator przesunięcia $\Gamma$ taki, że
\begin{align*}
\Gamma\left((h_0,\ldots,h_n,\ldots)\right)=(h_1,\ldots,h_{n+1},\ldots)\quad\text{dla dowolnego ciągu stałych }(h_n)_{n\in N_0}\subset X.
\end{align*}
Niech $\Gamma_k$ oznacza $k$-tą iterację odwzorowania $\Gamma$, a więc 
\begin{align*}
\Gamma_k((h_0,\ldots,h_n,\ldots))=(h_k,\ldots,h_{n+k},\ldots)\quad\text{dla }k\in N.
\end{align*}
Zauważmy, że
\begin{align*}
\begin{aligned}
\left(\upsilon\circ\Gamma_n\right)(x_0,x_1,\ldots)=\upsilon(x_n,x_{n+1},\ldots)\quad\textrm{ dla }n\in N.
\end{aligned}
\end{align*}
Mówimy, że $(x_n)_{n\in N_0}$ jest procesem Markowa (por. \cite{durrett}), gdy
\begin{align*}
E_{\mu}\left(\upsilon \circ\Gamma_n|x_0,x_1,\ldots, x_n ;x_n=x\right)=E_x(\upsilon) \quad\text{dla }x\in X.
\end{align*}
Własnością Markowa nazywamy następujący warunek
\begin{align}\label{wl:E_markow}
\begin{aligned}
{E}_{\mu}\left(\upsilon\circ\Gamma_n|\mathcal{F}_n\right)(\omega)
={E}_{x_n(\omega)}(\upsilon)\quad \text{dla }\omega\in \Omega
\end{aligned}
\end{align}
(por. rozdz. 3, \cite{tweedie}).

Załóżmy, że $D:\Omega\to N_0$ jest momentem stopu względem filtracji $(\mathcal{F}_n)_{n\in N_0}$, tzn. zmienną losową taką, że $\left\{D=n\right\}\in \mathcal{F}_n$ dla $n\in N_0$. Ponadto $\sigma$-ciało $\mathcal{F}_D$ jest dane wzorem 
\begin{align*}
\mathcal{F}_D=\left\{F\in \otimes_{i=1}^{\infty}B_X:\;\left\{D=n\right\}\cap F\in \mathcal{F}_n\right\}.
\end{align*}
Możemy rozważać zmienną losową $x_D$, która na zbiorze $\left\{D=n\right\}$ jest postaci $x_D=x_n$. 
Jeśli ciąg $(x_n)_{n\in N_0}$ spełnia na zbiorze $\left\{D<\infty\right\}$ warunek
\begin{align}\label{mocna_wl:E_markow}
\begin{aligned}
{E}_{\mu}\left(\upsilon\circ \Gamma_D|\mathcal{F}_D\right)(\omega)
={E}_{x_D(\omega)}(\upsilon)\quad\text{dla }\omega\in \Omega,
\end{aligned}
\end{align}
to mówimy, że spełnia on mocną własność Markowa (por. \cite{tweedie}).
\begin{stwierdzenie}
Każdy łańcuch Markowa z czasem dyskretnym spełnia mocną własność Markowa (stwierdzenie 3.4.6, \cite{tweedie}).
\end{stwierdzenie}

\section{Pojęcie asymptotycznej stabilności}
Operator Markowa $P$ nazywamy mieszającym, jeśli
\begin{align}\label{def:mieszjący}
\lim_{n\to\infty}\|P^n\mu_1-P^n\mu_2\|_{\mathcal{L}}=0\quad\textrm{dla }\mu_1,\mu_2\in M_{1}(X).
\end{align}
O mierze $\mu_*\in M(X)$ mówimy, że jest niezmiennicza względem operatora $P$, jeśli $P\mu_*=\mu_*$. 
Operator Markowa $P$ nazywamy asymptotycznie stabilnym, jeśli istnieje miara niezmiennicza $\mu_*\in M_{1}(X)$ taka, że
\begin{align}\label{def:asymptotyczna_stabiloność}
\lim_{n\to\infty}\|P^n\mu-\mu_*\|_{\mathcal{L}}=0\quad\textrm{dla dowolnej }\mu\in M_{1}(X).
\end{align}

\begin{stwierdzenie}\label{stwierdzenie:jedyność_miary}
Jeśli operator Markowa jest mieszający, tzn. spełnia warunek (\ref{def:mieszjący}), to istnieje co najwyżej jedna unormowana miara niezmiennicza względem tego operatora (uwaga 4.21, \cite{lasota}).
\end{stwierdzenie}

\begin{stwierdzenie}\label{stwierdzenie:Feller_daje_niezmienniczość}
Warunek (\ref{def:asymptotyczna_stabiloność}) nie musi implikować niezmienniczości miary $\mu_*$. Implikacja jest jednak prawdziwa, gdy operator $P$ jest fellerowski (twierdzenie 3.49, \cite{lasota}). 
\end{stwierdzenie}

\chapter{Opis modelu matematycznego}

W rozdziale drugim zarysujemy biologiczne tło badanego modelu, jak również opiszemy założenia przyjęte w pracy. 

\section{Inspiracja - prosty model cyklu komórkowego}
Poniżej przedstawimy krótko model cyklu komórkowego opisany prostym stochastycznym układem dynamicznym (por. \cite{lasotam}).

Ustalmy $T\in (0,\infty)$. 
Załóżmy, że każda komórka w rozważanej populacji składa się z~$j$~różnych substancji, których masy są opisane wektorem 
$y(t)=\left(y^1(t),\ldots,y^j(t)\right)$, gdzie $t\in[0,T]$ oznacza wiek komórki. Niech wektor $y(t)$ będzie dany wzorem $y(t)=\pi(x,t)$, gdzie $\pi:X\times[0,T]\to X$ jest pewną ustaloną funkcją oraz $\pi(x,0)=x$. Przykładem może być wektor $y(t)$ spełniający układ równań różniczkowych zwyczajnych 
$dy/dt=g(y(t))$ 
z warunkiem początkowym $y(0)=x$ i rozwiązaniem postaci $y(t)=\pi(x,t)$.

Niech $x_n$ oznacza wartość początkową substancji $x=y(0)$ w $n$-tym pokoleniu, a $t_n$ --~moment podziału komórki w $n$-tym pokoleniu. Ponadto niech
\begin{align}\label{rozklad_podzilow}
\mathbb{P}\left(t_n<t|x_n=x\right)=\int_0^tp(x,s)ds\quad \text{ dla }t\in[0,T], \;n\in N_0,
\end{align}
gdzie $\mathbb{P}$ jest miarą określoną na~pewnej przestrzeni mierzalnej $(\Omega, \mathcal{F})$. 
Wektor $y(t_n)=\pi(x_n,t_n)$ 
opisuje ilość wewnątrzkomórkowych substancji tuż przed mitozą w $n$-tym pokoleniu. Przyjmujemy, że komórka potomna składa się dokładnie z~połowy składników komórki macierzystej, a~zatem
\begin{align}\label{kom_potomna}
x_{n+1}=\frac{1}{2}\pi(x_n,t_n)\quad \text{ dla }n\in N_0.
\end{align}

Przejdźmy do matematycznego opisu biologicznych intuicji. 
Rozważmy przestrzeń polską $(X,\varrho)$ oraz przestrzeń probabilistyczną $(\Omega,\mathcal{F},\mathbb{P})$. Ustalmy $T\in(0,\infty)$. Niech $(t_n)_{n\in N_0}$ będzie ciągiem niezależnych zmiennych losowych o wartościach w $[0,T]$ oraz niech rozkład $t_{n+1}$ pod warunkiem $x_n=x$ będzie dany wzorem (\ref{rozklad_podzilow}). 
Załóżmy, że funkcja $p:X\times[0,T]\to[0,\infty)$ jest półciągła z~dołu, nieujemna i~unormowana, tzn. $\int_0^Tp(x,u)du=1$ dla $x\in X$. Ponadto przyjmijmy, że funkcja $S:X\times[0,T]\to X$ jest ciągła. Niech 
\begin{align*}
x_{n+1}=S(x_n,t_n)\quad\textrm{dla } n\in N_0.
\end{align*}
Wówczas operator Markowa $P:M_{fin}(X)\to M_{fin}(X)$ jest postaci 
\[P\delta_x(A)=\int_0^T 1_A(S(x,t))p(x,t)dt\quad\text{dla }x\in X, \;A\in B_X.\]

Powyższe założenia są spełnione m.in. przez modele z prac J.J.~Tysona i~K.B.~Hannsgena~\cite{tysonhannsgen} oraz A.~Murraya i~T.~Hunta \cite{murrayhunt}.

\section{Uogólniony model cyklu komórkowego. Założenia}
Niech $H$ będzie ośrodkową przestrzenią Banacha. O domkniętym podzbiorze $X$ przestrzeni $H$ możemy myśleć jak o przestrzeni polskiej $(X,\varrho)$, gdzie metryka $\varrho$ jest indukowana przez normę w przestrzeni $H$. 
Rozważmy dodatkowo przestrzeń probabilistyczną $(\Omega,\mathcal{F},\mathbb{P})$.  

Ustalmy $\varepsilon_*<\infty$. Weźmy $\varepsilon\in[0,\varepsilon_*]$ oraz $T\in (0,\infty)$. Rozważmy stochastycznie zaburzony układ dynamiczny postaci
\begin{align*}
x_{n+1}=S(x_n,t_{n+1})+H_{n+1}\quad \text{dla } n\in N_0,
\end{align*}
gdzie $(H_n)_{n\in N}$ jest ciągiem 
 niezależnych wektorów losowych o wartościach w $H$, z których wszystkie mają ten sam rozkład $\nu^{\varepsilon}$ spełniający warunek $\text{supp}\,\nu^{\varepsilon}\subset\bar{B}(0,\varepsilon)$. Ponadto ustalmy, że wektor $x_0$ nie zależy od $(H_n)_{n\in N}$.
 
Przyjmijmy następujące założenia definiujące model.
\begin{itemize}
\item[(I)] Niech $(t_n)_{n\in N_0}$ będzie ciągiem niezależnych zmiennych losowych na~$(\Omega,\mathcal{F},\mathbb{P})$ o~wartościach w~$[0,T]$. Rozkład $t_{n+1}$ pod warunkiem $x_n=x$ określmy wzorem
\begin{align*}
\mathbb{P}\left(t_{n+1}<t|x_n=x\right)=\int_0^tp(x,u)du\quad \text{dla }t\in[0,T],
\end{align*}
gdzie $p:X\times[0,T]\to[0,\infty)$ jest funkcją mierzalną, nieujemną i~unormowaną, tzn. $\int_0^T p(x,u)du=1$ dla wszystkich $x\in X$.
\item[(II)] Niech $S:X\times[0,T]\to X$ będzie funkcją ciągłą, spełniającą poniższy lokalny warunek Lipschitza
\begin{align}\label{zal:Lipschitz}
\varrho(S(x,t),S(y,t))\leq \lambda(x,t)\varrho(x,y)\quad\textrm{dla }x,y\in X,\: t\in[0,T],
\end{align}
gdzie $\lambda:X\times[0,T]\to[0,\infty)$ jest funkcją borelowską o własności
\begin{align}\label{def:Lambda1}
\Lambda_1:=\sup_{x\in X}\int_0^T\lambda(x,t)p(x,t)dt<1.
\end{align}
\item[(III)] Niech $\sup_{t\in[0,T]}\varrho(S(\bar{x},t),\bar{x})<\infty$ dla pewnego punktu $\bar{x}\in X$ oraz niech
\begin{align}\label{def:c}
c=\sup_{t\in[0,T]}\varrho(S(\bar{x},t),\bar{x})+\varepsilon_*.
\end{align}
\item[(IV)] Niech $p$ spełnia warunek Diniego, tj.
\begin{align*}
\int_0^T|p(x,t)-p(y,t)|dt\leq\omega(\varrho(x,y))\quad\textrm{dla }x,y\in X,
\end{align*}
gdzie $\omega:[0,\infty)\to[0,\infty)$ jest funkcją niemalejącą, wklęsłą i~taką, że $\omega(0)=0$ oraz
\begin{align}\label{wl:Dini_omega}
\int_0^{\sigma}\frac{\omega(t)}{t}dt<\infty\quad\textrm{dla pewnego }\sigma>0.
\end{align}
\item[(V)] Przyjmijmy, że funkcja $p$ jest ograniczona i zdefiniujmy stałe
\begin{align}\label{def:M_1}
M_1=\inf_{x\in X,\: t\in[0,T]}p(x,t),\\
\label{def:M_2}
M_2=\sup_{x\in X,\: t\in[0,T]}p(x,t).
\end{align}
Będziemy zakładać, że $M_1>0$.
\item[(VI)] Niech $\nu^{\varepsilon}$ będzie miarą borelowską na przestrzeni $H$ o~własności $\textrm{supp }\nu^{\varepsilon}\subset \bar{B}(0,\varepsilon)$. Połóżmy
$
\nu_x^{\varepsilon}(\cdot)=\nu^{\varepsilon}(\cdot-x)$ dla  $x\in X$.
Załóżmy, że 
\begin{align*}
S(x,t)+h\in X\quad\text{dla wszystkich}\;x\in X,\:t\in[0,T] \text{ oraz }h\in\textrm{supp }\nu^{\varepsilon}.
\end{align*}
\end{itemize}

\begin{uwaga}
Zauważmy, że dla $\varepsilon_*=0$ uogólniony model cyklu komórkowego sprowadza się do prostszej wersji modelu z podrozdziału 2.1.
\end{uwaga}

\begin{stwierdzenie}\label{stwierdzenie:Dini}
Niech $\omega:[0,\infty)\to[0,\infty)$ będzie funkcją niemalejącą, wklęsłą i taką, że $\omega(0)=0$. Ponadto niech $\omega$ spełnia warunek (\ref{wl:Dini_omega}). Dla wszystkich $t\in [0,T]$ oraz dowolnej stałej $\zeta<1$ zdefiniujmy 
\begin{align}\label{def:varphi}
\varphi(t)=\sum_{k=1}^{\infty}\omega(\zeta^kt).
\end{align}
Warunek Diniego (założenie (IV)) implikuje własności: 
\begin{align}\label{Dini:varphi<infty}
&\varphi(t)<\infty\quad\textrm{dla }t\in[0,T],\\
\label{Dini:lim_varphi=0}
\text{oraz }\qquad &\lim_{t\to 0}\varphi(t)=0.
\end{align}
\end{stwierdzenie}

\begin{proof}
Chcąc uzasadnić (\ref{Dini:varphi<infty}), wystarczy skorzystać z odpowiednich własności funkcji $\omega$ i~zauważyć, że
\begin{align*}
\infty>\int_0^{\sigma}\frac{\omega(t)}{t}dt=\sum_{k=0}^{\infty}\int_{\zeta^{k+1}\sigma}^{\zeta^k\sigma}\frac{\omega(t)}{t}dt\geq\sum_{k=0}^{\infty}\frac{\omega(\zeta^{k+1}\sigma)}{\zeta^k\sigma}(\zeta^k\sigma-\zeta^{k+1}\sigma)
\geq(1-\zeta)\varphi(\sigma).
\end{align*}

Ustalmy $\epsilon>0$ i $t_0\in[0,T]$. Dzięki własności (\ref{Dini:varphi<infty}) możemy tak wybrać $n_0\in N$, że
$
\sum_{k=n_0}^{\infty}\omega(\zeta^kt_0)<{\epsilon}/{2}$. 
Zauważmy, że funkcja $\omega$, która jest z założenia niemalejąca i nieujemna, a ponadto spełnia warunek (\ref{wl:Dini_omega}), jest również ciągła w~zerze. 
Istnieje więc $\delta\in(0,t_0)$ taka, że  
$\sum_{k=1}^{n_0}\omega(\zeta^kt)<{\epsilon}/{2}$  dla $t\in[0,\delta)$. 
Wówczas $\sum_{k=1}^{\infty}\omega(\zeta^kt)\leq \epsilon/2+\epsilon/2$ dla $t\in[0,\delta)$, co implikuje (\ref{Dini:lim_varphi=0}).
\end{proof}

Określmy funkcję $\Pi_{\varepsilon}:X\times B_X\to[0,1]$ wzorem
\begin{align}\label{def:Pi_varepsilon}
\Pi_{\varepsilon}(x,A)=\int_0^Tp(x,t)\nu^{\varepsilon}_{S(x,t)}(A) \,dt\quad\textrm{dla }x\in X\textrm{ oraz }A\in B_X.
\end{align}
Zauważamy, że odwzorowanie $\Pi_{\varepsilon}(x,\cdot):B_X\to R$ jest miarą probabilistyczną dla ustalonego $x\in X$ oraz funkcja $\Pi_{\varepsilon}(\cdot,A):X\to R$ jest mierzalna dla ustalonego $A\in B_X$. Ze~stwierdzenia \ref{stwierdzenie:operatory_regularne_z_funkcji_przejścia} wynika więc, że operator $P_{\varepsilon}:M_{fin}(X)\to M_{fin}(X)$ postaci
\begin{align}\label{def:P_varepsilon}
P_{\varepsilon}\mu(A)=\int_X\Pi_{\varepsilon}(x,A)\mu(dx)
\end{align}
jest regularnym operatorem Markowa, a operator do niego dualny $U_{\varepsilon}:B(X)\to B(X)$ wyraża się wzorem
\begin{align}\label{def:U_varepsilon}
U_{\varepsilon}f(x)=\int_Xf(y)\Pi_{\varepsilon}(x,dy).
\end{align}
Jak zauważyliśmy w uwadze \ref{uwaga:rozszerzenie_U}, operator dualny $U_{\varepsilon}$ można rozszerzyć do przestrzeni $\tilde{B}(X)$.

Aby udowodnić CTG i PIL, będziemy odpowiednio przyjmować dodatkowe założenia:
\begin{itemize}
\item[(II$^{\prime}$)] funkcja $S:X\times[0,T]\to X$ jest ciągła i spełnia lokalny warunek Lipschitza~(\ref{zal:Lipschitz}), gdzie $\lambda:X\times[0,T]\to[0,\infty)$ jest funkcją borelowską o własności
\begin{align}\label{def:Lambda2}
\Lambda_2:=\sup_{x\in X}\int_0^T\lambda^2(x,t)p(x,t)dt<1;
\end{align}

\item[(II$^{\prime\prime}$)] funkcja $S:X\times[0,T]\to X$ jest ciągła i spełnia lokalny warunek Lipschitza~(\ref{zal:Lipschitz}), 
gdzie $\lambda:X\times[0,T]\to[0,\infty)$ jest funkcją borelowską o własności
\begin{align}\label{def:Lambda2+delta}
\Lambda_{2+\delta}:=\sup_{x\in X}\int_0^T\lambda^{2+\delta}(x,t)p(x,t)dt<1.
\end{align}
\end{itemize}

Korzystając z nierówności H\"{o}ldera, łatwo pokazać, że warunek (II$^{\prime\prime}$) implikuje (II$^{\prime}$) oraz (II). Implikacje w drugą stronę nie są prawdziwe.

\chapter{Własności modelu matematycznego}\label{rozdzial:wlasnosci}

Wprowadzony zostanie model pomocniczy. Pokażemy pewne własności operatorów Markowa generujących model pomocniczy, a następnie uzasadnimy, dlaczego operator $P_{\varepsilon}$ również spełnia te własności.

\section{Model pomocniczy}

Niech funkcje  $T_{h}:X\times[0,T]\to X$, $h\in\bar{B}(0,\varepsilon)$, będą postaci
\begin{align*}
T_{h}(x,t)=S(x,t)+h\quad\text{dla } x\in X,\: t\in[0,T].
\end{align*}
Ustalmy ciąg stałych $(h_n)_{n\in N}\subset H$ takich, że $h_n\in\bar{B}(0,\varepsilon)$ dla $n\in N$ oraz ciąg niezależnych zmiennych losowych $(t_n)_{n\in N_0}$ o wartościach w $[0,T]$. Rozważmy stochastycznie zaburzony układ dynamiczny 
\begin{align*}
\tilde{x}_{n+1}=T_{h_{n+1}}(\tilde{x}_n,t_{n+1})=S(\tilde{x}_n,t_{n+1})+h_{n+1}\quad\textrm{dla } n\in N_0.
\end{align*}
Dla uproszczenia notacji będziemy pisać
\begin{align}\label{def:tilde_x}
\tilde{x}_{n+1}^{x_0}=T_{h_{n+1}}\left(T_{h_n}\left(\ldots T_{h_2}\left(T_{h_1}(x_0,t_1),t_2\right)\ldots\right),t_{n+1}\right)\quad\text{dla }x_0\in X.
\end{align}
Przyjmijmy, że funkcje $S$ i $p$ spełniają założenia (I)-(VI). Niech $h\in \bar{B}(0,\varepsilon)$. Zauważmy, że odwzorowanie $T_h:X\times [0,T]\to X$ jest ciągłe i spełnia lokalny warunek Lipschitza sformułowany w poprzednim rozdziale (zob. (\ref{zal:Lipschitz})). 
Zdefiniujmy funkcję $\Pi_h:X\times B_X\to [0,1]$ wzorem
\begin{align}\label{def:Pi_hn}
\Pi_{h}(x,A)=\int_0^T 1_A\left(T_{h}(x,t)\right)p(x,t)dt\quad \text{dla }x\in X,\:A\in B_X.
\end{align}
Zauważamy, że $\Pi_{h}(x,\cdot)$ jest miarą probabilistyczną dla ustalonego $x\in X$ oraz $\Pi_{h}(\cdot,A)$ jest funkcją mierzalną dla $A\in B_X$. Na mocy stwierdzenia \ref{stwierdzenie:operatory_regularne_z_funkcji_przejścia} operator $P_h$ postaci
\begin{align}
(P_{h}\mu)(A)=\int_X \Pi_{h}(x,A)\mu(dx)\quad \textrm{dla }A\in B_X,\: \mu\in M_{fin}(X)
\end{align}
jest regularnym operatorem Markowa, a operator do niego dualny $U_h$ możemy opisać wzorem
\begin{align}
(U_{h}f)(x)=\int_Xf(y)\Pi_{h}(x,dy)\quad\textrm{dla }x\in X,\: f\in B(X).
\end{align}
Zgodnie z uwagą \ref{uwaga:rozszerzenie_U} operator dualny $U_h$ można rozszerzyć do przestrzeni $\tilde{B}(X)$.

Ustalmy $\varepsilon\in[0,\varepsilon_*]$ dla $\varepsilon_*<\infty$. Przypomnijmy, że funkcja $\Pi_{\varepsilon}$ i operator $P_{\varepsilon}$ są odpowiednio dane wzorami (\ref{def:Pi_varepsilon}) i (\ref{def:P_varepsilon}). Zauważmy, że 
\begin{align}\label{Pi_epsilon-as-Pi_h} 
\Pi_{\varepsilon}(x,A)=\int_{\bar{B}(0,\varepsilon)}\Pi_h(x,A) \nu^{\varepsilon}(dh)\quad \text{dla }x\in X,\: A\in B_X
\end{align}
i wobec tego
\begin{align}\label{P_epsilon-as-P_h}
(P_{\varepsilon}\mu)(A)=\int_X\int_{\bar{B}(0,\varepsilon)}\Pi_h(x,A)\nu^{\varepsilon}(dh)\mu(dx)\quad \text{dla }A\in B_X,\: \mu\in M_{fin}(X).
\end{align}
Ponadto
\begin{align*}
\begin{aligned}
(U_{\varepsilon}f)(x)=\int_{\bar{B}(0,\varepsilon)}\int_Xf(y)\Pi_{h}(x,dy)\nu^{\varepsilon}(dh)\quad \text{dla }x\in X,\: f\in \tilde{B}(X).
\end{aligned}
\end{align*}

\section{Miary na trajektoriach}\label{sekcja:miary_na_trajektoriach}
Ustalmy $x\in X$, $h\in \bar{B}(0,\varepsilon)$ oraz taki ciąg stałych $(h_n)_{n\in N}\subset H$, że $h_n\in \bar{B}(0,\varepsilon)$ dla $n\in N$. Niech $A\in B_X$. Zdefiniujmy indukcyjnie następujące rozkłady jednowymiarowe
\begin{align}\label{def:1-wym_rozklady}
\begin{aligned}
&\Pi^0(x,A)=\delta_x(A)\\
&\Pi_{h}^1(x,A)=\Pi_{h}(x,A)\\
&\vdots\\
&\Pi_{h_1,\ldots,h_n}^n(x,A)=\int_X\Pi_{h_n}^1(y,A)\Pi_{h_1,\ldots,h_{n-1}}^{n-1}(x,dy).
\end{aligned}
\end{align}
Niech $n\in N_0$. Zauważmy, że $\Pi^n_{h_1,\ldots,h_n}(x,\cdot):B_X\to R$ jest miarą probabilistyczną dla dowolnego $x\in X$. Ponadto $\Pi^n_{h_1,\ldots, h_n}(\cdot,A):X\to R$ jest funkcją mierzalną dla $A\in B_X$. Na~mocy stwierdzenia~\ref{stwierdzenie:operatory_regularne_z_funkcji_przejścia} możemy rozważać operator Markowa $P^n_{h_1,\ldots,h_n}:M_{fin}(X)\to M_{fin}(X)$ postaci
\begin{align}
(P^n_{h_1,\ldots,h_n}\mu)(A)=\int_X\Pi_{h_1,\ldots,h_n}^n(x,A)\mu(dx)
\end{align}
oraz operator do niego dualny $U^n_{h_1,\ldots,h_n}:\tilde{B}(X)\to \tilde{B}(X)$ dany wzorem
\begin{align}
(U^n_{h_1,\ldots,h_n}f)(x)=\int_Xf(y)\Pi^n_{h_1,\ldots,h_n}(x,dy).
\end{align}

Określmy teraz rozkłady wielowymiarowe. Jeśli założymy, że $\Pi_{h_1,\ldots,h_n}^{1,\ldots,n}(x,\cdot)$ jest miarą na $X^n$ generowaną przez ciąg $\left(\Pi_{h_i}^1(x,\cdot)\right)_{i=1}^n$, to miarę $\Pi_{h_1,\ldots,h_{n+1}}^{1,\ldots,n+1}(x,\cdot)$ na $X^{n+1}$ możemy definiować jako jedyną miarę taką, że
\begin{align}\label{def:miara_produktowa}
\Pi_{h_1,\ldots,h_{n+1}}^{1,\ldots,n+1}(x,A\times B)=\int_A\Pi_{h_{n+1}}^1(z_n,B)\Pi_{h_1,\ldots,h_n}^{1,\ldots,n}(x,dz)\quad\text{dla }z=(z_1,\ldots,z_n),\: A\in \otimes_{i=1}^n B_{X},\: B\in B_X.
\end{align}
Ostatecznie otrzymujemy miarę $\Pi_{h_1,h_2,\ldots}^{\infty}(x,\cdot)$ na $X^{\infty}$. Należy zauważyć, że 
miary postaci (\ref{def:1-wym_rozklady}) są rozkładami marginalnymi $\Pi^{\infty}_{h_1,h_2,\ldots}(x,\cdot)$. Istnienie miary $\Pi^{\infty}_{h_1,h_2,\ldots}(x,\cdot)$ wynika z twierdzenia Kołmogorowa, mianowicie istnieje pewna przestrzeń probabilistyczna $(\Omega, \mathcal{F},\mathbb{P})$, na której możemy definiować proces stochastyczny $\xi^x$ o~rozkładzie $\Phi_{\xi^x}$ takim, że
\begin{align*}
\Phi_{\xi^x}(B)=\mathbb{P}\left(\left\{\omega\in\Omega:\xi^x(\omega)\in B\right\}\right)=\Pi^{\infty}_{h_1,h_2,\ldots}(x,B)\quad\textrm{dla }B\in \otimes_{i=1}^{\infty} B_{X}.
\end{align*}
Miara $\Pi^{\infty}_{h_1,h_2,\ldots}(x,\cdot)$ jest więc rozkładem niejednorodnego łańcucha Markowa $\xi^x$ o~ciągu funkcji przejścia $\left(\Pi^1_{h_i}\right)_{i\in N}$ i~rozkładzie początkowym $\delta_x$. Jeśli rozkład początkowy nie jest miarą $\delta_x$, tylko dowolną miarą $\mu\in M_{fin}(X)$, to piszemy
\begin{align}
\mathbb{P}_{\mu}(B)=\int_X\Pi^{\infty}_{h_1,h_2,\ldots}(x,B)\mu(dx)\quad\textrm{dla }B\in\otimes_{i=1}^{\infty} B_{X}.
\end{align}

Niech $n\in N_0$. Zauważmy następujące zależności: 
\begin{align}\label{istotna_obserwacja}
\begin{aligned}
(P^n_{\varepsilon}\mu)(A)&=\int_X\int_{\left(\bar{B}(0,\varepsilon)\right)^n}
\Pi^n_{h_1,\ldots,h_n}(x,A)\nu^{\varepsilon}(dh_1)\ldots\nu^{\varepsilon}(dh_n)\mu(dx)
\quad\text{dla } A\in B_X,\: \mu\in M_{fin}(X)
\end{aligned}
\end{align}
oraz 
\begin{align*}
\begin{aligned}
(U^n_{\varepsilon}f)(x)&=\int_{\left(\bar{B}(0,\varepsilon)\right)^n}
\int_X f(y)\Pi^n_{h_1,\ldots,h_n}(x,dy)\nu^{\varepsilon}(dh_1)\ldots\nu^{\varepsilon}(dh_n)
\quad\text{dla }x\in X,\:f\in\tilde{B}(X).
\end{aligned}
\end{align*}

\section{Własność Fellera}  
Niech $h,h_n\in\bar{B}(0,\varepsilon)$ dla $n\in N$ oraz $x,y\in X$. Weźmy $f\in\mathcal{L}$. Na mocy (\ref{def:Pi_hn}) oraz założeń~(II)~i~(IV) z~poprzedniego rozdziału otrzymujemy 
\begin{align*}
\begin{aligned}
&\left|\left\langle f,\Pi_h(x,\cdot)-\Pi_h(y,\cdot) \right\rangle\right|=\\
&=\left|\int_Xf(u)\Pi_h(x,du)-\int_Xf(v)\Pi_h(y,dv)\right|=\\
&=\left|\int_0^T f\left(T_h(x,t)\right)p(x,t)dt-\int_0^Tf\left(T_h(y,t)\right)p(y,t)dt\right|\leq\\
&\leq \int_0^T\left|f\left(T_h(x,t)\right)-f\left(T_h(y,t)\right)\right|p(x,t)dt+\int_0^T\left|f\left(T_h(y,t)\right)\right|\left|p(x,t)-p(y,t)\right|dt\leq\\
&\leq\int_0^T\varrho\left(T_h(x,t),T_h(y,t)\right)p(x,t)dt+\int_0^T\left|p(x,t)-p(y,t)\right|dt\leq\\
&\leq\Lambda_1\varrho(x,y)+\omega(\varrho(x,y)).
\end{aligned}
\end{align*}
Postępując analogicznie, wobec wzorów (\ref{def:1-wym_rozklady}) i (\ref{def:tilde_x}) otrzymujemy oszacowanie
\begin{align*}
&\left|\left\langle f,\Pi^n_{h_1,\ldots,h_n}(x,\cdot)-\Pi^n_{h_1,\ldots,h_n}(y,\cdot)\right\rangle\right|=\\
&=\left|\int_{[0,T]^n}f\left(\tilde{x}_n^x\right)p\left(\tilde{x}_{n-1}^x,t_n\right)\ldots p(x,t_1)\;dt_n\ldots dt_1
-\int_{[0,T]^n}f\left(\tilde{x}_n^y\right)p\left(\tilde{x}_{n-1}^y,t_n\right)\ldots p(y,t_1)\;dt_n\ldots dt_1\right|\leq\\
&\leq \int_{[0,T]^n}\left|f\left(\tilde{x}_n^x\right)-f\left(\tilde{x}_n^y\right)\right|p\left(\tilde{x}_{n-1}^x,t_n\right)\ldots p(x,t_1)\;dt_n\ldots dt_1+\\
&\quad+\int_{[0,T]^n}\left|f\left(\tilde{x}_n^y\right)\right|\left|p\left(\tilde{x}_{n-1}^x,t_n\right)\ldots p(x,t_1)-p\left(\tilde{x}_{n-1}^y,t_n\right)\ldots p(y,t_1)\right|\;dt_n\ldots dt_1\leq\\
&\leq \Lambda_1^n\varrho(x,y)
+\int_{[0,T]^n}\left|p\left(\tilde{x}_{n-1}^x,t_n\right)-p\left(\tilde{x}_{n-1}^y,t_n\right)\right|p\left(\tilde{x}_{n-2}^x,t_{n-1}\right)\ldots p(x,t_1)\;dt_n\ldots dt_1+\\
&\quad+\int_{[0,T]^n}\left|p\left(\tilde{x}_{n-1}^y,t_n\right)\right|\left|p\left(\tilde{x}_{n-2}^x,t_{n-1}\right)\ldots p(x,t_1)-p\left(\tilde{x}_{n-2}^y,t_{n-1}\right)\ldots p(y,t_1)\right|\;dt_n\ldots dt_1\leq\\
&\leq\ldots\leq \Lambda_1^n\varrho(x,y)+\omega\left(\Lambda_1^{n-1}\varrho(x,y)\right)+\ldots+\omega\left(\Lambda_1\varrho(x,y)\right)+\omega(\varrho(x,y))
\end{align*}
i stąd
\begin{align}\label{oszacowanie_II_IV}
\|\Pi^n_{h_1,\ldots, h_n}(x,\cdot)-\Pi^n_{h_1,\ldots,h_n}(y,\cdot)\|_{\mathcal{L}}\leq \Lambda_1^n\varrho(x,y)+\varphi(\varrho(x,y))+\omega(\varrho(x,y)),
\end{align}
gdzie funkcja $\varphi$ jest dana wzorem (\ref{def:varphi}). Powyższa nierówność oraz stwierdzenia \ref{tw_aleksandrowa}  i~\ref{stwierdzenie:Dini} implikują słabą ciągłość odwzorowania
\begin{align*}
X\ni x\mapsto \Pi^n_{h_1,\ldots, h_n}(x,\cdot)\in M_1(X).
\end{align*}
Ponieważ operator $P^n_{h_1,\ldots,h_n}$ jest regularnym operatorem Markowa, to zgodnie ze stwierdzeniem~\ref{stwierdzenie:kiedy_Markowa_jest_Fellera} jest on również fellerowski.

\begin{wniosek}\label{wn:P_epsilon-Fellera}
Wobec (\ref{istotna_obserwacja}) i niezależności oszacowania (\ref{oszacowanie_II_IV}) od $(h_n)_{n\in N}$ operator $P^n_{\varepsilon}$ posiada własność Fellera.
\end{wniosek}

\section{Miary o momentach skończonych}
 
Ustalmy punkt $\bar{x}\in X$, dla którego założenie (III) jest spełnione. Zdefiniujmy funkcję Lapunowa $V:X\to[0,\infty)$ wzorem
\begin{align*}
V(x)=\varrho(x,\bar{x})\quad \textrm{dla }x\in X.
\end{align*}

\begin{lemat}\label{lem:mom.skończ}
Niech $\Lambda_1,\Lambda_2, \Lambda_{2+\delta}<1$ będą odpowiednio dane wzorami (\ref{def:Lambda1}), (\ref{def:Lambda2}) i (\ref{def:Lambda2+delta}) oraz niech $c$ będzie dana wzorem (\ref{def:c}). Ustalmy $n\in N$ oraz ciąg $(h_i)_{i\in N}\subset H$ taki, że $h_i\in\bar{B}(0,\varepsilon)$ dla~$i\in N$. Jeśli $\mu\in M_1^j(X)$, to również $P^n_{h_1,\ldots,h_n}\mu\in M_1^j(X)$ dla $j\in\left\{1,2,2+\delta\right\}$. Ponadto 
\begin{align}\label{oszacowanie_mom_skończ}
\left\langle V^j,P^n_{h_1,\ldots,h_n}\mu \right\rangle\leq\left(\left(\Lambda_j^n\left\langle V^j,\mu\right\rangle\right)^{1/j} +\frac{c}{1-\Lambda_j^{1/j}}\right)^j\quad\text{dla }j\in\left\{1,2,2+\delta\right\},
\end{align}
gdzie $\Lambda_j$ i $c$ nie zależą od wyboru ciągu stałych $(h_i)_{i\in N}$.
\end{lemat}

\begin{proof}
Niech $\mu\in M_1^j(X)$ dla $j\in \left\{1,2,2+\delta\right\}$ oraz niech $h\in \bar{B}(0,\varepsilon)$. Z definicji funkcji przejścia~$\Pi_h$ (zob. (\ref{def:Pi_hn})) mamy
\begin{align*}
\begin{aligned}
\left(\left\langle V^j,P_h\mu\right\rangle\right)^{1/j}&=\left(\int_X\int_XV^j(y)\Pi_h(x,dy)\mu(dx)\right)^{1/j}=\\
&=\left(\int_X\int_0^T V^j(T_h(x,t))p(x,t)\,dt\,\mu(dx)\right)^{1/j}\leq\|V\circ T_h\|_{{L}^j(\varsigma)},
\end{aligned} 
\end{align*}
gdzie $\| \cdot \|_{L^j(\varsigma)}$ jest normą w przestrzeni $L^j(\varsigma)$ taką, że
\[\|f\|_{L^j(\varsigma)}=\left|\left\langle f^j,\varsigma\right\rangle\right|^{1/j}\] dla $f\in \tilde{B}(X\times[0,T])$ oraz $\varsigma\in M_{fin}(X\times[0,T])$ danej wzorem 
\[\varsigma(A)=\int_{X\times[0,T]} 1_{A}((x,t))p(x,t)\,dt\,\mu(dx) \quad\text{ dla }A\in B_{X}\otimes B_{[0,T]}.\] 
Zauważmy, że z założeń (I) i (II) otrzymujemy
\begin{align*}
\begin{aligned}
\left(V\circ T_h\right)(x,t)=\varrho\left(T_h(x,t),\bar{x}\right)\leq \varrho\left(T_h(x,t),T_h(\bar{x},t)\right)+\varrho\left(T_h(\bar{x},t),\bar{x}\right)\leq\lambda(x,t)V(x)+c,
\end{aligned}
\end{align*}
zatem z własności normy i odpowiedniego założenia o funkcji $\lambda$ (założenie (II), (II$^{\prime}$) lub (II$^{\prime\prime}$)) mamy
\begin{align*}
\begin{aligned}
\left|\left\langle V^j,P_h\mu\right\rangle\right|^{1/j}=
\left\|V\circ T_h\right\|_{L^j(\varsigma)}
&\leq\left(\int_{X\times[0,T]}\lambda^j(x,t)V^j(x)p(x,t)\,dt\,\mu(dx)\right)^{1/j}+c\\
&\leq\left(\Lambda_j\left\langle V^j,\mu\right\rangle\right)^{1/j} +c
\end{aligned}
\end{align*}
i wobec tego
\begin{align*}
\begin{aligned}
\left\langle V^j,P^n_{h_1,\ldots,h_n}\mu\right\rangle&\leq\left(\Lambda^{1/j}_j\left(\left\langle V^j,P^{n-1}_{h_1,\ldots,h_{n-1}}\mu\right\rangle\right)^{1/j} +c\right)^j\leq\\
&\leq \left(\Lambda^{2/j}_j\left(\left\langle V^j,P^{n-2}_{h_1,\ldots,h_{n-2}}\mu\right\rangle\right)^{1/j} +c\left(1+\Lambda_j^{1/j}\right)\right)^j\leq\\
&\leq\ldots\leq\left(\Lambda_j^{n/j}\left(\left\langle V^j,\mu\right\rangle\right)^{1/j}+\frac{c}{1-\Lambda_j^{1/j}}\right)^j,
\end{aligned}
\end{align*}
co kończy dowód, ponieważ z założenia $\Lambda_j^{1/j}<1$ oraz $c<\infty$.
\end{proof}

\newpage
\begin{wniosek}\label{wn:skończone_momenty_P_varepsilon_mu}
Wobec (\ref{istotna_obserwacja}) oraz niezależności oszacowania (\ref{oszacowanie_mom_skończ}) od $(h_i)_{i\in N}$ poniższe stwierdzenia wynikają bezpośrednio z lematu \ref{lem:mom.skończ}.
\begin{itemize}
\item Jeśli $\mu\in M_1^j(X)$, to również \[P^n_{\varepsilon}\mu\in M_1^j(X) \quad\text{dla  $j\in\left\{1,2,2+\delta\right\}$ oraz $n\in N$}. \]
\item Mamy 
\begin{align}\label{oszacowanie_mom_skończ_epsilon}
\begin{aligned}
\left\langle V^j,P^n_{\varepsilon}\mu \right\rangle\leq\left(\left(\Lambda_j^n\left\langle V^j,\mu\right\rangle\right)^{1/j} +{c}\left({1-\Lambda_j^{1/j}}\right)^{-1}\right)^j\quad\text{dla }j\in\left\{1,2,2+\delta\right\},\:n\in N.
\end{aligned}
\end{align}
\end{itemize}
\end{wniosek}

\chapter{Miara sprzęgająca} 

W tym rozdziale przedstawimy konstrukcję miary sprzęgającej.

\section{Definicja miary sprzęgającej}
\begin{definicja}\label{def:pojecie_miary_sprzegajacej}
Weźmy $x,y\in X$ oraz $h\in\bar{B}(0,\varepsilon)$. Niech funkcja przejścia $\Pi_h$ będzie postaci (\ref{def:Pi_hn}). Dla $\Pi^1_{h}(x,\cdot),\Pi^1_{h}(y,\cdot)\in M_1(X)$  możemy określić miarę probabilistyczną $C^1_{h}((x,y),\cdot)$ na $X^2$ taką, że
\begin{itemize}
\item $C^1_{h}((x,y),A\times X)=\Pi^1_{h}(x,A)\quad\textrm{dla }A\in B_X$,
\item $C^1_{h}((x,y),X\times B)=\Pi^1_{h}(y,B)\quad\textrm{dla }B\in B_X$.
\end{itemize}
Elementy z rodziny $\left\{C^1_{h}((x,y),\cdot):x,y\in X\right\}$ nazywamy miarami sprzęgającymi.
\end{definicja}
Klasyczną definicję miary sprzęgającej można znaleźć m.in. w \cite{villani} (definicja 1.1).

\section{Konstrukcja miary sprzęgającej} 
Ustalmy ciąg stałych $(h_i)_{i\in N}\subset\bar{B}(0,\varepsilon)$. Weźmy $\mu,\nu\in M_1^1(X)$ oraz $A,B\in B_X$. Rozważmy miary $b,b^n_{h_1,\ldots,h_n}\in M_{{fin}}\left(X^2\right)$ dla $n\in N$ takie, że
\begin{align*}
b(A\times X)=\mu(A)\text{,}\quad b(X\times B)=\nu(B)
\end{align*}
oraz 
\begin{align*}
b^n_{h_1,\ldots,h_n}(A\times X)=\left(P^n_{h_1,\ldots,h_n}\mu\right)(A)\text{,}\quad b^n_{h_1,\ldots,h_n}(X\times B)=\left(P^n_{h_1,\ldots,h_n}\nu\right)(B)\quad\text{dla }n\in N.
\end{align*}
Niech funkcja $\bar{V}:X^2\to[0,\infty)$ będzie dana wzorem
\begin{align*}
\bar{V}(x,y)=V(x)+V(y)\quad\text{dla $x,y\in X$.}
\end{align*}
Zauważmy, że z lematu \ref{lem:mom.skończ} otrzymujemy
\begin{align}\label{prop:barV} 
\left\langle\bar{V},b^n_{h_1,\ldots,h_n}\right\rangle\:\leq \:\Lambda_1\left\langle\bar{V},b^{n-1}_{h_1,\ldots,h_{n-1}}\right\rangle+2c\:\leq\:\Lambda_1^n\left\langle\bar{V},b\right\rangle+\frac{2c}{1-\Lambda_1}\quad\text{dla }n\in N.
\end{align}
Określmy funkcjonał liniowy 
\begin{align}\label{def:phi}
\phi(b)=\int_{X^2}\varrho(x,y)b(dx\times dy)\quad\text{dla }b\in M_{{fin}}^1\left(X^2\right).
\end{align}
Zgodnie z przyjętą notacją mamy
\begin{align}\label{prop:phi}
\phi(b)\leq\langle \bar{V},b\rangle.
\end{align}
Niech $x,y\in X$ oraz $A,B\in B_X$ będą dowolne. Zdefiniujmy sub-probabilistyczne funkcje przejścia na $X^2$ wzorami
\begin{align}\label{def:Qxy}
Q^1_{h}((x,y),A\times B)=\int_0^T \min\left\{p(x,t),p(y,t)\right\}\delta_{\left(T_{h}(x, t),T_{h}(y, t)\right)}(A\times B)dt\quad\text{dla }h\in\bar{B}(0,\varepsilon)
\end{align}
oraz
\begin{align}\label{def:Qxy^n}
Q^n_{h_1,\ldots,h_n}((x,y),A\times B)=\int_{X^2}Q^1_{h_n}((u,v),A\times B)Q^{n-1}_{h_1,\ldots,h_{n-1}}((x,y),du\times dv)\quad\text{dla } n\in N.
\end{align}
Miary generowane przez zdefiniowane wyżej funkcje przejścia będziemy oznaczać tym samym symbolem. Znaczenie symbolu będzie wynikać z kontekstu. Zauważmy, że
\begin{align*}
\begin{aligned}
Q^1_{h}((x,y),A\times X)&\leq\int_0^Tp(x,t)\delta_{T_{h}(x,t)}(A)dt
&=\int_0^T1_A\left(T_{h}(x,t)\right)p(x,t)dt
&=\Pi^1_{h}(x,A)
\end{aligned}
\end{align*}
i analogicznie 
$Q^1_{h}((x,y),X\times B)\leq \Pi^1_{h}(y,B)$ dla $h\in \bar{B}(0,\varepsilon)$. 
Podobnie, dla $n\in N$, mamy
\begin{align*}
Q^n_{h_1,\ldots,h_n}((x,y),A\times X)\leq\Pi^n_{h_1,\ldots,h_n}(x,A)\text{,}\\
Q^n_{h_1,\ldots,h_n}((x,y),X\times B)\leq\Pi^n_{h_1,\ldots,h_n}(y,B).
\end{align*}
Niech $b\in M_{{fin}}\left(X^2\right)$ oraz niech $Q^n_{h_1,\ldots,h_n}b$ oznacza miarę postaci
\begin{equation}\label{def:Qb}
\left(Q^n_{h_1,\ldots,h_n}b\right)(A\times B)=\int_{X^2}Q^n_{h_1,\ldots,h_n}((x,y),A\times B)\:b(dx\times dy)\quad\textrm{dla }A,B\in B_{X},\: n\in N\text{.}
\end{equation}
Zauważmy, że
\begin{align}\label{Q1Qnb}
\begin{aligned}
\left(Q^{n+1}_{h_1,\ldots,h_{n+1}}b\right)(A\times B)
&=\int_{X^2}Q_{h_1,\ldots,h_{n+1}}^{n+1}((x,y),A\times B)b(dx\times dy)=\\
&=\int_{X^2}\int_{X^2}Q^1_{h_{n+1}}((u,v),A\times B)Q^n_{h_1,\ldots,h_n}((x,y),du\times dv)b(dx\times dy)=\\
&=\int_{X^2}Q^1_{h_{n+1}}((u,v),A\times B)\left(Q^n_{h_1,\ldots,h_n}b\right)(du\times dv)
=\left(Q^1_{h_{n+1}}\left(Q_{h_1,\ldots,h_n}^nb\right)\right)(A\times B),
\end{aligned}
\end{align}
gdzie $A,B\in B_X$ i $n\in N$. Postępując analogicznie jak w (\ref{def:miara_produktowa}), możemy budować miary produktowe $Q^{1,\ldots,n}_{h_1,\ldots,h_n}((x,y),\cdot)$ na $X^n$ dla $x,y\in X$ oraz $n\in N$. Istnienie miary $Q^{\infty}_{h_1,h_2,\ldots}((x,y),\cdot)$ na całych trajektoriach wynika z~twierdzenia Kołmogorowa. Ponadto  
\[\left(Q_{h_1,h_2,\ldots}^{\infty}b\right)(A)
=\int_{X^2}Q^{\infty}_{h_1,h_2,\ldots}((x,y),A)\:b(dx\times dy)
\quad\text{ dla $b\in M_{{fin}}\left(X^2\right)$ oraz $A\in \otimes_{i=1}^{\infty} B_{X}$.}\] 

Zachodzi nierówność
\begin{align}\label{prop:phiQb}
\phi\left(Q^n_{h_1,\ldots,h_n}b\right)\leq \Lambda_1^n\phi(b)\quad\text{dla $n\in N$,\; $b\in M^1_{{fin}}\left(X^2\right)$,}
\end{align}
gdzie funkcjonał $\phi$ jest dany wzorem (\ref{def:phi}). Rzeczywiście,
\begin{align*}
\begin{aligned}
\phi(Q^n_{h_1,\ldots,h_n}b)
&=\int_{X^2}\int_{X^2}\varrho(u,v)Q^n_{h_1,\ldots,h_n}((x,y),du\times dv)\:b(dx\times dy)=\\
&=\int_{X^2}\int_{X^2}\int_0^T\int_{X^2}\varrho(u,v)\min\left\{p(\bar{u},t),p(\bar{v},t)\right\}\delta_{\left(T_{h_n}(\bar{u},t),T_{h_n}(\bar{v},t)\right)}(du\times dv)\:dt\:\times\\
&\quad \times Q^{n-1}_{h_1,\ldots,h_{n-1}}\left((x,y),d\bar{u}\times  d\bar{v}\right)\:b(dx\times dy)\leq\\
&\leq\int_{X^2}\int_{X^2}\int_0^T\varrho\left(T_{h_n}(\bar{u},t),T_{h_n}(\bar{v},t)\right)p(\bar{u},t)dt\:Q^{n-1}_{h_1,\ldots,h_{n-1}}\left((x,y),d\bar{u}\times d\bar{v}\right)\:b(dx\times dy)\leq\\
&\leq\int_{X^2}\int_{X^2}\int_0^T\varrho\left(\bar{u},\bar{v}\right)\lambda\left(\bar{u},t\right)p\left(\bar{u},t\right)dt\:Q^{n-1}_{h_1,\ldots,h_{n-1}}\left((x,y),d\bar{u}\times d\bar{v}\right)\:b(dx\times dy)\leq\\
&\leq \Lambda_1\int_{X^2}\int_{X^2}\varrho\left(\bar{u},\bar{v}\right)Q^{n-1}_{h_1,\ldots,h_{n-1}}\left((x,y),d\bar{u}\times d\bar{v}\right)\:b(dx\times dy)\leq\\
&\leq\ldots\leq \Lambda_1^n\phi(b).
\end{aligned}
\end{align*}

\begin{lemat}\label{lemma1}
Ustalmy $h\in\bar{B}(0,\varepsilon)$. Niech rodzina miar $\left\{Q^1_{h}((x,y),\cdot):x,y\in X\right\}$ na~$X^2$ będzie generowana przez funkcje przejścia postaci (\ref{def:Qxy}). Istnieje rodzina $\left\{R^1_{h}((x,y),\cdot):x,y\in X\right\}$ taka, że 
\[C^1_{h}((x,y),\cdot)=Q^1_{h_i}((x,y),\cdot)+R^1_{h}((x,y),\cdot)\quad\text{ dla $x,y\in X$}\]
oraz 
\begin{itemize}
\item[(i)] odwzorowanie $(x,y)\mapsto R^1_{h}((x,y),A\times B)$ jest mierzalne dla $A,B\in B_{X}$;
\item[(ii)] miary $R^1_{h}((x,y),\cdot)$ są nieujemne dla wszystkich $x,y\in X$;
\item[(iii)] miary $C^1_{h}((x,y),\cdot)$ są probabilistyczne dla dowolnych $x,y\in X$; 
\item[(iv)] dla wszystkich $A,B\in B_X$ oraz $x,y\in X$ mamy
\begin{align*}
C^1_{h}((x,y),A\times X)=\Pi^1_{h}(x,A), \quad
C^1_{h}((x,y),X\times B)=\Pi^1_{h}(y,B).
\end{align*}
\end{itemize}
\end{lemat}

\begin{proof}
Ustalmy $A,B\in B_{X}$. Niech
\begin{align*}
\begin{aligned}
R^1_{h}((x,y),A\times B)
=&\left(1-Q^1_{h}\left((x,y),X^2\right)\right)^{-1}\times\\
&\times\left(\Pi^1_{h}(x,A)-Q^1_{h}((x,y),A\times X)\right)
\left(\Pi^1_{h}(y,B)-Q^1_{h}((x,y),X\times B)\right),
\end{aligned}
\end{align*}
gdy $Q^1_{h}\left((x,y),X^2\right)<1$ oraz 
\begin{align*}
R^1_{h}\left((x,y),A\times B\right)=0, 
\end{align*} 
gdy $Q^1_{h}\left((x,y),X^2\right)=1$. 
Odwzorowanie można rozszerzyć do miary spełniającej wszystkie wymagane warunki $(i)-(iv)$.
\end{proof}

Wobec lematu $\ref{lemma1}$ dla $\left\{\Pi^1_{h}(x,\cdot):x\in X\right\}$ możemy konstruować rodzinę $\left\{C^1_{h}((x,y),\cdot):x,y\in X\right\}$ miar sprzęgających takich, iż $Q^1_{h}((x,y),\cdot)\leq C^1_{h}((x,y),\cdot)$ oraz $R^1_{h}((x,y),\cdot)$ są nieujemne dla dowolnych $x,y\in X$. 

Powtarzając konstrukcję miar na cylindrach i~całych trajektoriach opisaną w podrozdziale~3.2, otrzymujemy rodzinę miar probabilistycznych $\left\{C^{\infty}_{h_1,h_2\ldots}((x,y),\cdot):x,y\in X\right\}$ na $\left(X^2\right)^{\infty}$ o rozkładach brzegowych $\Pi_{h_1,h_2,\ldots}^{\infty}(x,\cdot)$ i $\Pi_{h_1,h_2,\ldots}^{\infty}(y,\cdot)$. 
Zauważmy, że każda miara $C^n_{h_1,\ldots,h_n}((x,y),\cdot)$, konstruowana według zasady (\ref{def:1-wym_rozklady}), jest rzutem na $n$-tą współrzędną miary $C^{\infty}_{h_1,h_2\ldots}((x,y),\cdot)$. Ponadto $\Big\{C^n_{h_1,\ldots,h_n}((x,y),\cdot):x,y\in X\Big\}$ pełni rolę rodziny miar sprzęgających dla $\left\{\Pi^n_{h_1,\ldots,h_n}(x,\cdot):x\in X\right\}$. Rzeczywiście, dla $x,y\in X$, $A\in B_X$, $n\in N$ oraz ciągu $(h_i)_{i\in N}$ mamy
\begin{align*}
\begin{aligned}
C^n_{h_1,\ldots,h_n}((x,y),A\times X)&=\int_{X^2}C^1_{h_n}((u,v),A\times X)C^{n-1}_{h_1,\ldots,h_{n-1}}((x,y),du\times dv)\\
&=\int_{X^2}\Pi^1_{h_n}(u,A)C^{n-1}_{h_1,\ldots,h_{n-1}}((x,y),du\times dv)=\ldots=\Pi^n_{h_1,\ldots,h_n}(x,A)
\end{aligned}
\end{align*}
i podobnie $C^n_{h_1,\ldots,h_n}((x,y),X\times B)=\Pi^n_{h_1,\ldots,h_n}(y,B)$.

Ustalmy $(x_0,y_0)\in X^2$ oraz $(h_n)_{n\in N}\subset\bar{B}(0,\varepsilon)$. Ciąg $\left(\left\{C^n_{h_1,\ldots,h_n}((x,y),\cdot):x,y\in X\right\}\right)_{n\in N}$ generuje niejednorodny łańcuch Markowa $\Psi$ na $X^2$ startujący z pary punktów $(x_0,y_0)$, z~kolei ciąg $\left(\left\{\hat{C}^n_{h_1,\ldots,h_n}((x,y,\theta),\cdot):x,y\in X, \theta\in\left\{0,1\right\}\right\}\right)_{n\in N}$ definiuje łańcuch Markowa $\hat{\Psi}$ na rozszerzonej przestrzeni $X^2\times\left\{0,1\right\}$ o~rozkładzie początkowym $\delta_{(x_0,y_0,1)}$. Dla $x,y\in X$, $A,B\in B_X$ oraz $n\in N$ mamy
\[\mathbb{P}\left(\hat{\Psi}_{n+1}\in A\times B\times\left\{1\right\}\:|\:\hat{\Psi}_n=(x,y,i),i\in\left\{0,1\right\}\right)=Q^n_{h_1,\ldots,h_n}\left((x,y),A\times B\right),\]
\[\mathbb{P}\left(\hat{\Psi}_{n+1}\in A\times B\times\left\{0\right\}\:|\:\hat{\Psi}_n=(x,y,i),i\in\left\{0,1\right\}\right)=R^n_{h_1,\ldots,h_n}\left((x,y),A\times B\right).\]
Powtarzając konstrukcję miar na cylindrach z podrozdziału 3.2 i~korzystając z~twierdzenia Kołmogorowa, otrzymujemy miarę $\hat{C}^{\infty}_{h_1,h_2,\ldots}\left((x_0,y_0,1),\cdot\right)$ na  $\left(X^2\times\left\{0,1\right\}\right)^{\infty}$, która jest stowarzyszona z łańcuchem Markowa $\hat{\Psi}$. 

Wartość oczekiwaną względem miar $C^{\infty}_{h_1,h_2,\ldots}((x_0,y_0),\cdot)$ oraz $\hat{C}^{\infty}_{h_1,h_2,\ldots}((x_0,y_0,1),\cdot)$ będziemy oznaczać symbolem $E_{x_0,y_0}$.

\chapter{Asymptotyczna stabilność i geometryczne tempo zbieżności}
W rozdziale piątym pokażemy, że operator $P_{\varepsilon}$ jest asymptotycznie stabilny, a~tempo zbieżności ciągu $(P^n_{\varepsilon}\mu)_{n\in N_0}$ dla $\mu\in M_1^1(X)$ do miary niezmienniczej w~normie Fortet-Mouriera jest geometryczne.

\section{Twierdzenia pomocnicze}\label{sec:tw_pomocnicze1}

Przypomnijmy, że stałe $\Lambda_1<1$ i $c<\infty$ są odpowiednio dane wzorami $(\ref{def:Lambda1})$ oraz $(\ref{def:c})$. Ustalmy $\varkappa\in(0,1-\Lambda_1)$ oraz zdefiniujmy zbiór
\begin{align}\label{def:K_kappa}
\begin{aligned}
K_{\varkappa}=\left\{(x,y)\in X^2: \: \bar{V}(x,y)<\varkappa^{-1}2c\right\}.
\end{aligned}
\end{align}
Niech $d:\left(X^2\right)^{\infty}\to N$ będzie momentem pierwszej wizyty w~zbiorze $K_{\varkappa}$, tzn.
\begin{align}\label{def:d}
\begin{aligned}
d\left(\left(x_n,y_n\right)_{n\in N_0}\right)=\inf\left\{n\in N:\: \left(x_n,y_n\right)\in K_{\varkappa}\right\}.
\end{aligned}
\end{align}
Gdy nie ma liczby naturalnej $n\in N$ takiej, że $\left(x_n,y_n\right)\in K_{\varkappa}$, przyjmujemy $d\left(\left(x_n,y_n\right)_{n\in N_0}\right)=\infty$.

\begin{twierdzenie}\label{theorem1}
Niech $\zeta\in(0,1)$. Istnieją stałe $C_1,C_2>0$ takie, że
\begin{align}
\begin{aligned}
E_{x_0,y_0}\left(\left(\Lambda_1+\varkappa\right)^{-\zeta d}\right)
\leq C_1\bar{V}\left(x_0,y_0\right)+C_2.
\end{aligned}
\end{align}
\end{twierdzenie} 

\begin{proof}
Ustalmy $\left(x_0,y_0\right)\in X^2$. Niech $\Psi=\left(\tilde{x}_n,\tilde{y}_n\right)_{n\in N_0}$ będzie łańcuchem Markowa startującym z~$\left(x_0,y_0\right)$ o ciągu funkcji przejścia $\left(C^1_{h_i}\right)_{i\in N}$. Przez $(\mathcal{F}_n)_{n\in N_0}\subset \mathcal{F}$ oznaczmy filtrację wyznaczoną przez łańcuch $\Psi$ w przestrzeni $\Omega$. Zdefiniujmy zbiór
\begin{align}\label{def:A_n}
{A}_n=\left\{\omega\in\Omega: \; \Psi_i=\left(\tilde{x}_i(\omega), \tilde{y}_i(\omega)\right)\notin K_{\varkappa}\; \text{ dla }\:i=1,\ldots,n-1\right\}\quad \text{dla } n\in N.
\end{align}
Zauważmy, że ${A}_{n+1}\subset {A}_n$ oraz ${A}_n\in \mathcal{F}_n$ dla $n\in N$. 
Własność $(\ref{prop:barV})$, jak również definicje zbiorów ${A}_n$ oraz $K_{\varkappa}$ (zob.  (\ref{def:A_n}) i~(\ref{def:K_kappa})) implikują, że następujące nierówności są $\mathbb{P}$-p.n. spełnione: 
\begin{align*}
1_{{A}_n}(\omega)E_{x_0,y_0}\left(\bar{V}\left(\tilde{x}_{n+1},\tilde{y}_{n+1}\right)\,|\,\mathcal{F}_n\right)(\omega)
&\leq 1_{{A}_n}(\omega)\left(\Lambda_1\bar{V}\left(\tilde{x}_n(\omega),\tilde{y}_n(\omega)\right)+2c\right)\leq\\ 	
&\leq 1_{{A}_n}(\omega)(\Lambda_1+\varkappa)\bar{V}\left(\tilde{x}_n(\omega),\tilde{y}_n(\omega)\right),
\end{align*}
wobec czego dostajemy
\begin{align*}
\begin{aligned}
\int_{A_n}\bar{V}(\tilde{x}_n,\tilde{y}_n)d\mathbb{P}
&\leq\int_{A_{n-1}}\bar{V}(\tilde{x}_n,\tilde{y}_n)d\mathbb{P}
=\int_{A_{n-1}}E_{x_0,y_0}
\left(\bar{V}(\tilde{x}_n,\tilde{y}_n)|\mathcal{F}_{n-1}\right)d\mathbb{P}\leq\\
&\leq\int_{A_{n-1}}\left(\Lambda_1\bar{V}(\tilde{x}_{n-1},\tilde{y}_{n-1})+2c\right)d\mathbb{P}
\leq(\Lambda_1+\varkappa)\int_{A_{n-1}}\bar{V}(\tilde{x}_{n-1},\tilde{y}_{n-1})d\mathbb{P}.
\end{aligned}
\end{align*}
Stosując powyższe oszacowanie $n$ razy, otrzymujemy
\begin{align}\label{oszacowanie:doP(A_n)}
\int_{A_n}\bar{V}(\tilde{x}_n,\tilde{y}_n)d\mathbb{P}\leq(\Lambda_1+\varkappa)^{n-1}\int_{A_1}\bar{V}(\tilde{x}_1,\tilde{y}_1)d\mathbb{P}\leq(\Lambda_1+\varkappa)^{n-1}\left(\Lambda_1\bar{V}(x_0,y_0)+2c\right).
\end{align}
Z definicji zbiorów $A_n$ i $K_{\varkappa}$ oraz z oszacowania (\ref{oszacowanie:doP(A_n)}) mamy
\begin{align*}
\begin{aligned}
\mathbb{P}(A_n)&\leq\int_{A_n}\varkappa(2c)^{-1}\bar{V}(\tilde{x}_n,\tilde{y}_n)d\mathbb{P}
&\leq\varkappa\left(2c(\Lambda_1+\varkappa)\right)^{-1}
(\Lambda_1+\varkappa)^n\left(\Lambda_1\bar{V}(x_0,y_0)+2c\right).
\end{aligned}
\end{align*}
Podstawmy  
\begin{align}\label{def:hat_c}
\hat{c}:=\varkappa\left(2c(\Lambda_1+\varkappa)\right)^{-1}\left(\Lambda_1\bar{V}(x_0,y_0)+2c\right).
\end{align} 
Wówczas $\mathbb{P}(A_n)\leq(\Lambda_1+\varkappa)^n\hat{c}$. 
Ustalmy $\zeta\in(0,1)$. Ponieważ zmienna losowa $d$ (zob. (\ref{def:d})) przyjmuje tylko wartości $n\in N$, otrzymujemy
\begin{align*}
\begin{aligned}
E_{x_0,y_0}\left((\Lambda_1+\varkappa)^{-\zeta d}\right)\leq\sum_{n=1}^{\infty}(\Lambda_1+\varkappa)^{-\zeta n}\mathbb{P}(A_n)
&\leq\sum_{n=1}^{\infty}(\Lambda_1+\varkappa)^{-\zeta n}(\Lambda_1+\varkappa)^n\hat{c}
&=\sum_{n=1}^{\infty}(\Lambda_1+\varkappa)^{(1-\zeta)n}\hat{c}.\\
\end{aligned}
\end{align*}
Szereg jest zbieżny, bo $\Lambda_1+\varkappa<1$. Podstawiając wartość $\hat{c}$ zgodnie z definicją (\ref{def:hat_c}) i dobierając odpowiednie stałe $C_1$ oraz $C_2$, otrzymujemy żądaną tezę.
\end{proof}
Zdefiniujmy zbiory 
\begin{align}\label{def:C_r}
C_r=\left\{(x,y)\in X^2:\; \varrho(x,y)<r\right\}\quad\text{dla }r>0.
\end{align}

\begin{lemat}\label{lemma2} 
Ustalmy $\Lambda\in(\Lambda_1,1)$. Niech $C_r$ będzie zbiorem postaci (\ref{def:C_r}) dla $r>0$. Załóżmy, że $b\in M_{{fin}}\left(X^2\right)$ jest taką miarą, że $\textrm{supp } b\subset C_r$. Istnieje $\bar{\gamma}>0$  o~własności
\begin{align}
\left(Q^n_{h_1,\ldots,h_n}b\right)\left(C_{{\Lambda}^nr}\right)\geq\bar{\gamma}^n\|b\|\quad \text{dla }n\in N.
\end{align}
\end{lemat}

\begin{proof}
Przypomnijmy, że zmienna $\tilde{x}^x_n$ jest dana wzorem (\ref{def:tilde_x}). Z~definicji $(\ref{def:Qb})$, $(\ref{def:Qxy^n})$ oraz (\ref{def:Qxy}) otrzymujemy
\begin{align*}
\begin{aligned}
&\left(Q^n_{h_1,\ldots,h_n}b\right)\left(C_{{\Lambda}^nr}\right)=\\
&=\int_{X^2}\int_{X^2}\int_0^T
\min\left\{p(u,t_n),p(v,t_n)\right\}
\delta_{\left(T_{h_n}(u,t_n),T_{h_n}(v,t_n)\right)}
\left(C_{{\Lambda}^nr}\right)\times\\ 
&\quad\times\:dt_n\:Q^{n-1}_{h_1,\ldots,h_{n-1}}((x,y),du\times dv) \: b(dx\times dy)=\\
&=\ldots=\\
&=\int_{X^2}
\Bigg(\int_{(0,T)^n}1_{C_{{\Lambda}^nr}}\left(\tilde{x}^x_n,\tilde{x}^y_n\right)
\min\left\{p(\tilde{x}^x_{n-1},t_n),p(\tilde{x}^y_{n-1},t_n)\right\}\ldots 
\min\left\{p(x,t_1),p(y,t_1)\right\}\:\times\\
&\quad\times\,dt_n\ldots dt_1\Bigg)b(dx\times dy).\\
\end{aligned}
\end{align*}
Zauważmy, że $1_{C_{{\Lambda}^nr}}\left(\tilde{x}^x_n,\tilde{x}^y_n\right)=1$ wtedy i tylko wtedy, gdy $(t_1,\ldots,t_n)\in\mathcal{T}_n$, gdzie 
\[\mathcal{T}_n:=\left\{(t_1,\ldots,t_n)\in(0,T)^n:\varrho\left(\tilde{x}^x_n,\tilde{x}^y_n\right)<{\Lambda}^nr\right\}.\]
Niech $\mathcal{T}_n':=(0,T)^n\backslash\mathcal{T}_n$. 
Z założenia (II) otrzymujemy 
\begin{align*}
\begin{aligned}
&\int_{\mathcal{T}_n'}
\varrho\left(\tilde{x}^x_n,\tilde{x}^y_n\right)
p\left(\tilde{x}^x_{n-1},t_n\right)\ldots p\left(x,t_1\right)
\:dt_n\ldots dt_1 \leq \Lambda_1^n\varrho(x,y)< \Lambda_1^nr\quad \text{dla }(x,y)\in C_r.
\end{aligned}
\end{align*}
Z drugiej strony z definicji $\mathcal{T}_n'$ dostajemy
\begin{align*}
\begin{aligned}
&{\Lambda}^n r \int_{\mathcal{T}_n'}p\left(\tilde{x}^x_{n-1},t_n\right)
\ldots p(x,t_1)dt_n\ldots dt_1
\leq\int_{\mathcal{T}_n'}\varrho\left(\tilde{x}^x_n,\tilde{x}^y_n\right)
p\left(\tilde{x}^x_{n-1},t_n\right)\ldots p(x,t_1)
dt_n\ldots dt_1,
\end{aligned}
\end{align*}
co implikuje
\begin{align*}
\begin{aligned}
\int_{\mathcal{T}_n'}p\left(\tilde{x}^x_{n-1},t_n\right)
\ldots p(x,t_1)dt_n\ldots dt_1
<\frac{\Lambda_1^n}{{\Lambda}^n}<1,
\end{aligned}
\end{align*}
a zatem mamy 
\begin{align*}
\int_{\mathcal{T}_n} p\left(\tilde{x}^x_{n-1},t_n\right)
\ldots p(x,t_1)dt_n\ldots dt_1\geq 1-\left(\frac{\Lambda_1}{\Lambda}\right)^n \geq \left(1-\frac{\Lambda_1}{\Lambda}\right)^n=:\gamma^n,
\end{align*}
gdzie ostatnia nierówność jest spełniona dla dużych $n\in N$. 
Przypomnijmy, że $M_1$ i $M_2$ są odpowiednio dane wzorami (\ref{def:M_1}) i (\ref{def:M_2}).  Korzystając z~założenia (V), mamy $|\mathcal{T}_n|\geq \left(\frac{\gamma}{M_2}\right)^n$, gdzie $|\mathcal{T}_n|$ oznacza miarę Lebesgue'a zbioru $\mathcal{T}_n$. Ostatecznie otrzymujemy
\begin{align*}
\begin{aligned}
&\left(Q^n_{h_1,\ldots,h_n}b\right)\left(C_{\Lambda^nr}\right)=\\
&=\int_{X^2}\left(\int_{\mathcal{T}_n}\min\left\{p(\tilde{x}^x_{n-1},t_n),p(\tilde{x}^y_{n-1},t_n)\right\}\ldots \min\left\{p(x,t_1),p(y,t_1)\right\}
dt_n\ldots dt_1\right)b(dx\times dy)\geq\\
&\geq\int_{X^2}M_1^n|\mathcal{T}_n|b(dx\times dy)
\geq M_1^n\left(\frac{\gamma}{M_2}\right)^n\|b\|.
\end{aligned}
\end{align*}
Podstawienie $\bar{\gamma}:=M_1 M_2^{-1}\gamma$ kończy dowód.
\end{proof}

\begin{twierdzenie}\label{theorem2}
Dla $\varkappa\in(0,1-\Lambda_1)$ istnieje liczba $n_0\in N$ taka, że
\begin{align}
\left\|Q_{h_1,h_2,\ldots}^{\infty}((x,y),\cdot)\right\|\geq\frac{1}{2}\bar{\gamma}^{n_0}\quad\text{dla $(x,y)\in K_{\varkappa}$,}
\end{align}
gdzie wielkość $\bar{\gamma}>0$ jest zdefiniowana w lemacie \ref{lemma2}.
\end{twierdzenie}

\begin{proof}
Zauważamy, że $\;\min\left\{u,v\right\}+|u-v|-u\geq 0$ dla dowolnych $u, v\in R$, zatem
\begin{align*}
\begin{aligned}
\int_0^T\left(\min\left\{p(x,t),p(y,t)\right\}+|p(x,t)-p(y,t)|-p(x,t)\right)dt\geq 0\quad\text{dla }(x,y)\in X^2.
\end{aligned}
\end{align*} 
Wobec definicji (\ref{def:Qxy}) otrzymujemy
\begin{align}\label{oszacowanie: Q_h}
\begin{aligned}
\left\|Q_{h}^1((x,y),\cdot)\right\|+\int_0^T|p(x,t)-p(y,t)|dt\geq 1\quad\text{dla }h\in\bar{B}(0,\varepsilon).
\end{aligned}
\end{align}
Niech $b\in M_{\text{fin}}\left(X^2\right)$. Korzystając kolejno z~(\ref{oszacowanie: Q_h}), warunku Diniego (założenie (IV)) i~nierówności Jensena, dostajemy
\begin{align*}
\begin{aligned}
\left\|Q^1_{h}b\right\|=\int_{X^2}Q^1_{h}\left((x,y),X^2\right)b(dx\times dy)
&=\int_{X^2}\left\|Q^1_{h}((x,y),\cdot)\right\|b(dx\times dy)\geq\\
&\geq\|b\|-\int_{X^2}\omega\left(\varrho(x,y)\right)b(dx\times dy)\geq\|b\|-\omega\left(\phi(b)\right).
\end{aligned}
\end{align*}
Następnie na mocy tożsamości $(\ref{Q1Qnb})$ mamy 
\begin{align*}
\begin{aligned}
\left\|Q^n_{h_1,\ldots,h_n}b\right\|&=\int_{X^2}Q_{h_n}^1((x,y),\cdot)\left(Q^{n-1}_{h_1,\ldots,h_{n-1}}b\right)(dx\times dy)\geq\\
&\geq\left\|Q^{n-1}_{h_1,\ldots,h_{n-1}}b\right\|-\omega\left(\phi\left(Q^{n-1}_{h_1,\ldots,h_{n-1}}b\right)\right)\geq\\
&\geq\ldots \geq\\
&\geq \|b\|-\omega(\phi(b))-\omega\left(\phi\left(Q^1_{h_1}b\right)\right)-\ldots-\omega\left(\phi\left(Q^{n-1}_{h_1,\ldots,h_{n-1}}b\right)\right).
\end{aligned}
\end{align*}
Korzystając z $(\ref{prop:phiQb})$ oraz własności funkcji $\omega$, otrzymujemy
\begin{align}\label{oszacowanie:Q^n|b|}
\begin{aligned}
\left\|Q^n_{h_1,\ldots,h_n}b\right\|\geq\|b\|-\sum_{i=1}^{n}\omega\left(\Lambda_1^{i-1}\phi(b)\right).
\end{aligned}
\end{align}
Ponieważ $\lim_{t\to 0}\varphi(t) = 0$ (stwierdzenie \ref{stwierdzenie:Dini}), to możemy wybrać liczbę $r>0$ taką, że dla $\varrho(x,y)<r$  
będzie zachodzić nierówność
\begin{align*}
\sum_{i=1}^{n}\omega\left(\Lambda_1^{i-1}\phi(b)\right)\leq\varphi\left(\Lambda_1^{-1}\phi(b)\right)<\frac{1}{2}\|b\|\quad\text{dla $n\in N$}.
\end{align*} 
Zatem, jeśli $\text{supp } b\subset C_r$, to wobec (\ref{oszacowanie:Q^n|b|}) otrzymujemy
\begin{align}\label{b}
\left\|Q_{h_1,h_2,\ldots}^{\infty}b\right\|\geq\frac{\|b\|}{2}.
\end{align}
Ustalmy $\varkappa\in(0,1-\Lambda_1)$. Łatwo zauważyć, że $K_{\varkappa}\subset C_{\varkappa^{-1} 2c}$. Dla   
$n_0:=\min\left\{n\in N:\: \Lambda_1^n(\varkappa)^{-1}2c<r\right\}$ zachodzi inkluzja $C_{\Lambda_1^{n_0}\varkappa^{-1} 2c}\subset C_r$. 
Pamiętając o następującej zależności: $Q_{h_1,\ldots,h_n,h_{n+1},\ldots,h_{n+m}}^{n+m}((x,y),\cdot)=\Big(Q^m_{h_{n+1},\ldots,h_{n+m}}Q_{h_1,\ldots,h_n}^n\Big)((x,y),\cdot)$ (por. (\ref{Q1Qnb})) oraz korzystając z~własności Markowa, mamy
\[Q_{h_1,h_2,\ldots}^{\infty}\left((x,y),X^2\right)= \left(Q^{\infty}_{h_{n_0+1},h_{n_0+2},\ldots} Q^{n_0}_{h_1,\ldots,h_{n_0}}\right)\left((x,y),X^2\right).\]
Wówczas, z $(\ref{b})$ i lematu $\ref{lemma2}$ otrzymujemy
\begin{align*}
\begin{aligned}
\left\|Q_{h_1,h_2,\ldots}^{\infty}((x,y),\cdot)\right\|
&=\left\|\left(Q^{\infty}_{h_{n_0+1},h_{n_0+2},\ldots} Q^{n_0}_{h_1,\ldots,h_{n_0}}\right)
((x,y),\cdot)\right\|
\geq\frac{1}{2}\left\|Q^{n_0}_{h_1,\ldots,h_{n_0}}((x,y),\cdot)\,|_{C_r}\right\|=\\
&=\frac{1}{2}Q^{n_0}_{h_1,\ldots,h_{n_0}}\left((x,y),C_r\right)
\geq\frac{1}{2}Q^{n_0}_{h_1,\ldots,h_{n_0}}\left((x,y),C_{\Lambda_1^{n_0}\varkappa^{-1} 2c}\right)
\geq\frac{1}{2}\bar{\gamma}^{n_0}
\end{aligned}
\end{align*}
dla $(x,y)\in K_{\varkappa}$, co kończy dowód.
\end{proof}

\begin{definicja}
Niech funkcja $\tau:(X^2\times\left\{0,1\right\})^{\infty}\to N$ będzie określona wzorem
\[\tau\left(\left(x_n,y_n,\theta_n\right)_{n\in N_0}\right)=\inf\left\{n\in N:\;\theta_k=1\;\text{ dla wszystkich }k\geq n\right\}.\]
Przyjmujemy, że $\inf\emptyset=\infty$.
\end{definicja}

\begin{twierdzenie}\label{theorem3}
Istnieją stałe $\tilde{q}\in(0,1)$ i $C_3>0$ takie, że
\[E_{x,y}\left(\tilde{q}^{-\tau}\right)\leq C_3\left(1+\bar{V}(x,y)\right)\quad \text{dla $(x,y)\in X^2$.}\]
\end{twierdzenie}

\begin{proof}
Ustalmy dowolne $\varkappa\in(0,1-\Lambda_1)$ oraz $(x,y)\in X^2$. Dla uproszczenia notacji będziemy pisać $\alpha=(\Lambda_1+\varkappa)^{-\frac{1}{2}}$. 
Niech $\hat{\Psi}=(\tilde{x}_n,\tilde{y}_n,\theta_n)_{n\in N_0}$ będzie niejednorodnym łańcuchem Markowa startującym z $(x,y,1)$ o~ciągu funkcji przejścia $\left(\hat{C}_{h_i}^{1}\right)_{i\in N}$ oraz niech $\hat{d}$ będzie momentem pierwszej wizyty łańcucha $\hat{\Psi}$ w~zbiorze $K_{\varkappa}$, tzn.
\[\hat{d}\left(\left(\tilde{x}_n,\tilde{y}_n,\theta_n\right)_{n\in N_0}\right)
=\inf\left\{n\in N:\:\left(\tilde{x}_n,\tilde{y}_n\right)\in K_{\varkappa}\right\}\]
(por. (\ref{def:d})).  Ponadto zdefiniujmy czasy kolejnych wizyt $d_n$ w~zbiorze~$K_{\varkappa}$ dla $n\in N$, które spełniają warunki
\[d_1=\hat{d},\quad d_{n+1}=d_n+\hat{d}\circ \Gamma_{d_n} \quad\text{dla }n\in N,\]
gdzie  $\Gamma_n:\left(X^2\times\left\{0,1\right\}\right)^{\infty}\to \left(X^2\times\left\{0,1\right\}\right)^{\infty}$ jest operatorem przesunięcia opisanym własnością $\Gamma_n\left(\left(\tilde{x}_k,\tilde{y}_k,\theta_k\right)_{k\in N_0}\right)=\left(\tilde{x}_{k+n},\tilde{y}_{k+n},\theta_{k+n}\right)_{k\in N_0}$.
Z twierdzenia $\ref{theorem1}$ wynika, że każda zmienna losowa $d_n$ jest $C_{h_1,h_2,\ldots}^{\infty}((x,y),\cdot)$-p.n. skończona. Mocna własność Markowa (por. (\ref{mocna_wl:E_markow})) implikuje
\begin{align}\label{wn:mocna_wl_Markowa}
E_{x,y}\left(\alpha^d\circ \Gamma_{d_n}|{\mathcal{F}}_{d_n}\right)
=E_{\left(\tilde{x}_{d_n},\tilde{y}_{d_n}\right)}\left(\alpha^d\right)
\quad \text{dla }n\in N,
\end{align}
gdzie $\mathcal{F}_{d_n}$ jest $\sigma$-ciałem generowanym przez zmienną losową $d_n$. 
Z równości (\ref{wn:mocna_wl_Markowa}), twierdzenia $\ref{theorem1}$ (dla $\zeta=1/2$) i definicji zbioru $K_{\varkappa}$ mamy 
\[E_{x,y}\left(\alpha^{d_{n+1}}\right)
=E_{x,y}\left(\alpha^{d_n}E_{\left(\tilde{x}_{d_n},\tilde{y}_{d_n}\right)}\left(\alpha^{\hat{d}}\right)\right)
\leq E_{x,y}\left(\alpha^{d_n}\right)\left(C_1\varkappa^{-1} 2c+C_2\right).\]
Ustalmy $\eta:=C_1\varkappa^{-1} 2c+C_2$. Wówczas otrzymujemy
\begin{align}\label{condition1}
E_{x,y}\left(\alpha^{d_{n+1}}\right)
\leq\eta^n E_{x,y}\left(\alpha^{\hat{d}}\right)
\leq\eta^n\left(C_1\bar{V}(x,y)+C_2\right).
\end{align}
Zdefiniujmy
\begin{align*}
\begin{aligned}
\hat{\tau}\left(\left(\tilde{x}_n,\tilde{y}_n,\theta_n\right)_{n\in N_0}\right)
=\inf\left\{n\in N:\; (\tilde{x}_n,\tilde{y}_n)\in K_{\varkappa},  \;\theta_k=1\:\text{ dla }k\geq n\right\}
\end{aligned}
\end{align*}
oraz
\begin{align*}
\begin{aligned} 
\sigma=\inf\left\{n\in N:\;\hat{\tau}=d_n\right\}.
\end{aligned}
\end{align*} 
Na mocy twierdzenia $\ref{theorem2}$ istnieje liczba $n_0\in N$ taka, że
\begin{align}\label{condition2}
\hat{C}_{h_1,h_2,\ldots}^{\infty}\left((x,y,1),\left\{\sigma>n\right\}\right)\leq\left(1-\frac{\bar{\gamma}^{n_0}}{2}\right)^n \quad\text{dla }n\in N.
\end{align}
Niech $p>1$. Korzystając z nierówności H\"{o}ldera i oszacowań $(\ref{condition1})$ oraz $(\ref{condition2})$, otrzymujemy
\begin{align*}
\begin{aligned}
E_{x,y}\left(\alpha^{\frac{\hat{\tau}}{p}}\right)
&\leq\sum_{k=1}^{\infty}E_{x,y}\left(\alpha^{\frac{d_k}{p}}1_{\{\sigma=k\}}\right)\leq\\
&\leq\sum_{k=1}^{\infty}\left(E_{x,y}\left(\alpha^{d_k}\right)\right)^{\frac{1}{p}}\left(\hat{C}_{h_1,h_2,\ldots}^{\infty}\left((x,y,1),\left\{\sigma=k\right\}\right)\right)^{(1-\frac{1}{p})}\leq\\
&\leq\left(C_1\bar{V}(x,y)+C_2\right)^{\frac{1}{p}}\eta^{-\frac{1}{p}}\sum_{k=1}^{\infty}\eta^{\frac{k}{p}}\left(1-\frac{1}{2}\bar{\gamma}^{n_0}\right)^{(k-1)(1-\frac{1}{p})}=\\
&=\left(C_1\bar{V}(x,y)+C_2\right)^{\frac{1}{p}}\eta^{-\frac{1}{p}}\left(1-\frac{1}{2}\bar{\gamma}^{n_0}\right)^{-(1-\frac{1}{p})}\sum_{k=1}^{\infty}\left(\left(\frac{\eta}{1-\frac{1}{2}\bar{\gamma}^{n_0}}\right)^{\frac{1}{p}}\left(1-\frac{1}{2}\bar{\gamma}^{n_0}\right)\right)^k.
\end{aligned}
\end{align*}
Dla odpowiednio dużych $p$ oraz $\tilde{q}=\alpha^{-\frac{1}{p}}$, mamy
\[E_{x,y}\left(\tilde{q}^{-\hat{\tau}}\right)=E_{x,y}\left(\alpha^{\frac{\hat{\tau}}{p}}\right)\leq\left(1+\bar{V}(x,y)\right)C_3,\]
gdzie $C_3$ jest pewną stałą. Ponieważ $\tau\leq\hat{\tau}$, dowód jest ukończony.
\end{proof}

\begin{lemat}\label{theorem_new}
Niech $f\in\mathcal{L}$. Istnieją stałe $q\in(0,1)$ oraz $C_5>0$ takie, że
\begin{align*}
\begin{aligned}
\int_{X^2}|f(u)-f(v)|\left(\Pi_{X^2}^*\Pi^*_n\hat{C}^{\infty}_{h_1,h_2\ldots}((x,y,1),\cdot)\right)(du\times dv)\leq q^nC_5\left(1+\bar{V}(x,y)\right)\quad\text{dla $x,y\in X$, $n\in N$, }
\end{aligned}
\end{align*}
gdzie $\Pi^*_n:\left(X^2\times\left\{0,1\right\}\right)^{\infty}\to X^2\times\left\{0,1\right\}$ są rzutami na $n$-te współrzędne oraz $\Pi^*_{X^2}:X^2\times\left\{0,1\right\}\to X^2$ jest projekcją na $X^2$.
\end{lemat}

\begin{proof}
Niech $n\in N$. Zdefiniujmy zbiory
\begin{align*}
\begin{aligned}
A_{\frac{n}{2}}=\left\{t\in\left(X^2\times\left\{0,1\right\}\right)^{\infty}:\:\tau(t)\leq\frac{n}{2}\right\},\\
B_{\frac{n}{2}}=\left\{t\in\left(X^2\times\left\{0,1\right\}\right)^{\infty}:\:\tau(t)>\frac{n}{2}\right\}.
\end{aligned}
\end{align*}
Zauważmy, że $A_{\frac{n}{2}}\cap B_{\frac{n}{2}}=\emptyset$ oraz $A_{\frac{n}{2}}\cup B_{\frac{n}{2}}=\left(X^2\times\left\{0,1\right\}\right)^{\infty}$, zatem
\[\hat{C}_{h_1,h_2,\ldots}^{\infty}((x,y,1),\cdot)=\hat{C}_{h_1,h_2,\ldots}^{\infty}((x,y,1),\cdot)|_{A_{\frac{n}{2}}}+\hat{C}_{h_1,h_2,\ldots}^{\infty}((x,y,1),\cdot)|_{B_{\frac{n}{2}}}\quad\text{dla $n\in N$.}\]
Wówczas
\begin{align*}
\begin{aligned}
&\int_{X^2}|f(u)-f(v)|\left(\Pi_{X^2}^*\Pi_n^*\hat{C}_{h_1,h_2,\ldots}^{\infty}((x,y,1),\cdot)|_{A_{\frac{n}{2}}}\right)(du\times dv)+\\
&\quad +\int_{X^2}|f(u)-f(v)|\left(\Pi_{X^2}^*\Pi_n^*
\hat{C}_{h_1,h_2,\ldots}^{\infty}((x,y,1),\cdot)|_{B_{\frac{n}{2}}}\right)(du\times dv)\leq\\
&\leq\int_{X^2}\varrho(u,v)\left(\Pi_{X^2}^*\Pi_n^*\hat{C}_{h_1,h_2,\ldots}^{\infty}((x,y,1),\cdot)|_{A_{\frac{n}{2}}}\right)(du\times dv)+2\hat{C}_{h_1,h_2,\ldots}^{\infty}\left((x,y,1),B_{\frac{n}{2}}\right).
\end{aligned}
\end{align*}
Stosując wielokrotnie własność $(\ref{prop:phiQb})$, otrzymujemy
\begin{align*}
\begin{aligned}
\int_{X^2}\varrho(u,v)
\left(\Pi_{X^2}^*\Pi_n^*\hat{C}_{h_1,h_2,\ldots}^{\infty}((x,y,1),\cdot)|_{A_{\frac{n}{2}}}\right)(du,dv)
&=\phi\left(\Pi_{X^2}^*\Pi_n^*\hat{C}_{h_1,h_2,\ldots}^{\infty}((x,y,1),\cdot)|_{A_{\frac{n}{2}}}\right)\leq\\
&\leq \Lambda_1^{\lfloor\frac{n}{2}\rfloor}\phi\left(\Pi_{X^2}^*\Pi_{\lfloor\frac{n+1}{2}\rfloor}^*\hat{C}_{h_1,h_2,\ldots}^{\infty}((x,y,1),\cdot)|_{A_{\frac{n}{2}}}\right).
\end{aligned}
\end{align*}
Dzięki $(\ref{prop:barV})$ oraz $(\ref{prop:phi})$ mamy
\begin{align*}
\begin{aligned}
\phi\left(\Pi_{X^2}^*\Pi_{\lfloor\frac{n+1}{2}\rfloor}^*\hat{C}_{h_1,h_2,\ldots}^{\infty}
((x,y,1),\cdot)|_{A_{\frac{n}{2}}}\right)
\leq \Lambda_1^{\lfloor\frac{n+1}{2}\rfloor}\bar{V}(x,y)+\frac{2c}{1-\Lambda_1}
\end{aligned}
\end{align*}
i wobec tego
\begin{align*}
\begin{aligned}
&\int_{X^2}|f(u)-f(v)|\left(\Pi_{X^2}^*\Pi^*_n
\hat{C}^{\infty}_{h_1,h_2\ldots}((x,y,1),\cdot)\right)(du\times dv)\leq\\
&\quad \leq \Lambda_1^{\lfloor\frac{n}{2}\rfloor}
\left(\Lambda_1^{\lfloor\frac{n+1}{2}\rfloor}\bar{V}(x,y)+\frac{2c}{1-\Lambda_1}\right)
+2\hat{C}_{h_1,h_2,\ldots}^{\infty}\left((x,y,1),B_{\frac{n}{2}}\right).
\end{aligned}
\end{align*}
Z nierówności Czebyszewa oraz twierdzenia $\ref{theorem3}$ otrzymujemy
\begin{align*}
\begin{aligned}
\hat{C}_{h_1,h_2,\ldots}^{\infty}\left((x,y,1),B_{\frac{n}{2}}\right)=\hat{C}_{h_1,h_2,\ldots}^{\infty}\left((x,y,1),\left\{\tau>\frac{n}{2}\right\}\right)
&=\hat{C}_{h_1,h_2,\ldots}^{\infty}
\left((x,y,1),\left\{\tilde{q}^{-\tau}\geq\tilde{q}^{-\frac{n}{2}}\right\}\right)\leq\\
&\leq\frac{E_{x,y}\left(\tilde{q}^{-\tau}\right)}{\tilde{q}^{-\frac{n}{2}}}\leq\tilde{q}^{\frac{n}{2}}C_3\left(1+\bar{V}(x,y)\right),
\end{aligned}
\end{align*}
gdzie $\tilde{q}\in(0,1)$ oraz $C_3>0$ są zadane w twierdzeniu $\ref{theorem3}$. 
Ostatecznie otrzymujemy
\[\int_{X^2}|f(u)-f(v)|\left(\Pi_{X^2}^*\Pi^*_n\hat{C}^{\infty}_{h_1,h_2\ldots}((x,y,1),\cdot)\right)(du\times dv)\leq \Lambda_1^{\lfloor\frac{n}{2}\rfloor}C_4\left(1+\bar{V}(x,y)\right)+2\tilde{q}^{\frac{n}{2}}C_3\left(1+\bar{V}(x,y)\right),\]
gdzie $C_4:=\max\left\{\Lambda_1^{1/2},(1-\Lambda_1)^{-1}2c\right\}$. Podstawiając $q:=\max\left\{\Lambda_1^{1/2},\tilde{q}^{1/2}\right\}$ oraz $C_5:=C_4+2C_3$, dostajemy tezę twierdzenia.
\end{proof}

\begin{wniosek}\label{wniosek_g}
Jeśli $g\in B(X)$ jest funkcją lipschitzowską ze stałą Lipschitza $L_g>0$, to na mocy lematu~\ref{theorem_new} istnieją $q\in(0,1)$ oraz $C_5>0$ takie, że
\begin{align*}
\begin{aligned}
\int_{X^2}|g(u)-g(v)|\left(\Pi_{X^2}^*\Pi^*_n\hat{C}^{\infty}_{h_1,h_2\ldots}((x,y,1),\cdot)\right)(du\times dv)\leq Gq^nC_5\left(1+\bar{V}(x,y)\right)\quad\text{dla $x,y\in X$, $n\in N$},
\end{aligned}
\end{align*}
gdzie $G:=\max\left\{L_g,\sup_{x\in X}|g(x)|\right\}$.
\end{wniosek}

\begin{twierdzenie}\label{theorem4}
Istnieją stałe $q\in(0,1)$ oraz $C_5>0$ takie, że
\begin{align}\label{oszacowanie_funkcji_przejścia_w_normie}
\left\|\Pi^n_{h_1,\ldots,h_n}(x,\cdot)-\Pi^n_{h_1,\ldots,h_n}(y,\cdot)\right\|_{\mathcal{L}}\leq q^nC_5\left(1+\bar{V}(x,y)\right)\quad\text{dla $x,y\in X$ oraz $n\in N$.}
\end{align}
\end{twierdzenie}

\begin{proof}
Twierdzenie jest konsekwencją lematu \ref{theorem_new}. Wystarczy zauważyć, że 
\begin{align*}
\begin{aligned}
\left\|\Pi^n_{h_1,\ldots,h_n}(x,\cdot)-\Pi^n_{h_1,\ldots,h_n}(y,\cdot)\right\|_{\mathcal{L}}&=\sup_{f\in\mathcal{L}}
\left|\int_{X}f(z)\left(\Pi^n_{h_1,\ldots,h_n}(x,\cdot)-\Pi^n_{h_1,\ldots,h_n}(y,\cdot)\right)(dz)\right|=\\
&=\sup_{f\in\mathcal{L}}\left|\int_{X^2}\left(f(z_1)-f(z_2)\right)
\left(\Pi_{X^2}^*\Pi_n^*\hat{C}_{h_1,h_2,\ldots}^{\infty}((x,y,1),\cdot)\right)\left(dz_1\times dz_2\right)\right|,
\end{aligned}
\end{align*}
co implikuje
\[\left\|\Pi^n_{h_1,\ldots,h_n}(x,\cdot)-\Pi^n_{h_1,\ldots,h_n}(y,\cdot)\right\|_{\mathcal{L}}\leq q^nC_5\left(1+\bar{V}(x,y)\right).\]
Dowód został ukończony.
\end{proof}

\section{Oszacowanie tempa zbieżności w modelu}

\begin{twierdzenie}\label{GTZ}
Ustalmy $\varepsilon\in[0,\varepsilon_*]$ dla $\varepsilon_*<\infty$.  Niech $X$ będzie domkniętym podzbiorem przestrzeni Banacha $H$.  Niech $\mu\in M_1^1(X)$. Przy założeniach (I)-(VI) istnieje jedyna miara niezmiennicza $\mu_*\in M_1^1(X)$ oraz stałe $C:=C(\mu)>0$ i $q\in[0,1)$ takie, że
\begin{align}\label{wl:GTZ}
\left\|P^n_{\varepsilon}\mu-\mu_*\right\|_{\mathcal{L}}\leq Cq^n\quad\textrm{ dla }n\in N.
\end{align}
\end{twierdzenie}

\begin{proof}

Operator Markowa $P_{\varepsilon}$ można wyrazić za pomocą funkcji przejścia $\Pi^1_{h}$ (zob. (\ref{Pi_epsilon-as-Pi_h}) oraz (\ref{P_epsilon-as-P_h})). Wówczas
\[P^n_{\varepsilon}\mu(\cdot)=\int_X \int_{\bar{B}(0,\varepsilon)}\ldots\int_{\bar{B}(0,\varepsilon)}\Pi^n_{h_1,\ldots,h_n}(x,\cdot)\nu^{\varepsilon}(dh_1)\ldots\nu^{\varepsilon}(dh_n)\mu(dx).\]
Ponadto z wniosków \ref{wn:P_epsilon-Fellera} i \ref{wn:skończone_momenty_P_varepsilon_mu} wiemy, że $P_{\varepsilon}$ jest operatorem Fellera, a $P_{\varepsilon}\mu\in M_1^1(X)$ dla $\mu\in M_1^1(X)$ oraz 
\begin{align}\label{prop:<V,P^n_epsilon_mu>}
\left\langle V,P^n_{\varepsilon}\mu\right\rangle \leq \Lambda_1^n\left\langle V,\mu\right\rangle +\frac{c}{1-\Lambda_1}.
\end{align}
Niech $\mu_1,\mu_2\in M_{1}^1(X)$ oraz $f\in\mathcal{L}$. Wówczas 
\begin{align*}
\begin{aligned}
&\left|\left\langle f,P^n_{\varepsilon}\mu_1-P^n_{\varepsilon}\mu_2\right\rangle \right|=\\
&=\Bigg|\int_X \int_{\left(\bar{B}(0,\varepsilon)\right)^n}\int_X
f(z)\Pi^n_{h_1,\ldots,h_n}(x,dz)\nu^{\varepsilon}(dh_1)\ldots\nu^{\varepsilon}(dh_n)\mu_1(dx)-\\
&\quad 
-\int_X\int_{\left(\bar{B}(0,\varepsilon)\right)^n}\int_X
f(z)\Pi^n_{h_1,\ldots,h_n}(y,dz)
\nu^{\varepsilon}(dh_1)\ldots\nu^{\varepsilon}(dh_n)\mu_2(dy)\Bigg|=\\
&=\Bigg|\int_X\left(\int_X\int_{\left(\bar{B}(0,\varepsilon)\right)^n}\int_X
f(z)\Pi^n_{h_1,\ldots,h_n}(x,dz)
\nu^{\varepsilon}(dh_1)\ldots\nu^{\varepsilon}(dh_n)\mu_1(dx)\right)\mu_2(dy)-\\  
&\quad -\int_X\left(\int_X\int_{\left(\bar{B}(0,\varepsilon)\right)^n}\int_Xf(z)
\Pi^n_{h_1,\ldots,h_n}(y,dz)
\nu^{\varepsilon}(dh_1)\ldots\nu^{\varepsilon}(dh_n)\mu_2(dy)\right)\mu_1(dx)\Bigg|.\\
\end{aligned}
\end{align*}
Zauważmy, że oszacowanie (\ref{oszacowanie_funkcji_przejścia_w_normie}) z twierdzenia \ref{theorem4} nie zależy od wyboru $(h_i)_{i\in N}$, zatem
\begin{align*} 
\begin{aligned}
&\left|\left\langle f,P^n_{\varepsilon}\mu_1-P^n_{\varepsilon}\mu_2\right\rangle \right|\leq\\
&\leq\int_{X^2}\Bigg(\int_{\left(\bar{B}(0,\varepsilon)\right)^n}
\left|\int_Xf(z)\Pi^n_{h_1,\ldots,h_n}(x,dz)-\int_Xf(z)\Pi^n_{h_1,\ldots,h_n}(y,dz)\right|
\nu^{\varepsilon}(dh_1)\ldots\nu^{\varepsilon}(dh_n)\Bigg)\mu_1(dx)\mu_2(dy)\leq\\
&\leq\int_{X^2}\Bigg(\int_{\left(\bar{B}(0,\varepsilon)\right)^n}\left\|\Pi^n_{h_1,\ldots,h_n}(x,\cdot)-\Pi^n_{h_1,\ldots,h_n}(y,\cdot)\right\|_{\mathcal{L}}\;\nu^{\varepsilon}(dh_1)\ldots\nu^{\varepsilon}(dh_n)\Bigg)\mu_1(dx)\mu_2(dy)\leq\\
&\leq q^nC_5\int_{X^2}\Bigg(\int_{\left(\bar{B}(0,\varepsilon)\right)^n}
\left(1+\bar{V}(x,y)\right)\nu^{\varepsilon}(dh_1)\ldots\nu^{\varepsilon}(dh_n)\Bigg)\mu_1(dx)\mu_2(dy)=\\
&=q^nC_5\int_{X^2}\left(1+\bar{V}(x,y)\right)\mu_1(dx)\mu_2(dy),
\end{aligned}
\end{align*}
gdzie $\mu_1,\mu_2\in M^1_{1}(X)$. Ponieważ funkcja $f\in\mathcal{L}$ była dowolna, to
\begin{align*}
\begin{aligned}
\left\|P^n_{\varepsilon}\mu_1-P^n_{\varepsilon}\mu_2\right\|_{\mathcal{L}}\leq q^nC_5\int_X\int_X\left(1+\bar{V}(x,y)\right)\mu_1(dx)\mu_2(dy).
\end{aligned}
\end{align*}
Operator Markowa $P_{\varepsilon}$ jest więc mieszający.

Podstawmy $\mu_1:=\mu\in M_1^1(X)$ oraz $\mu_2:=P^m_{\varepsilon}\mu\in M_1^1(X)$  dla $m\in N$.  
Mamy
\begin{align*}
\begin{aligned}
\left\|P^n_{\varepsilon}\mu-P^{n+m}_{\varepsilon}\mu\right\|_{\mathcal{L}}
&\leq q^nC_5\int_X\int_X\left(1+\bar{V}(x,y)\right)\mu(dx)P^m_{\varepsilon}\mu(dy)=\\
&=q^nC_5\left(1+\left\langle V,\mu\right\rangle+\left\langle V,P^m_{\varepsilon}\mu\right\rangle\right)\leq q^nC_6
\end{aligned}
\end{align*}
dla pewnej stałej $C_6$, zatem $\left(P^n_{\varepsilon}\mu\right)_{n\in N_0}$ jest ciągiem Cauchy'ego w przestrzeni $\left(M_1(X),\|\cdot\|_{\mathcal{L}}\right)$. Z~zupełności przestrzeni $\left(M_1(X), \|\cdot\|_{\mathcal{L}}\right)$ (stwierdzenie \ref{stwierdzenie:zupelność}) wynika, że ciąg $\left(P^n_{\varepsilon}\mu\right)_{n\in N_0}$ jest w~niej zbieżny. Zdefiniujmy $\mu_*(\mu)$ jako słabą granicę ciągu $\left(P^n_{\varepsilon}\mu\right)_{n\in N_0}$, przy $n\to\infty$. Miara $\mu_*(\mu)$ jest niezmiennicza, bo $P_{\varepsilon}$ spełnia własność Fellera (stwierdzenie \ref{stwierdzenie:Feller_daje_niezmienniczość}). Ponadto miara $\mu_*:=\mu_*(\mu)$ jest jedyna (stwierdzenie \ref{stwierdzenie:jedyność_miary}).

Weźmy ciąg funkcji $(V_k)_{k\in N_0}$, $V_k:X\to R$, zadanych wzorem $V_k(y)=\min\left\{k,V(y)\right\}$ dla $k\in N_0$ oraz $y\in X$. Ustalamy $x\in X$. Z pierwszej części dowodu 
wiemy już, że $\lim_{n\to\infty}\left\|P^n_{\varepsilon}\delta_x-\mu_*\right\|_{\mathcal{L}}=0$, co jest równoważne temu, że 
$P^n_{\varepsilon}\delta_x\xrightarrow{w}\mu_*$, gdy $n\to\infty$ (stwierdzenie \ref{tw_aleksandrowa}). Dla wszystkich $k\in N_0$ mamy $V_k\in C(X)$, zatem
\[\lim_{n\to\infty}\langle V_k,P^n_{\varepsilon}\delta_x\rangle=\langle V_k,\mu_*\rangle.\]
Zauważmy, że zgodnie z (\ref{prop:<V,P^n_epsilon_mu>})  otrzymujemy 
\[\left\langle V_k,P^n_{\varepsilon}\delta_x\right\rangle 
=\Lambda_1^n\left\langle V_k,\delta_x\right\rangle +(1-\Lambda_1)^{-1}c 
\leq \Lambda_1^nV_k(x)+\left(1-\Lambda_1\right)^{-1}c\quad\text{dla $n\in N_0$}.\]
Ponadto 
\[\left\langle V_k,\mu_*\right\rangle =\lim_{n\to\infty}\left\langle V_k,P^n_{\varepsilon}\delta_x\right\rangle \leq (1-\Lambda_1)^{-1}c,\]
więc ciąg $\left(\left\langle V_k,\mu_*\right\rangle \right)_{k\in N_0}$ jest wspólnie ograniczony. 
Ponieważ funkcje $(V_k)_{k\in N_0}$ są nieujemne i niemalejące, możemy skorzystać z twierdzenia Lebesgue'a o zbieżności monotonicznej i otrzymać
\[\langle V,\mu_*\rangle=\lim_{k\to\infty}\langle V_k,\mu_*\rangle.\]
Pokazaliśmy więc, że $\mu_*\in M_1^1(X)$.

Pamiętając, że $\bar{V}(x,y)=V(x)+V(y)$, asymptotyczna stabilność operatora $P_{\varepsilon}$, jak również geometryczne tempo zbieżności do jedynej miary niezmienniczej $\mu_*\in M_1^1(X)$ wynikają z~poniższego oszacowania
\begin{align*}
\begin{aligned}
\left\|P^n_{\varepsilon}{\mu}-\mu_*\right\|_{\mathcal{L}}
\leq\int_X\int_X q^nC_5\left(1+\bar{V}(x,y)\right)\mu_*(dy)\mu(dx)
\leq q^nC,
\end{aligned}
\end{align*}
gdzie $C:=C_5 \left( 1+ \langle V,\mu \rangle + \langle V,\mu_* \rangle \right)$. 
Ponieważ stała $C$ zależy tylko i wyłącznie od miary $\mu\in M_1^1(X)$, dowód jest ukończony.
\end{proof}

\chapter{Centralne twierdzenie graniczne}

Celem tego rozdziału jest udowodnienie CTG dla uogólnionego modelu cyklu komórkowego.

\section{Centralne twierdzenia graniczne dla stacjonarnych łańcuchów Markowa}

Dany jest łańcuch Markowa $(x_n)_{n\in N_0}$ generowany przez funkcję przejścia $\Pi$, startujący z miary niezmienniczej $\mu_*$, o wartościach w przestrzeni mierzalnej $(\mathcal{X},\mathcal{B})$. Niech $L^2(\mu_*)$ będzie przestrzenią funkcji $g:\mathcal{X}\to R$ całkowalnych z~kwadratem, dla których $\|g\|^2_{L^2(\mu_*)}=\int_{\mathcal{X}}g^2d\mu_*<\infty$. Podprzestrzeń przestrzeni $L^2(\mu_*)$, której elementy spełniają własność $\left\langle g,\mu_*\right\rangle=0$, będziemy oznaczać symbolem $L_0^2(\mu_*)$. Niech $g\in L_0^2(\mu_*)$. Zdefiniujmy
\begin{align*}
\begin{aligned}
S_n^*=\frac{g(x_0)+\ldots +g(x_{n-1})}{\sqrt{n}}\quad\text{dla }n\in N.
\end{aligned}
\end{align*}

\begin{twierdzenie}\label{tw:wn_woodr}
Niech $g\in L_0^2(\mu_*)$. Jeśli spełniony jest warunek
\begin{align}\label{condition_woodroofe_original}
\begin{aligned}
\sum_{n=1}^{\infty}n^{-3/2}\left(\int_X\left(\sum_{k=0}^{n-1}\int_{\mathcal{X}} g(y) \Pi^k(x,dy) \right)^2\mu_*(dx)\right)^{1/2}<\infty,
\end{aligned}
\end{align}
to istnieje skończona granica
\begin{align*}
\begin{aligned}
\sigma^2:=\sigma^2(g)=\lim_{n\to\infty}E_{\mu_*}\left(\left(S^*_n\right)^{2}\right)
\end{aligned}
\end{align*}
oraz ciąg zmiennych losowych $\left({S}^*_n\right)_{n\in N}$ zbiega według rozkładu do zmiennej losowej o rozkładzie normalnym $\mathcal{N}\left(0,\sigma^2\right)$ (wniosek 1, \cite{woodr}).
\end{twierdzenie}

\section{Zastosowanie do badanego modelu}\label{zastosowanie_woodr}
Ustalmy $\varepsilon\in[0,\varepsilon_*]$ dla $\varepsilon_*<\infty$. 
\begin{lemat}\label{second_moment}
Niech $\mu\in M_1^2(X)$. Jeśli $\mu_*$ jest granicą ciągu $\left(P^n_{\varepsilon}\mu\right)_{n\in N_0}$ zbieżnego w normie Fortet-Mouriera, to $\mu_*\in M_1^2(X)$.
\end{lemat}
\begin{proof}
Niech $\mu\in M_1^2(X)$. Ustalmy $x\in X$ oraz $n\in N_0$. Przypomnijmy, że stałe $\Lambda_2<1$ i~$c<\infty$ są odpowiednio zadane wzorami (\ref{def:Lambda2}) oraz (\ref{def:c}). Z wniosku \ref{wn:skończone_momenty_P_varepsilon_mu} wiemy, że $P^n_{\varepsilon}\mu\in M_1^2(X)$. Ponadto (por. (\ref{oszacowanie_mom_skończ_epsilon})) mamy
\[\left\langle V^2,P^n_{\varepsilon}\mu\right\rangle 
\leq \left(\left(\Lambda_2^n\left\langle V^2,\mu\right\rangle\right)^{1/2}+c\left(1-\Lambda_2^{1/2}\right)^{-1}\right)^{2}.\] 
Niech $\tilde{V}_k(y)=\min\left\{k,V^2(y)\right\}$ dla wszystkich $y\in X$ oraz $k\in N_0$. Z twierdzenia \ref{GTZ} oraz stwierdzenia \ref{tw_aleksandrowa} wiemy, że $P^n_{\varepsilon}\mu \xrightarrow{w} \mu_*$, gdy $n\to\infty$. Zauważmy, że $\tilde{V}_k\in C(X)$ dla dowolnych $k\in N_0$, zatem
\[\left\langle \tilde{V}_k,\mu_*\right\rangle=\lim_{n\to\infty}\left\langle \tilde{V}_k,P^n_{\varepsilon}\mu\right\rangle \leq c^2\left(1-\sqrt{\Lambda_2}\right)^{-2}\]
i wobec tego ciąg $\left(\left\langle \tilde{V}_k,\mu_*\right\rangle \right)_{k\in N_0}$ jest wspólnie ograniczony. 
Ponieważ $\left(\tilde{V}_k\right)_{k\in N_0}$ jest niemalejącym ciągiem funkcji nieujemnych, to z~twierdzenia Lebesgue'a o~zbieżności monotonicznej otrzymujemy
\[\left\langle V^2,\mu_*\right\rangle =\lim_{k\to\infty}\left\langle \tilde{V}_k,\mu_*\right\rangle\leq c^2(1-\sqrt{\Lambda_2})^{-2}\]
i stąd $\mu_*\in M_1^2(X)$.
\end{proof}

Ustalmy łańcuch Markowa $(x_n)_{n\in N_0}$ generowany przez funkcję przejścia $\Pi_{\varepsilon}$ oraz rozkład początkowy $\mu\in M_1^2(X)$. 
Niech $S_n^{\mu}$ będzie zmienną losową zadaną wzorem 
\begin{align}\label{def:eta_n^mu}
S_n^{\mu}=\frac{g(x_0)+\ldots+g(x_{n-1})}{\sqrt{n}}\quad\text{dla }n\in N
\end{align}
oraz niech $\Phi_{S_n^{\mu}}$ oznacza rozkład $S_n^{\mu}$. W szczególności niech zmienne $S_n^*$ oraz $S_n^x$ będą określone dla łańcuchów Markowa o~funkcji przejścia $\Pi_{\varepsilon}$ i rozkładach początkowych $\mu_*$ oraz $\delta_x$. Załóżmy, że funkcja $g\in B(X)$ spełnia warunek Lipschitza ze~stałą $L_g$ oraz $\left\langle g,\mu_*\right\rangle =0$. Łatwo zauważyć, że $\left\langle g^2,\mu_*\right\rangle<\infty$.

\begin{lemat}\label{lem:condition_woodroofe}
Niech $g\in B(X)$ będzie funkcją spełniającą warunek Lipschitza ze stałą $L_g$. Ponadto niech $\left\langle g,\mu_*\right\rangle =0$. Wtedy spełniony jest warunek (\ref{condition_woodroofe_original}) z twierdzenia \ref{tw:wn_woodr}, czyli
\begin{align}\label{condition_woodroofe}
\begin{aligned}
\sum_{n=1}^{\infty}n^{-3/2}\left(\int_X\left(\sum_{k=0}^{n-1}\left\langle g,P^{k}_{\varepsilon}\delta_x\right\rangle \right)^2\mu_*(dx)\right)^{1/2}<\infty.
\end{aligned}
\end{align}
\end{lemat}
\begin{proof}
Stosując lemat \ref{theorem_new} i wniosek \ref{wniosek_g}, otrzymujemy 
\begin{align*}
\begin{aligned}
\sum_{k=0}^{n-1}\left\langle g,P^{k}_{\varepsilon}\delta_x\right\rangle &=\sum_{k=0}^{n-1}\left(\left\langle g,P^{k}_{\varepsilon}\delta_x\right\rangle -\left\langle g,\mu_*\right\rangle \right)=\\
&=\sum_{k=0}^{n-1}\int_X\left(\int_Xg(z)\left(\Pi^k_{\varepsilon}(x,\cdot)-\Pi^k_{\varepsilon}(y,\cdot)\right)(dz)\right)\mu_*(dy)=\\
&=\sum_{k=0}^{n-1}\int_X\left(\int_{X^2}\left(g(z_1)-g(z_2)\right)\left(\Pi^*_{X^2}\Pi^*_k\hat{C}^{\infty}_{h_1,h_2,\ldots}((x,y,1),\cdot)\right)\left(dz_1\times dz_2\right)\right)\mu_*(dy)\leq\\
&\leq \sum_{k=0}^{n-1}Gq^nC_5\int_{X^2}\left(1+\bar{V}(x,y)\right)\mu_*(dy),
\end{aligned}
\end{align*}
dla $x\in X$, $n\in N$, a zatem 
\[\sum_{k=0}^{n-1}\left\langle g,P^{k}_{\varepsilon}\delta_x\right\rangle \leq GC_5\frac{1-q^n}{1-q}\int_{X^2}\left(1+\bar{V}(x,y)\right)\mu_*(dy)\leq C_9\left(1+V(x)\right),\]
gdzie $C_9:=GC_5(1-q)^{-1}\left(1+\left\langle V,\mu_*\right\rangle\right)$. Pamiętając, że $\mu_*\in M_1^2(X)$ (lemat \ref{second_moment}), otrzymujemy, że lewa strona wyrażenia $(\ref{condition_woodroofe})$ jest nie większa niż
\begin{align*}
\begin{aligned}
\sum_{n=1}^{\infty}n^{-3/2}\left(C_9^2\:\left\langle 1+2V+V^2,\mu_*\right\rangle\: \right)^{1/2}<\infty. 
\end{aligned}
\end{align*} 
Otrzymaliśmy tezę.
\end{proof}

Spełnione są założenia twierdzenia \ref{tw:wn_woodr}, zatem $\sigma^2=\lim_{n\to\infty}E_{\mu_*}\left(\left(S_n^*\right)^2\right)$ oraz $\Phi_{S_n^*}$ zbiega słabo do rozkładu normalnego $\mathcal{N}\left(0,\sigma^2\right)$, gdy $n\to\infty$. 
CTG jest więc spełnione dla łańcuchów Markowa o funkcji przejścia $\Pi_{\varepsilon}$ i~rozkładzie początkowym zadanym przez miarę niezmienniczą~$\mu_*$.

\section{Centralne twierdzenie graniczne dla uogólnionego modelu cyklu komórkowego}
\begin{twierdzenie}\label{CTG}
Dana jest przestrzeń polska $(X,\varrho)$. Niech $\mu\in M_1^2(X)$. Przyjmijmy, że $\Phi_{S_n^{\mu}}$ jest rozkładem zmiennej losowej $S_n^{\mu}$ określonej wzorem (\ref{def:eta_n^mu}). Założenia (I)-(VI) oraz (II') implikują istnienie skończonej granicy $\sigma^2=\lim_{n\to\infty}E_{\mu_*}\left(\left(S_n^*\right)^2\right)$ oraz słabą zbieżność ciągu rozkładów $(\Phi_{S_n^{\mu}})_{n\in N}$ do rozkładu normalnego $\mathcal{N}\left(0,\sigma^2\right)$, przy $n\to\infty$.
\end{twierdzenie}

\begin{proof}
Niech $x,y\in X$ oraz $f\in\mathcal{L}$ będą dowolne. Ponieważ w podrozdziale 6.1 ustaliliśmy już, że łańcuchy Markowa o funkcji przejścia $\Pi_{\varepsilon}$, startujące z miary niezmienniczej $\mu_*\in M_1^2(X)$ spełniają CTG, to wystarczy pokazać, że $\left|\Phi_{S_n^{\mu}}-\Phi_{S_n^*}\right|\xrightarrow{w} 0$, gdy $n\to\infty$. Równoważnie (stwierdzenie~\ref{tw_aleksandrowa}), wystarczy dowieść, że 
\begin{align}\label{prop:xy}
\begin{aligned}
\lim_{n\to\infty}\left|\left\langle f,\Phi_{S_n^x}-\Phi_{S_n^y}\right\rangle\right|=0\quad\text{dla }f\in\mathcal{L}.
\end{aligned}
\end{align}
Wówczas z twierdzenia Lebesgue'a o zbieżności ograniczonej otrzymamy
\begin{align*}
\begin{aligned}
\lim_{n\to\infty}\left|\left\langle f,\Phi_{S_n^{\mu}}\right\rangle-\left\langle f,\Phi_{S_n^*}\right\rangle\right|
\leq\lim_{n\to\infty}\int_X\int_X\left|\left\langle f,\Phi_{S_n^x}\right\rangle-\left\langle f,\Phi_{S_n^y}\right\rangle\right|\mu(dx)\mu_*(dy)=0
\end{aligned}
\end{align*}
 i dowód będzie ukończony. Pozostaje wykazać zbieżność $(\ref{prop:xy})$. Zauważmy, że
\begin{align}\label{prop:pre-estimate}
\begin{aligned}
\left|\left\langle f,\Phi_{S_n^x}\right\rangle-\left\langle f,\Phi_{S_n^y}\right\rangle\right|
=&\Bigg|\int_{X^n}f\left(\frac{g(u_1)+\ldots+g(u_{n})}{\sqrt{n}}\right)\Pi_{\varepsilon}^{1,\ldots,n}\left(x,du_1\times\ldots\times du_n\right)\,-\\
& - \int_{X^n}f\left(\frac{g(v_1)+\ldots+g(v_{n})}{\sqrt{n}}\right)\Pi_{\varepsilon}^{1,\ldots,n}\left(y,dv_1\times\ldots\times dv_n\right)  \Bigg|,
\end{aligned}
\end{align}
gdzie $\Pi_{\varepsilon}^{1,\ldots,n}(z,\cdot)=
\int_{\bar{B}(0,\varepsilon)}\ldots\int_{\bar{B}(0,\varepsilon)}\Pi_{h_1,\ldots,h_n}^{1,\ldots,n}(z,\cdot)\nu^{\varepsilon}(dh_1)\ldots\nu^{\varepsilon}(dh_n)$ dla $z\in X$ jest miarą na $X^n$. Dalej mamy
\begin{align}\label{estimate}
\begin{aligned}
&\left|\int_{X^{2n}}
\left(f\left(\frac{g(u_1)+\ldots+g(u_n)}{\sqrt{n}}\right) 
- f\left(\frac{g(v_1)+\ldots+g(v_n)}{\sqrt{n}}\right)\right)
\Pi_{h_1,\ldots,h_n}^{1,\ldots,n}\left(x,\,\times_{i=1}^ndu_i\right)
\Pi_{h_1,\ldots,h_n}^{1,\ldots,n}\left(y,\times_{i=1}^ndv_i\right)  \right|\leq\\
&\leq\int_{X^{2n}}
\left| f\left(\frac{g(u_1)+\ldots+g(u_{n})}{\sqrt{n}}\right) 
- f\left(\frac{g(v_1)+\ldots+g(v_{n})}{\sqrt{n}}\right) \right|\,\times\\
&\quad\times\,\left( \Pi^*_{X^{2n}}\Pi^*_{1,\ldots,n}\hat{C}^{\infty}_{h_1,h_2\ldots}((x,y,1),\cdot)\right)\left(\times_{i=1}^n \left(du_i\times dv_i\right)\right),
\end{aligned}
\end{align}
gdzie, dla dowolnego $n\in N$, odwzorowanie $\Pi^*_{1,\ldots,n}:\left(X^2\times\left\{0,1\right\}\right)^{\infty}\to\left(X^2\times\left\{0,1\right\}\right)^n$ jest rzutem na pierwsze $n$ współrzędnych oraz $\Pi^*_{X^{2n}}:\left(X^2\times\left\{0,1\right\}\right)^n\to X^{2n}$ jest rzutem na przestrzeń $X^{2n}$. 
Ponieważ funkcja $f$ spełnia warunek Lipschitza ze stałą $1$, możemy oszacować $(\ref{estimate})$ z~góry w~następujący sposób
\begin{align*}
\begin{aligned}
&\frac{1}{\sqrt{n}}\int_{X^{2n}}\left(\left|g(u_1)-g(v_1)\right|+\ldots+\left|g(u_n)-g(v_n)\right|\right)
\left(\Pi^*_{X^{2n}}\Pi^*_{1,\ldots,n}\hat{C}^{\infty}_{h_1,h_2\ldots}((x,y,1),\cdot)\right)\left(\times_{i=1}^n(du_i\times dv_i)\right)=\\
&= \frac{1}{\sqrt{n}}\sum_{i=1}^n\int_{X^2}\left|g(u_i)-g(v_i)\right|
\left(\Pi^*_{X^{2}}\Pi^*_{i}\hat{C}^{\infty}_{h_1,h_2\ldots}((x,y,1),\cdot)\right)
\left(du_i\times dv_i\right)\leq\\
&\leq \frac{G}{\sqrt{n}}\sum_{i=1}^nq^iC_5\left(1+\bar{V}(x,y)\right)=\\
&=n^{-\frac{1}{2}}GC_5q\frac{1-q^n}{1-q}\left(1+\bar{V}(x,y)\right),
\end{aligned}
\end{align*}
gdzie powyższa nierówność wynika z lematu \ref{theorem_new} oraz wniosku \ref{wniosek_g}. 
Ponieważ wyrażenie jest niezależne od $(h_n)_{n\in N}$, to jest również dobrym oszacowaniem górnym dla $(\ref{prop:pre-estimate})$. Przechodząc z~$n$~do nieskończoności, otrzymujemy $(\ref{prop:xy})$. Dowód został ukończony.
\end{proof}

\chapter{Prawo iterowanego logarytmu}

W ostatnim rozdziale, adaptując metodę martyngałową (por. \cite{heyde_scott}, \cite{bms}, \cite{kom_szar}) oraz korzystając z~odpowiednich własności miary sprzęgającej, udowodnimy PIL dla uogólnionego modelu cyklu komórkowego.

\section{Metoda martyngałowa}
Niech proces $(M_n)_{n\in N_0}$ określony na $(\Omega,\mathcal{F},\mathbb{P})$ będzie martyngałem względem swojej filtracji naturalnej $(\mathcal{F}_n)_{n\in N_0}$. Określmy proces $(Z_n)_{n\in N_0}$ w następujący sposób: $Z_0=M_0=0\;$ $\mathbb{P}$-p.n. oraz $Z_n=M_n-M_{n-1}$ dla $n\in N$. Ponadto zdefiniujmy $s_n^2=EM_n^2<\infty$.

Rozważmy przestrzeń $(\mathcal{C},\tilde{\varrho})$ wszystkich funkcji $f:[0,1]\to R$ z metryką $\tilde{\varrho}$ postaci
\begin{align*}
\begin{aligned}
\tilde{\varrho}\left(f_1,f_2\right)=\sup_{t\in[0,1]}\left|f_1(t)-f_2(t)\right|
\quad \text{dla }f_1,f_2\in \mathcal{C}.
\end{aligned}
\end{align*}
Zdefiniujmy również zbiór $\mathcal{K}$ funkcji $f\in \mathcal{C}$ absolutnie ciągłych i~takich, że~ $f(0)=0$ oraz $\int_0^1\left(f'(t)\right)^2dt\leq 1$. Niech funkcja $F:[0,\infty)\to N_0$ będzie postaci $F(s)=\sup\left\{n\in N_0:s_n^2\leq s\right\}$. Ponadto niech ciąg funkcji losowych $(\eta_n)_{n\in N_0}$, $\eta_n:[0,1]\to R$, będzie dany wzorem 
\begin{align*}
\begin{aligned}
\eta_n(t)=\frac{M_k+\left(s_n^2t-s_k^2\right)\left(s_{k+1}^2-s_k^2\right)^{-1}Z_{k+1}}{\sqrt{2s_n^2\log\log s_n^2}}\quad\text{dla }n>F(e),
\end{aligned}
\end{align*}
gdzie $1\leq k\leq n-1$, $\;s_k^2\leq s_n^2 t\leq s_{k+1}^2$, oraz niech $\eta_n(t)=0$ dla $n\leq F(e)$.

Rodzinę $\mathcal{G}\subset \mathcal{C}$ nazywamy relatywnie zwartą, jeśli domknięcie zbioru $\left\{f:[0,1]\to R \;:\;f\in\mathcal{G}\right\}$ jest zbiorem zwartym w $(\mathcal{C},\tilde{\varrho})$.

\begin{twierdzenie}\label{tw:heyde-scott}
Niech $(M_n)_{n\in N_0 }$ będzie martyngałem całkowalnym z kwadratem oraz niech $s_n^2=EM_n^2$ dla $n\in N_0$. Jeśli $\lim_{n\to\infty}s_n^2=\infty$ oraz spełnione są warunki:
\begin{align}\label{th1}
\begin{aligned}
&\sum_{n=1}^{\infty}s_n^{-4}E\left(Z_n^41_{\left\{|Z_n|<\beta s_n\right\}}\right)< \infty\quad\text{dla pewnej wartości }\beta>0,
\end{aligned}
\end{align}
\begin{align}\label{th2}
\begin{aligned}
\sum_{n=1}^{\infty}s_n^{-1}E\left(|Z_n|1_{\left\{|Z_n|\geq\vartheta s_n\right\}}\right)<\infty\quad\text{dla wszystkich }\vartheta>0,
\end{aligned}
\end{align}
\begin{align}\label{th3}
\begin{aligned}
s_n^{-2}\sum_{k=1}^nZ_k^2\to 1\quad \mathbb{P}\text{-p.n., gdy }\;n\to\infty,
\end{aligned}
\end{align}
to ciąg $(\eta_n)_{n> F(e)}$ jest relatywnie zwarty w $\mathcal{C}$, a zbiór jego punktów granicznych pokrywa się ze~zbiorem $\mathcal{K}$ (twierdzenie 1, \cite{heyde_scott}).
\end{twierdzenie}

\section{Zastosowanie do badanego modelu}
Ustalmy $\varepsilon\in[0,\varepsilon_*]$ dla $\varepsilon_*<\infty$. 
Niech $\mu_*$ będzie słabą granicą ciągu $(P^n_{\varepsilon}\mu)_{n\in N_0}$ oraz jedyną miarą niezmienniczą operatora $P_{\varepsilon}$, której istnienie wynika z twierdzenia \ref{GTZ}. Weźmy łańcuchy Markowa $(x_n)_{n\in N_0},(y_n)_{n\in N_0}$ o~wartościach w~przestrzeni stanów $X$, funkcji przejścia $\Pi_{\varepsilon}$ oraz miarach początkowych $\mu,\nu\in M_1^{2+\delta}(X)$ dla $\delta>0$. Załóżmy, że $g\in B(X)$ jest funkcją Lipschitza ze~stałą $L_g>0$. Ponadto niech $\left\langle g,\mu_*\right\rangle =0$ (w przeciwnym przypadku kładziemy $\tilde{g}=g-\left\langle g,\mu_*\right\rangle $).

Z nierówności Minkowskiego w $L^{2+\delta}\left(P^n_{\varepsilon}\mu\right)$ dla $n\in N_0$ oraz z wniosku \ref{wn:skończone_momenty_P_varepsilon_mu} otrzymujemy
\begin{align}\label{minkowski}
\begin{aligned}
\left(E_{\mu}\left(|g(x_n)|^{2+\delta}\right)\right)^{1/(2+\delta)}&=\left(\int_X|g(x)|^{2+\delta}P^n_{\varepsilon}\mu(dx)\right)^{1/(2+\delta)}\leq\\
&\leq|g(\bar{x})|+L_g\left(\left\langle V^{2+\delta},P^n_{\varepsilon}\mu\right\rangle\right)^{1/(2+\delta)}\leq\\
&\leq |g(\bar{x})|+L_g\left( \left(\left\langle V^{2+\delta},\mu\right\rangle\right)^{1/(2+\delta)}+{c}{\left(1-\Lambda_{2+\delta}^{1/(2+\delta)}\right)^{-1}}\right)<\infty.
\end{aligned}
\end{align}
Zatem $\sup_{n\in N_0}E_{\mu}\left(|g(x_n)|^{2+\delta}\right)<\infty$. Niech $x\in X$. Zauważmy, że założenie $\left\langle g,\mu_*\right\rangle=0$ oraz równość  (\ref{P_epsilon-as-P_h}) implikują
\begin{align}\label{finite_chi}
\begin{aligned}
\sum_{i=0}^{\infty}\left|U^i_{\varepsilon}g(x)\right|
&=\sum_{i=0}^{\infty}\left|\left\langle g,P^i_{\varepsilon}\delta_x\right\rangle 
-\left\langle g,P^i_{\varepsilon}\mu_*\right\rangle \right|\leq\\
&\leq\sum_{i=0}^{\infty}\int_{\left(\bar{B}(0,\varepsilon)\right)^i}
\left|\left\langle g,P^i_{h_1,\ldots,h_i}\delta_x\right\rangle 
-\left\langle g,P^i_{h_1,\ldots,h_i}\mu_*\right\rangle \right|\nu^{\varepsilon}
(dh_1)\ldots\nu^{\varepsilon}(dh_i).
\end{aligned}
\end{align}
Z wniosku \ref{wniosek_g} otrzymujemy 
\begin{align}\label{finite_chi_cd}
\begin{aligned}
&\left|\left\langle g,P^i_{h_1,\ldots,h_i}\delta_x\right\rangle 
-\left\langle g,P^i_{h_1,\ldots,h_i}\mu_*\right\rangle \right|\leq\\
&\leq\int_X\int_{X^2}|g(u)-g(v)|
\left(\Pi^*_{X^2}\Pi^*_i\hat{C}^{\infty}_{h_1,h_2,\ldots}((x,y,1),\cdot)\right)
(du\times dv)\;\mu_*(dy)\leq\\
&\leq q^iGC_5\left(1+V(x)+\left\langle V,\mu_*\right\rangle \right),
\end{aligned}
\end{align}
gdzie $G:=\max\left\{L_g,\sup_{x\in X}|g(x)|\right\}$. 
Wobec (\ref{finite_chi}) i (\ref{finite_chi_cd}) dostajemy 
\begin{align}\label{fact:finite_chi}
\begin{aligned}
\sum_{i=0}^{\infty}\left|U^i_{\varepsilon}g(x)\right|
\leq (1-q)^{-1}GC_5\left(1+V(x)+\left\langle V,\mu_*\right\rangle \right)<\infty,
\end{aligned}
\end{align}
dzięki czemu możemy definiować funkcję
\begin{align}\label{def:chi}
\begin{aligned}
\chi(x)=\sum_{i=0}^{\infty}U^i_{\varepsilon}g(x)\quad \text{dla }x\in X.
\end{aligned}
\end{align}

\begin{lemat}\label{lem:chi-chi}
Niech funkcja $\chi$ będzie dana wzorem (\ref{def:chi}). Wówczas
\begin{align*}
\begin{aligned}
|\chi(x)-\chi(y)|\leq \frac{GC_5}{1-q}(1+V(x)+V(y))\quad\text{dla }x,y\in X.
\end{aligned}
\end{align*}
\end{lemat}
\begin{proof}
Ustalmy $x,y\in X$. 
Stosując $(\ref{P_epsilon-as-P_h})$ oraz wniosek \ref{wniosek_g}, otrzymujemy
\begin{align*}
\begin{aligned}
&|\chi(x)-\chi(y)|=\\
&=\left|\sum_{i=0}^{\infty}U^i_{\varepsilon}g(x)
-\sum_{i=0}^{\infty}U^i_{\varepsilon}g(y)\right|\leq\\
&\leq\sum_{i=0}^{\infty}\left|\left\langle g,P^i_{\varepsilon}\delta_x\right\rangle 
-\left\langle g,P^i_{\varepsilon}\delta_y\right\rangle \right|\leq\\
&\leq\sum_{i=0}^{\infty}\int_{\left(\bar{B}(0,\varepsilon)\right)^i}\int_{X^2}|g(u)-g(v)|
\left(\Pi^*_{X^2}\Pi^*_i\hat{C}^{\infty}_{h_1,h_2,\ldots}((x,y,1),\cdot)\right)
(du\times dv)\;\nu^{\varepsilon}(dh_1)\ldots\nu^{\varepsilon}(dh_i)\leq\\
&\leq \sum_{i=0}^{\infty}q^iGC_5(1+V(x)+V(y))
=(1-q)^{-1}GC_5(1+V(x)+V(y)).
\end{aligned}
\end{align*}
Dowód został ukończony.
\end{proof}

Wprowadźmy zmienne losowe
\begin{align}\label{def:M_n}
\begin{aligned}
M_n=\chi(x_n)-\chi(x_0)+\sum_{i=0}^{n-1}g(x_i)\quad \text{dla }n\in N_0
\end{aligned}
\end{align}
oraz 
\begin{align}\label{def:Z_n}
\begin{aligned}
Z_n=\chi(x_n)-\chi(x_{n-1})+g(x_{n-1})\quad\text{dla }n\in N,\quad Z_0=0.
\end{aligned}
\end{align}

\begin{lemat}
Proces $(M_n)_{n\in N_0}$ postaci (\ref{def:M_n}) jest martyngałem w przestrzeni $(X^{\infty},\otimes_{i=1}^{\infty}B_X, \mathbb{P}_\mu)$ względem swojej filtracji naturalnej $(\mathcal{F}_n)_{n\in N_0}$.
\end{lemat}

\begin{proof}
Własność Markowa implikuje 
\begin{align}\label{Markow_pil}
E_{\mu}\left(g\circ \Gamma_n|\mathcal{F}_n\right)(\omega)
=E_{x_n(\omega)}\left(g\right)=(U_{\varepsilon}g)(x_n(\omega)),
\end{align}
zatem
\begin{align*}
\begin{aligned}
E_{\mu}\left(M_{n+1}|\mathcal{F}_n\right)&=E_{\mu}\left(\chi(x_{n+1})|\mathcal{F}_n\right)-\chi(x_0)+\sum_{i=0}^ng(x_i)=\\
&=\sum_{i=0}^{\infty}U_{\varepsilon}(U^i_{\varepsilon}g)(x_n)-\chi(x_0)+\sum_{i=0}^ng(x_i)=\\
&=\sum_{i=1}^{\infty}U^i_{\varepsilon}g(x_n)+U_{\varepsilon}^0g(x_n)-\chi(x_0)+\sum_{i=0}^{n-1}g(x_i)=M_n,
\end{aligned}
\end{align*}
co kończy dowód.
\end{proof}

\begin{lemat}\label{prop:Z_n}
Niech $(Z_n)_{n\in N_0}$ będzie postaci (\ref{def:Z_n}). Wówczas $E_{\mu_*}Z_1^2<\infty$.
\end{lemat}

\begin{proof} 
Ustalmy dowolne $n\in N_0$. Niech $\mu\in M_1^{2+\delta}(X)$. Korzystając z własności Markowa (por.~(\ref{Markow_pil})), otrzymujemy
\begin{align}\label{estim:Z}
\begin{aligned}
E_{\mu}\left(Z_{n+1}^2\right)=E_{P^n_{\varepsilon}\mu}\left(Z_1^2\right)&=\int_XE\left(\left(\chi(x_1)-\chi(x_0)+g(x_0)\right)^2\,|\,x_0=x\right)P_{\varepsilon}^n\mu(dx)\leq\\
&\leq\int_XE\left(2\chi^2(x_1)+2(\chi-g)^2(x_0)\,|\,x_0=x\right)\:P_{\varepsilon}^n\mu(dx)=\\
&= \int_X\left(2E(\chi^2(x_1)|x_0=x)+2(\chi-g)^2(x)\right)\:P_{\varepsilon}^n\mu(dx)\leq\\
&\leq 2\int_X\left(U_{\varepsilon}\chi^2\right)(x)P^n_{\varepsilon}\mu(dx)
+4\int_X\chi^2(x)P^n_{\varepsilon}\mu(dx)+4\int_Xg^2(x)P^n_{\varepsilon}\mu(dx)=\\
&=2\int_X\chi^2(x)P^{n+1}_{\varepsilon}\mu(dx)+4\int_X\chi^2(x)P^n_{\varepsilon}\mu(dx)+4E_{\mu}\left(g(x_n)\right)^2.
\end{aligned}
\end{align}
Postępując analogicznie jak w $(\ref{minkowski})$, otrzymujemy, że ostatni składnik w (\ref{estim:Z}) jest skończony. Wystarczy więc pokazać, że wyrażenie $\left\langle \chi^2,P^n_{\varepsilon}\mu\right\rangle $ jest skończone dla $n\in N_0$. Zauważmy, że
\begin{align}\label{estim:chi}
\begin{aligned}
\int_X\chi^2(x)P^n_{\varepsilon}\mu(dx)
&=\int_X\left(\left(\chi(x)-\chi(\bar{x})\right)+\chi(\bar{x})\right)^2
P^n_{\varepsilon}\mu(dx)\leq\\
&\leq 2\chi^2\left(\bar{x}\right)+2\int_X\left(\chi(x)-\chi(\bar{x})\right)^2
P^n_{\varepsilon}\mu(dx).
\end{aligned}
\end{align}
Na mocy (\ref{fact:finite_chi}) i (\ref{def:chi}) pierwszy składnik w oszacowaniu (\ref{estim:chi}) jest skończony. Skończoność drugiego składnika wynika z~lematu \ref{lem:chi-chi}. Istotnie,
\begin{align}\label{estim:chi2}
\begin{aligned}
2\int_X\left(\chi(x)-\chi(\bar{x})\right)^2P^n_{\varepsilon}\mu(dx)
&\leq 2\int_X(1-q)^{-2}G^2C_5^2(1+V(x))^2P^n_{\varepsilon}\mu(dx)\leq\\
&\leq 4G^2C_5^2(1-q)^{-2}\left(1+\left\langle V^2,P^n_{\varepsilon}\mu\right\rangle \right)<\infty.
\end{aligned}
\end{align}
Podsumowując, z (\ref{minkowski}) i (\ref{estim:Z})-(\ref{estim:chi2}) otrzymujemy
\begin{align*}
\begin{aligned}
E_{P^n_{\varepsilon}\mu}\left(Z_1^2\right)\leq\tilde{C}\left(1+\left\langle V^2,P^{n+1}_{\varepsilon}\mu\right\rangle +\left\langle V^2,P^n_{\varepsilon}\mu\right\rangle \right)
\end{aligned}
\end{align*}
i na mocy wniosku \ref{wn:skończone_momenty_P_varepsilon_mu} mamy
\begin{align}\label{supEZn-finite}
\begin{aligned}
\sup_{n\in N_0}E_{P^n_{\varepsilon}\mu}\left(Z_1^2\right)
\leq\bar{C}\left(1+\left\langle V^2,\mu\right\rangle\right).
\end{aligned}
\end{align}
Łatwo zauważyć, że $x\mapsto E_x\left(Z_1^2\wedge k\right)$ jest odwzorowaniem ciągłym i ograniczonym dla dowolnych $k\in N_0$, zatem z twierdzenia \ref{GTZ} mamy 
\[\lim_{n\to\infty}E_{P^n_{\varepsilon}\mu}(Z_1^2\wedge k)=E_{\mu_*}\left(Z_1^2\wedge k\right).\]
Wobec powyższego, a także wobec oszacowania (\ref{supEZn-finite}) ciąg $\left(E_{\mu_*}\left(Z_1^2\wedge k\right)\right)_{k\in N_0}$ jest wspólnie ograniczony. Stosując twierdzenie Lebesgue'a o~zbieżności monotonicznej, otrzymujemy 
\begin{align*}
\lim_{k\to\infty}E_{\mu_*}\left(Z_1^2\wedge k\right)=E_{\mu_*}\left(Z_1^2\right)<\infty,
\end{align*}
co kończy dowód.
\end{proof}

Ustalmy
\begin{align}\label{def:sigma2}
\begin{aligned}
\sigma^2=E_{\mu_*}Z_1^2.
\end{aligned}
\end{align}

\begin{lemat}\label{lem:sn/n=sigma}
Weźmy $\mu\in M_1^{2+\delta}(X)$. Niech zmienne losowe $M_n$ będą postaci $(\ref{def:M_n})$ oraz niech $s^2_n=E_{\mu}M_n^2<\infty$ dla $n\in N_0$. Wówczas
\[\lim_{n\to\infty}\frac{s_n^2}{n}=\sigma^2.\]
\end{lemat}

\begin{proof}
Postępując analogicznie jak w dowodzie lematu \ref{prop:Z_n} (zob. (\ref{estim:Z})-(\ref{supEZn-finite})), otrzymujemy
\begin{align}\label{supEZn-delta-finite}
\begin{aligned}
\sup_{n\in N_0}E_{\mu}\left|Z_n\right|^{2+\delta}<\infty.
\end{aligned}
\end{align}
Mamy zatem 
\begin{align*}
\begin{aligned}
\sup_{n\in N_0}E_{\mu}\left(Z_n^21_{\left\{|Z_n|^2\geq k\right\}}\right)
\leq\sup_{n\in N_0}E_{\mu}\left(|Z_n|^{2+\delta}\left(|Z_n|^2\right)^{-\delta/2} 1_{\left\{|Z_n|^2\geq k\right\}}\right)
\leq k^{-\delta/2}\sup_{n\in N_0}E_{\mu}|Z_n|^{2+\delta}\to 0,
\end{aligned}
\end{align*}
gdy $k\to\infty$. Ponieważ zmienne $\left(Z_1^2\wedge k\right)$ są ciągłe i ograniczone oraz $P_{\varepsilon}$ ma własność Fellera, co wynika z wniosku \ref{wn:P_epsilon-Fellera}, to z~twierdzenia~\ref{GTZ} otrzymujemy 
\[\lim_{n\to\infty}E_{P^n_{\varepsilon}\mu}\left(Z_1^2\wedge k\right)
= E_{\mu_*}\left(Z_1^2\wedge k\right) 
\quad\text{dla }k\in N_0.\]
Zauważmy, że ciąg $\left(E_{\mu_*}\left(Z_1^2\wedge k\right)\right)_{k\in N_0}$ jest wspólnie ograniczony (por. (\ref{supEZn-finite})) i~dlatego z~twierdzenia Lebesgue'a  o~zbieżności monotonicznej mamy
\[\lim_{k\to\infty}E_{\mu_*}\left(Z_1^2\wedge k\right)=E_{\mu_*}Z_1^2=\sigma^2,\]
wobec czego
\begin{align*}
\lim_{n\to\infty}E_{\mu}Z_{n+1}^2=\lim_{n\to\infty}E_{P^n_{\varepsilon}\mu}Z_1^2
=E_{\mu_*}Z_1^2=\sigma^2.
\end{align*}
Ostatecznie z ortogonalności różnic martyngałowych otrzymujemy
\begin{align*}
\begin{aligned}
\lim_{n\to\infty}\frac{s_n^2}{n}=\lim_{n\to\infty}\frac{E_{\mu}M_n^2}{n}=\lim_{n\to\infty}\frac{\sum_{i=1}^nE_{\mu}Z_i^2}{n}=\sigma^2,
\end{aligned}
\end{align*}
co kończy dowód.
\end{proof}

\begin{stwierdzenie}
Wariancja $\sigma^2=E_{\mu_*}Z_1^2$ jest zgodna z wariancją graniczną rozkładu normalnego z~twierdzenia \ref{CTG}, tzn. $\sigma^2=\lim_{n\to\infty}E_{\mu_*}\left(\left(S_n^*\right)^2\right)$.
\end{stwierdzenie}
\begin{proof}
Zauważmy, że 
\begin{align}\label{1}
\begin{aligned}
\lim_{n\to\infty} E_{\mu_*}\left(\left(S_n^*\right)^2\right)
&=\lim_{n\to\infty} E_{\mu_*}\left(n^{-1}\left(\sum_{i=0}^{n-1}g(x_i)\right)^2\right)=\\
&=\lim_{n\to\infty} E_{\mu_*}\left(n^{-1}\left(M_n+\chi(x_0)-\chi(x_n)\right)^2\right)=\\
&=\lim_{n\to\infty} E_{\mu_*}\left(n^{-1}M_n^2\right)
+\lim_{n\to\infty} 2n^{-1/2}E_{\mu_*}\left(\left(n^{-1/2}M_n\right)
\left(\chi(x_0)-\chi(x_n)\right)\right)+\\
&+\lim_{n\to\infty} n^{-1}E_{\mu_*}\left(\chi(x_0)-\chi(x_n)\right)^2.
\end{aligned}
\end{align}
Na mocy lematu \ref{lem:sn/n=sigma} mamy $\lim_{n\to\infty} E_{\mu_*}\left(n^{-1}M_n^2\right)=E_{\mu_*}Z_1^2$. Wobec lematu \ref{lem:chi-chi} otrzymujemy
\begin{align}\label{2}
\begin{aligned}
E_{\mu_*}\left(\chi(x_0)-\chi(x_n)\right)^2
&=\int_X\int_X\left(\chi(u)-\chi(v)\right)^2P^n_{\varepsilon}\delta_u(dv)\mu_*(du)\leq\\
&\leq\int_X\int_XG^2C_5^2(1-q)^{-2}(1+V(u)+V(v))^2P^n_{\varepsilon}\delta_u(dv)\mu_*(du)\leq\\
&\leq C_0\int_X\int_X(1+V^2(u)+V^2(v))P^n_{\varepsilon}\delta_u(dv)\mu_*(du)\leq\\
&\leq C_0\int_X \left(1+V^2(u)+\left\langle V^2,P^n_{\varepsilon}\delta_u\right\rangle\right)\mu_*(du),
\end{aligned}
\end{align}
gdzie $C_0$ jest pewną stałą. Zgodnie z oszacowaniem (\ref{oszacowanie_mom_skończ_epsilon}) mamy
\begin{align}\label{3}
E_{\mu_*}\left(\chi(x_0)-\chi(x_n)\right)^2
\leq \tilde{C}_0\int_X \left(1+V^2(u)\right)\mu_*(du)<\infty,
\end{align}
gdzie $\tilde{C}_0$ jest pewną stałą. Stosując nierówność H\"{o}ldera, otrzymujemy
\begin{align*}
E_{\mu_*}\left|\left(n^{-1/2}M_n\right)
\left(\chi(x_0)-\chi(x_n)\right)\right|
\leq\left(E_{\mu_*}\left(n^{-1}M_n^2\right)\right)^{1/2}
\left(E_{\mu_*}\left( \chi(x_0)-\chi(x_n)\right)^2\right)^{1/2}<\infty,
\end{align*}
zatem
\begin{align}\label{4}
\lim_{n\to\infty} 2n^{-1/2}E_{\mu_*}\left(\left(n^{-1/2}M_n\right)
\left(\chi(x_0)-\chi(x_n)\right)\right)=0.
\end{align}
Podsumowując powyższe oszacowania (\ref{1})-(\ref{4}), otrzymujemy
\begin{align*}
\lim_{n\to\infty} E_{\mu_*}\left(\left(S_n^*\right)^2\right)=E_{\mu_*}Z_1^2,
\end{align*}
co kończy dowód.
\end{proof}

\begin{lemat}\label{lem1}
Niech $(Z_n)_{n\in N}$ będzie postaci (\ref{def:Z_n}). Wówczas mamy
\begin{align}\label{lem_sigma2}
\begin{aligned}
\frac{1}{n}\sum_{l=1}^nZ_l^2\to\sigma^2 \quad \mathbb{P}_{\mu}\text{-p.n., $\:$gdy }n\to\infty,
\end{aligned}
\end{align}
co implikuje warunek $(\ref{th3})$ dla $\sigma^2>0$.
\end{lemat}
\begin{proof}
Odwzorowania
\begin{align}\label{functions}
\begin{aligned}
x\mapsto E_x\left(\left|{\lim\inf}_{n\to\infty}
\left(\frac{1}{n}\sum_{l=1}^nZ_l^2\right)-\sigma^2\right|\wedge 1\right)\\
x\mapsto E_x\left(\left|{\lim\sup}_{n\to\infty}\left(\frac{1}{n}\sum_{l=1}^nZ_l^2\right)-\sigma^2\right|\wedge 1\right)
\end{aligned}
\end{align}
są w oczywisty sposób ograniczone. Istotne jest pokazanie ich ciągłości. Przypuścimy, że są ciągłe. Wówczas mamy
\begin{align}\label{convergence}
\begin{aligned}
&E_{\mu}\left(\left|{\lim\inf}_{n\to\infty}
\left(\frac{1}{n}\sum_{l=1}^nZ_l^2\right)-\sigma^2\right|\wedge 1\right)=\\
&=\int_X E_x\left(\left|{\lim\inf}_{n\to\infty}
\left(\frac{1}{n}\sum_{l=1}^nZ_l^2\right)-\sigma^2\right|\wedge 1\right)\mu(dx)=\\
&=\int_XE_x\left(\left|{\lim\inf}_{n\to\infty}
\left(\frac{1}{n}\sum_{l=1}^nZ_l^2\right)-\sigma^2\right|\wedge 1\right)P^m_{\varepsilon}\mu(dx)
\to
E_{\mu_*}\left(\left|{\lim\inf}_{n\to\infty}
\left(\frac{1}{n}\sum_{l=1}^nZ_l^2\right)-\sigma^2\right|\wedge 1\right),
\end{aligned}
\end{align}
gdy $m\to\infty$. 
Miara $\mu_*$ jest jedyna, więc ergodyczna. Wobec tego na mocy twierdzenia ergodycznego Birkhoffa, które mówi, że
\begin{align*}
\begin{aligned}
E_{\mu_*}\left(\left|{\lim\inf}_{n\to\infty}\left(\frac{1}{n}\sum_{l=1}^nZ_l^2\right)-\sigma^2\right|\wedge 1\right)=0,
\end{aligned}
\end{align*}
otrzymujemy
\begin{align*}
\begin{aligned}
E_{\mu}\left(\left|{\lim\inf}_{n\to\infty}\left(\frac{1}{n}\sum_{l=1}^nZ_l^2\right)-\sigma^2\right|\wedge 1\right)=0,
\end{aligned}
\end{align*}
co z kolei implikuje
\begin{align}\label{estim:inf}
\begin{aligned}
{\lim\inf}_{n\to\infty}\left(\frac{1}{n}\sum_{l=1}^nZ_l^2\right)=\sigma^2\quad \mathbb{P}_{\mu}\text{-p.n.}
\end{aligned}
\end{align}
W sposób analogiczny otrzymujemy 
\begin{align}\label{estim:sup}
\begin{aligned}
{\lim\sup}_{n\to\infty}\left(\frac{1}{n}\sum_{l=1}^nZ_l^2\right)=\sigma^2\quad \mathbb{P}_{\mu}\text{-p.n.}
\end{aligned}
\end{align}
Ostatecznie z (\ref{estim:inf}) i (\ref{estim:sup}) wynika (\ref{lem_sigma2}).

Aby ukończyć dowód, musimy pokazać, że funkcje postaci $(\ref{functions})$ są ciągłe, co jest niezbędne, by zachodziła zbieżność w $(\ref{convergence})$. Zauważmy, że
\begin{align}\label{min_Hnk}
\begin{aligned}
E_x\left(\left|{\lim\inf}_{n\to\infty}
\left(\frac{1}{n}\sum_{l=1}^nZ_l^2\right)-\sigma^2\right|\wedge 1\right)
&=\lim_{n\to\infty}\lim_{k\to\infty}
E_x\left(\left|\min_{j\in\{n,\ldots,n+k\}}
\left(\frac{1}{j}\sum_{l=1}^jZ_l^2-\sigma^2\right)\right|\wedge 1\right)=\\
&=\lim_{n\to\infty}\lim_{k\to\infty}H_{n,k}(x),
\end{aligned}
\end{align}
gdzie
\begin{align}
\begin{aligned}
H_{n,k}(x):=E_x\left(\left|\min_{j\in\{n,\ldots,n+k\}}
\left(\frac{1}{j}\sum_{l=1}^jZ_l^2-\sigma^2\right)\right|\wedge 1\right).
\end{aligned}
\end{align}
Niech dane będzie odwzorowanie 
\begin{align}\label{def:g}
\begin{aligned}
\psi_{n,k}(y_0,\ldots,y_{n+k})&=
\left| \min_{j\in\{n\ldots,n+k\}} 
\left( \frac{1}{j}
\left(\sum_{l=1}^{j}\left(\chi(y_l)-\chi(y_{l-1})+g(y_{l-1})\right)^2
\wedge j\left(1+\sigma^2\right)\right)
-\sigma^2 \right) \right|.
\end{aligned}
\end{align}
Z definicji różnic martyngałowych $(Z_n)_{n\in N}$ (zob. (\ref{def:Z_n})) oraz z własności
\begin{align*}
\begin{aligned}
 \left|\min_{j\in\{n,\ldots,n+k\}}\left(\frac{a_j}{j}-b\right)\right|\wedge 1
=\left|\min_{j\in\{n,\ldots,n+k\}}\left(\frac{1}{j}(a_j\wedge j(1+b))-b\right)\right|
\end{aligned}
\end{align*}
mamy
\begin{align}\label{H_as_g}
\begin{aligned} 
H_{n,k}(x)=E_x\left(\psi_{n,k}(x_0, x_1, \ldots, x_{n+k})\right).
\end{aligned}
\end{align}

Przypomnijmy, że $(x_n)_{n\in N_0}$ oraz $(y_n)_{n\in N_0}$ są łańcuchami Markowa o~wartościach w~przestrzeni stanów $X$, funkcji przejścia $\Pi_{\varepsilon}$ i rozkładach początkowych $\mu,\nu\in M^{2+\delta}_1(X)$. W~szczególności rolę rozkładów początkowych mogą pełnić $\delta_x$ oraz $\delta_y$. Ze względu na konieczność zastosowania modelu pomocniczego ustalmy $(h_i)_{i\in N}\subset\bar{B}(0,\varepsilon)$ i rozważmy niejednorodne łańcuchy Markowa $(\tilde{x}_n)_{n\in N_0}$, $(\tilde{y}_n)_{n\in N_0}$ o ciągu funkcji przejścia $\left(\Pi^1_{h_i}\right)_{i\in N}$ (zob. (\ref{def:Pi_hn})) oraz rozkładach początkowych $\delta_x$ i $\delta_y$. Zauważmy, że wobec (\ref{Pi_epsilon-as-Pi_h}) mamy 
\begin{align*}
\begin{aligned}
\Pi_{\varepsilon}(x,\cdot)=\int_{\bar{B}(0,\varepsilon)}\Pi^1_{h}(x,\cdot)\nu^{\varepsilon}(dh)
\end{aligned}
\end{align*}
i stąd
\begin{align}
\begin{aligned}
E_x\left(\psi(x_0,\ldots,x_{n+k})\right)=\int_{\left(\bar{B}(0,\varepsilon)\right)^{n+k}}E_x\left(\psi(\tilde{x}_0,\ldots,\tilde{x}_{n+k})\right)\;\nu^{\varepsilon}(dh_1)\ldots\nu^{\varepsilon}(dh_{n+k}).
\end{aligned}
\end{align}
Wprowadźmy funkcję pomocniczą
\begin{align}\label{barH_as_psi}
\begin{aligned}
\tilde{H}_{n,k}(x)=E_x\left(\psi_{n,k}(\tilde{x}_0, \tilde{x}_1, \ldots, \tilde{x}_{n+k})\right).
\end{aligned}
\end{align}
Przypomnijmy, że $E_{x,y}$ jest wartością oczekiwaną względem miary sprzęgającej $C^{\infty}_{h_1,h_2,\ldots}((x,y),\cdot)$ lub $\hat{C}^{\infty}_{h_1,h_2,\ldots}((x,y,1),\cdot)$. 
Zauważmy, że możemy pisać 
\begin{align}\label{barH-barH}
\begin{aligned}
\lim_{n\to\infty}\lim_{k\to\infty}
\left|\tilde{H}_{n,k}(x)-\tilde{H}_{n,k}(y)\right|
&=\lim_{n\to\infty}\lim_{k\to\infty}
\left|E_x\left(\psi(\tilde{x}_0,\ldots,\tilde{x}_{n+k})\right)-E_y\left(\psi(\tilde{y}_0,\ldots,\tilde{y}_{n+k})\right)\right|\leq\\
&\leq 
\lim_{n\to\infty}\lim_{k\to\infty}E_{x,y}
\left|\psi(\tilde{x}_0,\ldots,\tilde{x}_{n+k})-\psi(\tilde{y}_0,\ldots,\tilde{y}_{n+k})\right|.
\end{aligned}
\end{align}
Łatwo pokazać, że
\begin{align*}
\begin{aligned}
 \min_{j\in\{1,\ldots,n\}}\left(f_j(x_j)\right)
-\min_{j\in\{1,\ldots,n\}}\left(f_j(y_j)\right)\leq
\max_{i\in\left\{1,\ldots,n\right\}}\left|f_i(x_i)-f_i(y_i)\right|
\end{aligned}
\end{align*}
dla $f_i:X\to R$ oraz $x_i, y_i\in X$, $i\in\left\{1,\ldots,n\right\}$. Zatem mamy
\begin{align}\label{estim:min_na_max}
\begin{aligned}
&\lim_{n\to\infty}\lim_{k\to\infty}\left|\tilde{H}_{n,k}(x)-\tilde{H}_{n,k}(y)\right|\leq\\
&\leq \lim_{n\to\infty}\lim_{k\to\infty}E_{x,y}
\Bigg(\max_{i\in\{n,\ldots,n+k\}}\frac{1}{i}\sum_{l=1}^{i}\Big|\left(\chi(\tilde{x}_l)-\chi(\tilde{x}_{l-1})+g(\tilde{x}_{l-1})\right)^2\wedge i(1+\sigma^2)-\\
&\quad
-\left(\chi(\tilde{y}_l)-\chi(\tilde{y}_{l-1})+g(\tilde{y}_{l-1})\right)^2\wedge i(1+\sigma^2)\Big|\Bigg).
\end{aligned}
\end{align}
Zauważmy, że prawa strona (\ref{estim:min_na_max}) jest równa wyrażeniu
\begin{align*}
\begin{aligned}
&\lim_{n\to\infty}\lim_{k\to\infty}E_{x,y}
\Bigg(\max_{i\in\{n,\ldots,n+k\}}\frac{1}{i}\sum_{l=k_0}^{i}\Big|\left(\chi(\tilde{x}_l)-\chi(\tilde{x}_{l-1})+g(\tilde{x}_{l-1})\right)^2\wedge i(1+\sigma^2)-\\
&\quad
-\left(\chi(\tilde{y}_l)-\chi(\tilde{y}_{l-1})+g(\tilde{y}_{l-1})\right)^2\wedge i(1+\sigma^2)\Big|\Bigg)
\end{aligned}
\end{align*}
dla dowolnego $k_0\in N$. 
Ponieważ rozważane funkcje są ograniczone, otrzymujemy
\begin{align}\label{estim:barH-barH_further}
\begin{aligned}
&\left|\tilde{H}_{n,k}(x)-\tilde{H}_{n,k}(y)\right|\leq\\
&\leq 
E_{x,y}\left(
\max_{i\in\{n,\ldots,n+k\}}\frac{1}{i}\sum_{l=k_0}^{i}\left|\left(\chi(\tilde{x}_l)-\chi(\tilde{x}_{l-1})+g(\tilde{x}_{l-1})\right)-\left(\chi(\tilde{y}_l)-\chi(\tilde{y}_{l-1})+g(\tilde{y}_{l-1})\right)\right|2i(1+\sigma^2)\right)\leq\\
&\leq 2(1+\sigma^2)E_{x,y}
\left(\sum_{l=k_0}^{n+k}|\chi(\tilde{x}_l)-\chi(\tilde{y}_l)|+|\chi(\tilde{x}_{l-1})-\chi(\tilde{y}_{l-1})|+|g(\tilde{x}_{l-1})-g(\tilde{y}_{l-1})|\right).
\end{aligned}
\end{align}
Dalej mamy
\begin{align*}
\begin{aligned}
&E_{x,y}\left|\chi(\tilde{x}_l)-\chi(\tilde{y}_l)\right|\leq\\
&\leq \sum_{i=0}^{\infty}E_{x,y}
\left|U_{\varepsilon}^ig(\tilde{x}_l)-U_{\varepsilon}^ig(\tilde{y}_l)\right|=\\
&=\sum_{i=0}^{\infty}\int_{X^2}\left|U^i_{\varepsilon}g(u)-U^i_{\varepsilon}g(v)\right|\left(\Pi^*_{X^2}\Pi^*_l\hat{C}^{\infty}_{h_1,h_2,\ldots}((x,y,1),\cdot)\right)(du\times dv)=\\
&=\sum_{i=0}^{\infty}\int_{X^2}
\left|\left\langle g,P^i_{\varepsilon}\delta_u\right\rangle
-\left\langle g,P^i_{\varepsilon}\delta_v\right\rangle\right|
\left(\Pi^*_{X^2}\Pi^*_l\hat{C}^{\infty}_{h_1,h_2,\ldots}((x,y,1),\cdot)\right)(du\times dv)=\\
&=\sum_{i=0}^{\infty}\int_{X^2}
\left|\int_{\left(\bar{B}(0,\varepsilon)\right)^i}\int_{X^2}
g(z)\left(\Pi^i_{h_{l+1},\ldots,{h}_{l+i}}(u,\cdot)-\Pi^i_{{h}_{l+1},\ldots,{h}_{l+i}}(v,\cdot)\right)(dz)
\:\nu^{\varepsilon}(d{h}_{l+1})\ldots\nu^{\varepsilon}(d{h}_{l+i})\right|\times\\
&\quad\times
\left(\Pi^*_{X^2}\Pi^*_l\hat{C}^{\infty}_{h_1,h_2,\ldots}((x,y,1),\cdot)\right)(du\times dv)
\leq\\
&\leq\sum_{i=0}^{\infty}\int_{\left(\bar{B}(0,\varepsilon)\right)^i}\int_{X^2}\int_{X^2}
\left|g(z_1)-g(z_2)\right|
\left(\Pi^*_{X^2}\Pi^*_i\hat{C}^{\infty}_{{h}_{l+1},{h}_{l+2},\ldots}((u,v,1),\cdot)\right)(dz_1\times dz_2)\times\\
&\quad\times\left(\Pi^*_{X^2}\Pi^*_l\hat{C}^{\infty}_{h_1,h_2,\ldots}((x,y,1),\cdot)\right)(du\times dv)
\:\nu^{\varepsilon}(d{h}_{l+1})\ldots\nu^{\varepsilon}(d{h}_{l+i}).
\end{aligned}
\end{align*}
Przypomnijmy, że miary $C^n_{h_1,\ldots,h_n}((x,y),\cdot)$ są $n$-tymi współrzędnymi miary $C^{\infty}_{h_1,h_2,\ldots}((x,y),\cdot)$ na całych trajektoriach dla $x,y\in X$, $n\in N$. Zgodnie z regułą konstrukcji miar jednowymiarowych na kolejnych współrzędnych (por. (\ref{def:Pi_hn}), (\ref{Q1Qnb})) mamy
\begin{align*}
\begin{aligned}
C^n_{h_1,\ldots,h_n}((x,y),\cdot)
=\int_XC^1_{h_n}((z_1,z_2),\cdot)C^{n-1}_{h_1,\ldots,h_{n-1}}((x,y),dz_1\times dz_2),
\end{aligned}
\end{align*}
a zatem
\begin{align*}
\begin{aligned}
&E_{x,y}\left|\chi(\tilde{x}_l)-\chi(\tilde{y}_l)\right|
\leq\\
&\leq
\sum_{i=0}^{\infty}\int_{\left(\bar{B}(0,\varepsilon)\right)^i}\int_{X^2}
\left|g(z_1)-g(z_2)\right|
\left(\Pi^*_{X^2}\Pi^*_{l+i}\hat{C}^{\infty}_{h_1,h_2,\ldots}((x,y,1),\cdot)\right)(dz_1\times dz_2)
\:\nu^{\varepsilon}(d{h}_{l+1})\ldots\nu^{\varepsilon}(d{h}_{l+i}).
\end{aligned}
\end{align*}
Na mocy wniosku \ref{wniosek_g} otrzymujemy więc
\begin{align}\label{estim:chi(x_l)-chi(y_l)_bar}
\begin{aligned}
E_{x,y}\left|\chi(\tilde{x}_l)-\chi(\tilde{y}_l)\right|
&\leq \sum_{i=0}^{\infty}
\int_{\left(\bar{B}(0,\varepsilon)\right)^i}
C_5Gq^{l+i}(1+V(x)+V(y))
\:\nu^{\varepsilon}(d{h}_{l+1})\ldots\nu^{\varepsilon}(d{h}_{l+i})\leq\\
&\leq C_5G(1+V(x)+V(y))\sum_{i=l}^{\infty}q^i.
\end{aligned}
\end{align}
Ponadto mamy 
\begin{align}\label{estim:g-g}
\begin{aligned}
E_{x,y}\left|g(\tilde{x}_{l-1})-g(\tilde{y}_{l-1})\right|
\leq C_5Gq^{l-1}(1+V(x)+V(y)).
\end{aligned}
\end{align}
Dzięki (\ref{estim:chi(x_l)-chi(y_l)_bar}) i (\ref{estim:g-g}) możemy szacować (\ref{estim:barH-barH_further}) z góry wyrażeniem
\begin{align*}
\begin{aligned}
&2\left(1+\sigma^2\right)C_5G(1+V(x)+V(y))\sum_{l=k_0}^{n+k}
\left(\sum_{i=l}^{\infty}q^i+\sum_{i=l-1}^{\infty}q^i+q^{l-1}\right)\leq\\
&\quad\leq 4\left(1+\sigma^2\right)C_5G(1+V(x)+V(y))\sum_{l=k_0}^{n+k}\left(q^{l-1}\sum_{i=0}^{\infty}q^i\right)=\\
&\quad=\frac{4}{1-q}\left(1+\sigma^2\right)C_5G(1+V(x)+V(y))\sum_{l=k_0}^{n+k}q^{l-1}.
\end{aligned}
\end{align*}
Oszacowanie jest niezależne od $(h_i)_{i\in N}$, zatem
\begin{align*}
\begin{aligned}
\lim_{n\to\infty}\lim_{k\to\infty}\left|H_{n,k}(x)-H_{n,k}(y)\right|
\leq \frac{4}{1-q}\left(1+\sigma^2\right)C_5G(1+V(x)+V(y))\sum_{l=k_0}^{\infty}q^{l-1}.
\end{aligned}
\end{align*}
Ponieważ $k_0$ możemy wybrać w sposób dowolny, to w szczególności możemy ustalić $k_0\in N$ takie, że wartość szeregu $\sum_{l=k_0}^{\infty}q^{l-1}$ jest bliska zeru, co implikuje 
\[\lim_{n\to\infty}\lim_{k\to\infty}\left|H_{n,k}(x)-H_{n,k}(y)\right|=0\quad\text{dla wszystkich $x,y\in X$.}\]
Dowód został ukończony.
\end{proof}

\begin{lemat}\label{lem2}
Niech $\sigma^2>0$ oraz $\beta>0$. Przy założeniach (I)-(VI) oraz (II$^{\prime\prime}$)  różnice martyngałowe $(Z_n)_{n\in N}$ spełniają warunki (\ref{th1}) oraz (\ref{th2}).
\end{lemat}

\begin{proof}

Niech $\mu\in M_1^{2+\delta}(X)$. Zauważmy, że 
\begin{align*}
\begin{aligned}
\sum_{n=1}^{\infty}s_n^{-4}E_{\mu}\left(Z_n^41_{\left\{|Z_n|<\beta s_n\right\}}\right)
&\leq\sum_{n=1}^{\infty}s_n^{-4}\beta^{2-\delta}s_n^{2-\delta}E_{\mu}|Z_n|^{2+\delta}
\leq \beta^{2-\delta}\sup_{n\in N}E_{\mu}|Z_n|^{2+\delta}\sum_{n=1}^{\infty}s_n^{-2-\delta}.
\end{aligned}
\end{align*}
Zgodnie z (\ref{supEZn-delta-finite}) mamy $\sup_{n\geq 1}E_{\mu}|Z_n|^{2+\delta}<\infty$, a z lematu \ref{lem:sn/n=sigma} otrzymujemy $\sum_{n=1}^{\infty}s_n^{-2-\delta}<\infty$, co~dowodzi warunku (\ref{th1}).

Aby wykazać warunek (\ref{th2}), wystarczy zauważyć, że 
\begin{align*}
\begin{aligned}
\sum_{n=1}^{\infty}s_n^{-1}E_{\mu}\left(|Z_n|1_{\left\{|Z_n|\geq \vartheta s_n\right\}}\right)&\leq\sum_{n=1}^{\infty}s_n^{-1}E_{\mu}\left(\frac{|Z_n|^{2+\delta}}{(\vartheta s_n)^{1+\delta}}\right)\leq \vartheta^{-1-\delta}\sup_{n\in N}E_{\mu}|Z_n|^{2+\delta}\sum_{n=1}^{\infty}s_n^{-2-\delta}<\infty.
\end{aligned}
\end{align*}
Dowód został ukończony.
\end{proof}

\section{Prawo iterowanego logarytmu dla uogólnionego modelu cyklu komórkowego}

Niech $\mathcal{C}$ będzie przestrzenią funkcji ciągłych $f:[0,1]\to R$ z~normą supremum. Ponadto niech $\mathcal{K}$ będzie podprzestrzenią przestrzeni $\mathcal{C}$ składającą się z funkcji $f$ absolutnie ciągłych i~takich, że $f(0)=0$ oraz $\int_0^1\left(f'(t)\right)^2dt\leq 1$.

Ustalmy funkcję $g\in B(X)$ spełniającą warunek Lipschitza ze stałą $L_g>0$. Ponadto niech $\left\langle g,\mu_*\right\rangle=0$, gdzie $\mu_*$ jest miarą niezmienniczą względem operatora $P_{\varepsilon}$.

Dana jest przestrzeń probabilistyczna $(\Omega, \mathcal{F},\mathbb{P})$. Ciągi zmiennych losowych $(M_n)_{n\in N_0}$ i $(Z_n)_{n\in N_0}$ są odpowiednio postaci (\ref{def:M_n}) i (\ref{def:Z_n}), a~wariancja $\sigma^2$ jest dana wzorem (\ref{def:sigma2}). Załóżmy, że $\sigma^2>0$ i zdefiniujmy ciąg
\begin{align}\label{def:Xi}
\Xi_n(t)=\frac{\sum_{i=1}^kg(x_i)+(nt-k)g(x_{k+1})}{\sigma\sqrt{2n\log\log n}}
\end{align}
dla $k\leq nt\leq k+1$, $\:k\in\{1,\ldots, n-1\}$, $\:t>0$, $\:n>e$, oraz $\Xi_n(t)=0$ w przeciwnym przypadku.

\begin{twierdzenie}\label{PIL}
Niech $(X,\varrho)$ będzie przestrzenią polską oraz $(x_n)_{n\in N_0}$ -- łańcuchem Markowa o~funkcji przejścia $\Pi_{\varepsilon}$ i rozkładzie początkowym $\mu\in M_1^{2+\delta}(X)$. Załóżmy, że warunki (I)-(VI) są spełnione oraz warunek (II) jest wzmocniony do (II$^{\prime\prime}$). Wówczas ciąg $(\Xi_n)_{n>e}$ dany wzorem (\ref{def:Xi}) jest $\mathbb{P}_{\mu}$-p.n. relatywnie zwarty w~zbiorze $\mathcal{C}$, a~zbiór jego punktów granicznych pokrywa się ze~zbiorem $\mathcal{K}$.
\end{twierdzenie}

\begin{proof} 
Zdefiniujmy ciąg
\begin{align*}
\begin{aligned}
\eta_n(t)=\frac{M_k+\left(s_n^2t-s_k^2\right)\left(s_{k+1}^2-s_k^2\right)^{-1}Z_{k+1}}{\sigma\sqrt{2n\log\log n}}
\end{aligned}
\end{align*}
dla $s_k^2\leq s_n^2t\leq s_{k+1}^2$, $k\in\{1,\ldots, n-1\}$, $t>0$, $n>e$, oraz $\eta_n(t)=0$ w przeciwnym przypadku. 
Ponieważ $\lim_{n\to\infty}s_n^2/n=\sigma^2$ (lemat \ref{lem:sn/n=sigma}), to mamy
\begin{align*}
\begin{aligned}
\lim_{n\to\infty}\frac{\sqrt{2s^2_n\log\log s_n^2}}{\sigma\sqrt{2n\log\log n}}=1.
\end{aligned}
\end{align*}
Na mocy lematów \ref{lem1} i \ref{lem2} spełnione są założenia twierdzenia \ref{tw:heyde-scott}, a więc ciąg $\left(\eta_n(t)\right)_{n> e}$ jest relatywnie zwarty w przestrzeni $\mathcal{C}$, a zbiór jego punktów granicznych pokrywa się ze zbiorem~$\mathcal{K}$. 

Niech $t\in(0,1]$ oraz $n>e$. Jeśli $k\leq nt\leq k+1$, to
\begin{align*}
\begin{aligned}
\frac{k\sigma^2}{s_k^2}s_k^2\leq \frac{n\sigma^2}{s_n^2}ts_n^2\leq\frac{(k+1)\sigma^2}{s_{k+1}^2}s_{k+1}^2.
\end{aligned}
\end{align*}
Ustalmy
\begin{align}
\begin{aligned}
\hat{\eta}_n(t)=\frac{M_k+(nt-k)Z_{k+1}}{\sigma\sqrt{2n\log\log n}}\quad\text{dla }k\in N,\:k\leq nt\leq k+1.
\end{aligned}
\end{align}
Ponieważ $\lim_{n\to\infty}n\sigma^2/s_n^2= 1$, to dla odpowiednio dużego $n\in N$ oraz $\tilde{\epsilon}>0$ mamy
\begin{align*}
\begin{aligned}
(1-\tilde{\epsilon})s_k^2\leq (1+\tilde{\epsilon})s_n^2 t\leq (1+\tilde{\epsilon})^2(1-\tilde{\epsilon})^{-1}s_{k+1}^2.
\end{aligned}
\end{align*}
Wobec tego istnieje $t_*\in\left[t(1-\tilde{\epsilon})(1+\tilde{\epsilon})^{-1},t(1+\tilde{\epsilon})(1-\tilde{\epsilon})^{-1}\right]$ o własności $s_k^2\leq s_n^2t_*\leq s_{k+1}^2$. Z~drugiej strony długości przedziałów $\left[s_k^2/s_n^2,s_{k+1}^2/s_n^2\right]$ dla $k\in \left\{1,\ldots,n-1\right\}$ zbiegają do zera, gdy $n\to\infty$, stąd istnieje $t_n>0$ taki, że $\hat{\eta}_n(t)=\eta_n(t_n)$ oraz $\lim_{n\to\infty}t_n=t$. Ponieważ ciąg $(\eta_n(t))_{n> e}$ jest relatywnie zwarty w  $\mathcal{C}$ oraz zbiór jego punktów granicznych pokrywa się z $\mathcal{K}$, to $(\hat{\eta}_n(t))_{n>e}$ również spełnia te własności.

Ustalmy $\tilde{\epsilon}>0$ oraz zdefiniujmy zbiór
\begin{align}
\begin{aligned}
\mathcal{A}_n
=\left\{\omega\in\Omega:\;\frac{\left|M_n(\omega)-\sum_{i=1}^{n-1}g(x_i(\omega))\right|}{\sqrt{n}}
\geq \frac{\tilde{\epsilon}}{2}\right\}
\cup\left\{\omega\in\Omega:\;\frac{\left|Z_n(\omega)-g(x_n(\omega))\right|}{\sqrt{n}}
\geq\frac{\tilde{\epsilon}}{2}\right\}\;\text{dla }n\in N.
\end{aligned}
\end{align}
Zauważmy, że dla $\epsilon>0$ oraz $\xi,\zeta\geq 0$ możemy wybrać $\delta_1>0$ taką, że $(\xi+\zeta)^{2+\delta}\leq (2+\epsilon)\left(\xi^{2+\delta} + \zeta^{2+\delta}\right)$ dla $\delta<\delta_1$. Korzystając z tej własności oraz z nierówności Markowa, otrzymujemy
\begin{align*}
\begin{aligned}
\mathbb{P}_{\mu}
\left(\omega\in\Omega:\;\frac{\left|M_n(\omega)-\sum_{i=1}^{n-1}g(x_i(\omega))\right|}{\sqrt{n}}
\geq\frac{\tilde{\epsilon}}{2}\right)
&=\mathbb{P}_{\mu}\left(\omega\in\Omega:\;\frac{\left|\chi(x_n(\omega))-\chi(x_0(\omega))\right|}{\sqrt{n}}
\geq\frac{\tilde{\epsilon}}{2}\right)\leq\\
&\leq \left(\frac{2}{\tilde{\epsilon}}\right)^{2+\delta}(2+\epsilon)
\frac{E_{\mu}|\chi(x_n)|^{2+\delta}+E_{\mu}|\chi(x_0)|^{2+\delta}}{n^{1+\delta/2}}.
\end{aligned}
\end{align*}
Na mocy lematu \ref{lem:chi-chi} mamy
\begin{align*}
\begin{aligned}
E_{\mu}|\chi(x_n)|^{2+\delta}&=\int_X|\chi(x)|^{2+\delta}P^n_{\varepsilon}\mu(dx)\leq\\
&\leq (2+\epsilon)|\chi(\bar{x})|^{2+\delta}+(2+\epsilon)\int_X|\chi(x)-\chi(\bar{x})|^{2+\delta}P^n_{\varepsilon}\mu(dx)\leq\\
&\leq (2+\epsilon)\left|\chi(\bar{x})\right|^{2+\delta}+(2+\epsilon)^2G^{2+\delta}C_5^{2+\delta}
(1-q)^{-(2+\delta)}\left(1+\left\langle V^{2+\delta},P^n_{\varepsilon}\mu\right\rangle \right).
\end{aligned}
\end{align*}
Zauważmy, że z własności (\ref{fact:finite_chi}) oraz wniosku \ref{wn:skończone_momenty_P_varepsilon_mu} oba składniki powyższego oszacowania są skończone. Otrzymujemy więc
\begin{align}\label{estim:An1}
\begin{aligned}
\mathbb{P}_{\mu}\left(\omega\in\Omega:\;
\frac{\left|M_n(\omega)-\sum_{i=1}^{n-1}g(x_i(\omega))\right|}{\sqrt{n}}\geq\frac{\tilde{\epsilon}}{2}\right)
\leq \frac{c_1}{n^{1+\delta/2}},
\end{aligned}
\end{align}
gdzie $c_1$ jest pewną stałą niezależną od $n$. Podobnie, 
\begin{align}\label{estim:An2}
\begin{aligned}
\mathbb{P}_{\mu}\left\{\omega\in\Omega:\;\frac{|Z_n(\omega)-g(x_n(\omega))|}{\sqrt{n}}\geq\frac{\tilde{\epsilon}}{2}\right\}
=\mathbb{P}_{\mu}\left\{\omega\in\Omega:\;\frac{|\chi(x_{n+1}(\omega))-\chi(x_{n}(\omega))|}{\sqrt{n}}\geq\frac{\tilde{\epsilon}}{2}\right\}
\leq \frac{c_2}{n^{1+\delta/2}},
\end{aligned}
\end{align}
gdzie stała $c_2$ również nie zależy od $n$. Oszacowania (\ref{estim:An1}) oraz (\ref{estim:An2}) implikują zbieżność szeregu $\sum_{n=1}^{\infty}\mathbb{P}_{\mu}(\mathcal{A}_n)$.

Ostatecznie z lematu Borela-Cantelli'ego poniższe oszacowanie zachodzi $\mathbb{P}_{\mu}$-p.n.:
\begin{align*}
\begin{aligned}
{\lim\sup}_{n\to\infty}{\sup}_{0< t\leq 1}\left|\frac{M_k-(nt-k)Z_{k+1}}{\sigma\sqrt{2n\log\log n}}-\frac{\sum_{i=1}^kg(x_i)+(nt-k)g(x_{k+1})}{\sigma\sqrt{2n\log\log n}}\right|<\tilde{\epsilon},
\end{aligned}
\end{align*}
gdzie $k\leq nt\leq k+1$, wobec czego mamy 
\begin{align}
{\lim\sup}_{n\to\infty}\sup_{0<t\leq 1} \left|\hat{\eta}_n(t)-\Xi_n(t)\right|\leq \tilde{\epsilon}.
\end{align}
Ponieważ $\tilde{\epsilon}>0$ wybraliśmy w sposób dowolny, dowód jest ukończony.
\end{proof}

\end{document}